\documentclass[a4paper,11pt]{report}
\usepackage{amssymb}
\usepackage{amsmath}
\usepackage{amscd}
\usepackage[francais]{babel}
\usepackage[T1]{fontenc}
\usepackage[all,cmtip]{xy}
\usepackage{fancyhdr}
\usepackage{hyperref}

\def\RM{\mathbf{R}}
\def\CM{\mathbf{C}}
\def\ZM{\mathbf{Z}}
\def\NM{\mathbf{N}}
\def\QM{\mathbf{Q}}
\def\FM{\mathbf{F}}
\def\gr{\mathrm{gr}}
\def\SM{\mathcal{S}}
\def\AM{\mathcal{A}}
\def\MM{\mathbf{M}}
\def\Tr{\mathrm{Tr}}
\def\dim{\mathrm{dim}\,}
\def\codim{\mathrm{codim}}
\def\End{\mathrm{End}}
\def\car{\mathrm{car}}
\def\Hom{\mathrm{Hom}}
\def\Aut{\mathrm{Aut}}
\def\Int{\mathrm{Int}}
\def\Id{\mathrm{Id}}
\def\Ver{\mathrm{Ver}}
\def\Im{\mathrm{Im\, }}

\def\App{{\it Application. }}
\def\exo{{\it Exercice. }}
\def\ex{{\it Exemple. }}
\def\exs{{\it Exemples.\\ }}
\def\rmq{{\it Remarque. }}
\def\rmqs{{\it Remarques.\\ }}
\def\dem{{\it D\'{e}monstration. }}

\newtheorem{defi}{D\'{e}finition}[chapter]
\newtheorem{thm}{Th\'{e}or\`{e}me}[chapter]
\newtheorem{coro}[thm]{Corollaire}
\newtheorem{prop}[thm]{Proposition}

\newtheorem{ppts}[thm]{Propri\'{e}t\'{e}s}
\newtheorem{lemme}[thm]{Lemme}

\def \findem {\hfill{$\square$}}

\begin{document}

\thispagestyle{empty}
\begin{center}
\vglue 15pc

{\Huge \sc Groupes finis} \vglue 4pc

{\LARGE Jean-Pierre {\sc Serre}}

\vglue 3pc
\end{center}

\newpage
\thispagestyle{empty}
\begin{center}
\vglue 15pc

{\large Cours \`{a} l'\'Ecole Normale Sup\'{e}rieure de Jeunes
Filles, $\mit 1978/1979$}

\

\

r\'{e}dig\'{e} par Martine  {\sc Buhler} et Catherine {\sc Goldstein}\\
(Montrouge, $\mit 1979$)

\

\

r\'{e}vis\'{e} et transcrit en \LaTeX\ par\\
Nicolas {\sc Billerey}, Olivier {\sc Dodane} et Emmanuel {\sc
Rey}\\
(Strasbourg -- Paris, $\mit 2004$)

\vglue 3pc
\end{center}

\newpage
\tableofcontents

\

{\bf Bibliographie} \hfill{\bf \pageref{biblio}}

\

{\bf Index} \hfill{\bf \pageref{index}}

\chapter{Pr\'{e}liminaires}
Ce chapitre est essentiellement constitu\'{e} de rappels sur la
th\'{e}orie g\'{e}n\'{e}rale des groupes. La lettre $G$ d\'{e}signe
un groupe.

\section{Actions de groupes}\label{action}

\begin{defi}
On dit que le groupe $G$ \emph{op\`{e}re \`{a} gauche} sur un
ensemble $X$ si l'on s'est donn\'{e} une application
$$
\left\{\!\!
 \begin{array}{rcl}
  G \times X & \longrightarrow & X\\
  (g,x) & \longmapsto & g.x
 \end{array}
\right.
$$
v\'{e}rifiant les conditions:
\begin{enumerate}
\item[(1)] $g.(g'.x)=(gg').x$ pour tout $x\in X$ et tout couple
$(g,g')\in G\times G$.

\item[(2)] $1.x=x$ pour tout $x\in X$, o\`{u} $1$ est l'\'{e}l\'{e}ment neutre
de $G$.
\end{enumerate}

\end{defi}

\rmq La donn\'{e}e d'une action \`{a} gauche de $G$ sur $X$
\'{e}quivaut \`{a} la donn\'{e}e d'un homomorphisme $\tau$ de $G$
dans le groupe $\SM_X$ des permutations de $X$ d\'{e}fini pour tout
$g\in G$ et tout $x\in X$ par $\tau(g)(x)=g.x.$

\bigskip On aurait une d\'{e}finition analogue pour les op\'{e}rations \`{a}
droite.

\bigskip
Le groupe $G$ d\'{e}coupe alors $X$ en \emph{orbites}: deux
\'{e}l\'{e}ments $x$ et $y$ de $X$ sont dans la m\^{e}me orbite si
et seulement s'il existe $g\in G$ tel que $x=g.y$. L'ensemble des
orbites est le quotient de $X$ par $G$ et est not\'{e} $G\backslash
X$ dans le cas d'une action \`{a} gauche (et $X/G$ dans le cas d'une
action \`{a} droite).

\begin{defi}
On dit que $G$ agit transitivement sur X si $G\backslash X$ est
r\'{e}duit \`{a} un \'{e}l\'{e}ment.
\end{defi}

En particulier, le groupe $G$ agit transitivement sur chaque
orbite.

\begin{defi}
Soit $x\in X$; on appelle \emph{stabilisateur de $x$} (ou
\emph{fixateur} de $x$) et on note $H_x$ le sous-groupe de $G$
form\'{e} des \'{e}l\'{e}ments $g\in G$ qui fixent $x$ (i.e. tels
que $g.x=x$).
\end{defi}\label{stabilisateur}

\rmq Si $G$ op\`{e}re transitivement sur $X$ et si $x\in X$, on a
une bijection de $G/H_x$ sur $X$ donn\'{e}e par $g H_x \longmapsto
g.x$, o\`{u} $G/H_x$ est l'ensemble des classes \`{a} gauche de $G$
modulo $H_x$. Si $x'\in X$, il existe $g\in G$ tel que $x'=g.x$.
Alors $H_{x'}=g{H_x}g^{-1}$. Donc changer de point de base revient
\`{a} remplacer le stabilisateur de $x$ par un de ses conjugu\'{e}s.
Inversement, si $H$ est un sous-groupe de $G$, alors $G$ agit
transitivement sur $G/H$ et $H$ stabilise la classe de $1$. Ainsi la
donn\'{e}e de $X$ sur lequel $G$ op\`{e}re transitivement revient
\`{a} celle d'un sous-groupe de $G$, d\'{e}termin\'{e} \`{a}
conjugaison pr\`{e}s.

\bigskip\ex Soit $X$ une droite affine d\'{e}finie sur un corps $K$ et soit
$G$ le groupe des similitudes
$$G=\left\{x\mapsto ax+b,\, a\in K^*,\,b\in K\right\}.$$

Le groupe $G$ op\`{e}re transitivement sur $X$. Si $x\in X$, le
stabilisateur de $x$ est le groupe des homoth\'{e}ties centr\'{e}es
en $x$.

\bigskip\App Soit $G$ un groupe \emph{fini}, dont on note $|G|$ l'ordre.
Soit $X$ un ensemble o\`{u} $G$ op\`{e}re. On a $X=\coprod_{i\in
I}{Gx_i}$ o\`{u} les $Gx_i$ sont les orbites (2 \`{a} 2 disjointes)
sous l'action de $G$, les $x_i$ formant un syst\`{e}me de
repr\'{e}sentants des \'{e}l\'{e}ments de $G\backslash X$. On a vu
que $Gx_i$ est en bijection avec $G/{H_{x_i}}$, donc
$|Gx_i|=|G|.{|H_{x_i}|}^{-1}$. On en d\'{e}duit $|X|=\sum_{i\in
I}{|G|.{|H_{x_i}|}^{-1}}$ puis $|X|.{|G|}^{-1}=\sum_{i\in
I}{{|H_{x_i}|}^{-1}}$.

\bigskip\noindent{\it Cas particulier.} Le groupe $G$ op\`{e}re sur
lui-m\^{e}me par automorphismes int\'{e}rieurs; on a une
application:
$$
\left\{\!\!
 \begin{array}{rcl}
  G & \longrightarrow & \SM_G\\
  x & \longmapsto & {\rm int}_x
 \end{array}
\right.$$ o\`{u} ${\rm int}_x(y)=xyx^{-1}={}^xy$. Les orbites sont
les classes de conjugaison\label{classeconjugaison}. Le
stabilisateur d'un \'{e}l\'{e}ment $x$ de $G$ est l'ensemble des
\'{e}l\'{e}ments de $G$ qui commutent \`{a} $x$ (on l'appelle
\emph{centralisateur}\label{centralisateur} de $x$ et on le note
$C_G(x)$). On a $1=\sum_{i\in I}{{|C_{G}(x_i)|}^{-1}}$ o\`{u}
${(x_i)}_{i\in I}$ est un syst\`{e}me de repr\'{e}sentants des
classes de conjugaison. Pour $x_i=1$, on a $C_{G}(x_i)=G$ et donc
$\sup_{i\in I}{|C_{G}(x_i)|}=|G|$.

\bigskip\exo\\
$(i)$ Si $h$ est un entier $\geqslant 1$, montrer qu'il n'y a qu'un
nombre fini de d\'{e}compositions $1=\sum_{i=1}^{h}{\frac{1}{n_i}}$
avec $n_i\in \ZM$, $n_i\geqslant 1$. [Par exemple, si $h=3$, les
seuls $n_i$ possibles sont $(3,3,3)$, $(2,4,4)$ et $(2,3,6)$.]

$(ii)$ En d\'{e}duire que, si un groupe fini $G$ a un nombre de
classes de conjugaison \'{e}gal \`{a} $h$, l'ordre de $G$ est
major\'{e} par une constante $N(h)$ ne d\'{e}pendant que de $h$. (On
peut prendre $N(h)$ de la forme ${c_1}^{{c_2}^h}$, o\`{u} $c_1$,
$c_2$ sont des constantes $>0$. J'ignore si l'on peut faire beaucoup
mieux.)

\section{Sous-groupes normaux; sous-groupes caract\'{e}ristiques; groupes simples}\label{grpe simple}

\begin{defi}
On dit qu'un sous-groupe $H$ de $G$ est \emph{normal} (ou
\emph{invariant}) si pour tout $x\in G$ et tout $h\in H$, on a
$xhx^{-1}\in H$.
\end{defi}\label{normal}

Cela revient \`{a} dire que le sous-groupe $H$ est stable par tout
automorphisme int\'{e}rieur. Une telle situation se d\'{e}crit par
une suite exacte:
$$\xymatrix{\{1\} \ar[r] & H \ar[r] & G \ar[r] & G/H \ar[r] & \{1\}}.$$

\rmq Si $H$ est un sous-groupe de $G$, il existe un plus grand
sous-groupe de $G$ dans lequel $H$ est normal, \`{a} savoir
l'ensemble des $g\in G$ tels que $gHg^{-1}=H$. On l'appelle le
\emph{normalisateur de $H$ dans $G$}\label{normalisateur}, et on le
note $N_G(H)$. On dit qu'une partie de $G$ \emph{normalise} $H$ si
elle est contenue dans $N_G(H)$.

\begin{defi}
On dit qu'un sous-groupe $H$ de $G$ est \emph{caract\'{e}ristique}
s'il est stable par tout automorphisme de $G$.
\end{defi}
\label{caracteristique}

Un tel sous-groupe est normal.

\bigskip\ex Le \emph{centre}\label{centre} de $G$ (ensemble des
\'{e}l\'{e}ments qui commutent \`{a} tous les \'{e}l\'{e}ments de
$G$) est un sous-groupe caract\'{e}ristique. Il en est de m\^{e}me
du \emph{groupe d\'{e}riv\'{e}} de $G$, ainsi des sous-groupes
$D^nG$, $C^iG$ et $\Phi(G)$ d\'{e}finis au chap. \ref{chap3}.

\begin{defi}
On dit qu'un groupe $G$ est \emph{simple} lorsqu'il a exactement
deux sous-groupes normaux: $\{1\}$ et $G$.
\end{defi}

\exs $(1)$ Les seuls groupes \emph{ab\'{e}liens} simples sont les
groupes cycliques d'ordre premier, c'est-\`{a}-dire les groupes
$\ZM/{p\ZM}$ avec $p$ premier.

$(2)$ Le groupe altern\'{e} $\AM_n$ est simple si $n\geqslant 5$.

$(3)$ Le groupe $\mathbf{PSL}_n(\FM_q)$ est simple pour
$n\geqslant 2$ sauf dans le cas $n=2$ et $q=2$ ou $3$.

\section{Filtrations et th\'{e}or\`{e}me de Jordan-H\"{o}lder}
\label{filtrations}
\begin{defi}
Une \emph{filtration} du groupe $G$ est une suite finie
${(G_i)}_{0\leqslant i \leqslant n}$ de sous-groupes telle que
$$G_0=\{1\}\subset G_1\subset \cdots \subset G_i \subset \cdots \subset G_n=G$$
avec $G_i$ normal dans $G_{i+1}$, pour $0\leqslant i\leqslant
n-1$.

On appelle \emph{gradu\'{e}} de $G$ (associ\'{e} \`{a} la filtration
${(G_i)}_{0\leqslant i \leqslant n}$) et on note $\gr(G)$ la suite
des $\gr_i(G)=G_i/G_{i-1}$, pour $1\leqslant i\leqslant n$.
\end{defi}

\begin{defi}
Une filtration ${(G_i)}_{0\leqslant i \leqslant n}$ de $G$ est dite
\emph{de Jordan-H\"{o}lder} si $G_i/G_{i-1}$ est simple pour tout
$1\leqslant i\leqslant n$.
\end{defi}

\begin{prop}
Si $G$ est fini, $G$ poss\`{e}de une suite de Jordan-H\"{o}lder.
\end{prop}

Si $G=\{1\}$, on a la suite de Jordan-H\"{o}lder triviale ($n=0$).
Si $G$ est simple, on prend $n=1$. Si $G$ n'est pas simple, on
raisonne par r\'{e}currence sur l'ordre de $G$. Soit $N\subset G$,
$N$ normal dans $G$ d'ordre maximal. Alors $G/N$ est simple, car
sinon il existerait $M$ normal dans $G$ contenant strictement $N$ et
distinct de $G$. Comme $|N|<|G|$, on peut appliquer l'hypoth\`{e}se
de r\'{e}currence et si ${(N_i)_{0\leqslant i\leqslant n}}$ est une
suite de Jordan-H\"{o}lder pour $N$, alors $(N_0, \cdots, N_{n-1},
N, G)$ en est une pour $G$.~\findem

\bigskip\rmq Si $G$ est infini, il peut ne pas poss\'{e}der de suite de
Jordan-H\"{o}lder: c'est par exemple le cas de $\ZM$.

\begin{thm}[Jordan-H\"{o}lder]
Soit $G$ un groupe fini et soit ${(G_i)}_{0\leqslant i \leqslant n}$
une suite de Jordan-H\"{o}lder de $G$. Le gradu\'{e} de $G$, \`{a}
permutation pr\`{e}s des indices, ne d\'{e}pend pas de la suite
choisie.
\end{thm}

Il suffit de montrer que si $S$ est un groupe simple fix\'{e} et si
$n\left(G,(G_i),S\right)$ est le nombre de $j$ tels que
$G_j/G_{j-1}$ est isomorphe \`{a} $S$, alors
$n\left(G,(G_i),S\right)$ ne d\'{e}pend pas de la suite $(G_i)$.

On commence par une remarque: si $H$ est un sous-groupe de $G$, une
filtration $(G_i)$ sur $G$ induit une filtration $(H_i)$ sur $H$
d\'{e}finie par $H_i=G_i\cap H$. De m\^{e}me, si $N$ est normal, on
a une filtration sur $G/N$ d\'{e}finie par $(G/N)_i=G_i/(G_i\cap
N)$. La suite exacte
$$\xymatrix{
\{1\} \ar[r] & N \ar[r] & G \ar[r] & G/N \ar[r] & \{1\} }$$ se
conserve par filtration:
$$\xymatrix{
\{1\} \ar[r] & N_i/{N_{i-1}} \ar[r] & G_i/G_{i-1} \ar[r] &
(G/N)_i/(G/N)_{i-1} \ar[r] & \{1\}}$$ d'o\`{u} finalement la suite
exacte
$$\xymatrix{
\{1\} \ar[r] & \gr_i(N) \ar[r] & \gr_i(G) \ar[r] & \gr_i(G/N)
\ar[r] & \{1\}. }$$

Si la filtration initiale est de Jordan-H\"{o}lder, $\gr_i(G)$ est
simple pour tout $i$, donc $\gr_i(N)$ est isomorphe \`{a} $\{1\}$ ou
\`{a} $\gr_i(G)$. Par r\'{e}indexation, on peut donc obtenir une
filtration de Jordan-H\"{o}lder sur $N$ et de m\^{e}me sur $G/N$.

Cette remarque permet de d\'{e}montrer le th\'{e}or\`{e}me; on a en
effet deux possibilit\'{e}s: soit $\gr_i(N)=\{1\}$ et
$\gr_i(G/N)=\gr_i(G)$, soit $\gr_i(N)=\gr_i(G)$ et
$\gr_i(G/N)=\{1\}$. On en d\'{e}duit une partition de $I=\{0, \dots,
n\}$ en deux parties: $I_1=\left\{i,\, \gr_i(N)=\{1\}\right\}$ et
$I_2=\left\{i,\,\gr_i(N)=\gr_i(G)\right\}$.\\
On raisonne alors par r\'{e}currence sur l'ordre de $G$. Si
$G=\{1\}$, il n'y a pas de probl\`{e}me. Sinon, on peut toujours
supposer que $G$ n'est pas simple. Soit alors $N$ un sous-groupe
normal tel que $|N|<|G|$ et $|G/N|<|G|$. L'hypoth\`{e}se de
r\'{e}currence s'applique \`{a} $N$ et $G/N$:
$n\left(N,{(N_i)}_{i\in I_2},S\right)$ et
$n\left(G/N,\left((G/N)_i\right)_{i\in I_1},S\right)$ sont
ind\'{e}pendants de la filtration. Or
$$n\left(G,{(G_i)}_{i\in I},S\right)=n\left(N,{(N_i)}_{i\in
I_2},S\right)+n\left(G/N,\left((G/N)_i\right)_{i\in I_1},S\right),$$
donc $n\left(G,{(G_i)}_{i\in I},S\right)$ est ind\'{e}pendant de la
filtration choisie.~\findem

\bigskip\App On retrouve ainsi l'unicit\'{e} de la d\'{e}composition d'un entier
en produit de facteurs premiers. En effet, si $n=p_1^{h_1}\cdots
p_k^{h_k}$, on a pour $\ZM/n\ZM$ la filtration de Jordan-H\"{o}lder
suivante:
$$\ZM/n\ZM\supset p_1 \ZM/n\ZM \supset p_1^2 \ZM/n\ZM\supset\cdots\supset p_1^{h_1} \ZM/n\ZM\supset\cdots$$

Donc $\ZM/p_i\ZM$ appara\^{i}t $h_i$ fois dans le gradu\'{e},
d'o\`{u} l'unicit\'{e}.

\bigskip\exs
$(1)$ Filtration de $\SM_3$: $\AM_3$ est normal dans $\SM_3$ et
$\AM_3$ est cyclique d'ordre $3$. D'o\`{u} la filtration
$$\{1\}\subset \AM_3 \subset \SM_3.$$

$(2)$ Filtration de $\SM_4$: $\AM_4$ est normal dans $\SM_4$ avec
$(\SM_4:\AM_4)=2$. Dans $\AM_4$, il existe un sous-groupe normal
$D$ de type $(2,2)$: $D=\{1, \sigma_1, \sigma_2, \sigma_3\}$ avec
\begin{eqnarray*}
\sigma_1 & = & (a,b)(c,d),\\
\sigma_2 & = & (a,c)(b,d),\\
\sigma_3 & = & (a,d)(b,c).
\end{eqnarray*}

On a donc la filtration
$$\{1\} \subset \{1, \sigma_i\} \subset D \subset \AM_4 \subset \SM_4.$$

L'ordre des quotients successifs est\, $2,2,3,2$. Le choix de $i$
\'{e}tant arbitraire, il n'y a pas unicit\'{e} de la filtration.

$(3)$ Filtration de $\SM_n$ pour $n\geqslant 5$: le groupe $\AM_n$
\'{e}tant simple, on a la filtration
$$\{1\}\subset \AM_n \subset \SM_n.$$

\chapter{Th\'{e}or\`{e}mes de Sylow}

Soit $p$ un nombre premier et soit $G$ un groupe fini.

\section{D\'{e}finitions}\label{sylow}

\begin{defi}
On dit que $G$ est un \emph{$p$-groupe} si l'ordre de $G$ est une
puissance de $p$. Si $G$ est d'ordre $p^n m$ avec $m$ premier \`{a}
$p$, on dit qu'un sous-groupe $H$ de $G$ est un \emph{$p$-Sylow} de
$G$ si $H$ est d'ordre $p^n$.
\end{defi}

\rmqs $(1)$ Soit $S$ un sous-groupe de $G$; $S$ est un $p$-Sylow de
$G$ si et seulement si $S$ est un $p$-groupe et $(G:S)$ est premier
\`{a} $p$.

$(2)$ Tout conjugu\'{e} d'un $p$-Sylow de $G$ est un $p$-Sylow de
$G$.

\bigskip\ex Soit $K$ un corps fini de caract\'{e}ristique
$p$ \`{a} $q=p^f$ \'{e}l\'{e}ments. Soit $G=\mathbf{GL}_n(K)$ le
groupe des matrices inversibles $n\times n$ \`{a} coefficients dans
$K$. Ce groupe est isomorphe \`{a} $\mathbf{GL}(V)$ o\`{u} $V$ est
un espace vectoriel sur $K$ de dimension $n$. On remarque que
l'ordre de $G$ est le nombre de bases d'un espace vectoriel de
dimension $n$ sur $K$, soit:
$$|G|= (q^n-1)(q^n-q)\cdots (q^n-q^{n-1})=
q^{n(n-1)/2}\prod_{i=1}^{n}{(q^i-1)} = p^{fn(n-1)/2} m,$$ o\`{u}
$m=\prod_{i=1}^{n}{(q^i-1)}$ est premier \`{a} $q$, donc \`{a} $p$.\\
Consid\'{e}rons d'autre part le groupe $P$ constitu\'{e} des
matrices triangulaires sup\'{e}rieures \`{a} coefficients diagonaux
\'{e}gaux \`{a} $1$. C'est un sous-groupe de $G$ d'ordre
$|P|=q^{n(n-1)/2}=p^{fn(n-1)/2}$. Donc $P$ est un $p$-Sylow de $G$.

\section{Existence des $p$-Sylow}\label{2.2}

Le but de cette section est de d\'{e}montrer le premier
th\'{e}or\`{e}me de Sylow:

\begin{thm}\label{sylow 1}
Tout groupe fini poss\`{e}de au moins un $p$-Sylow.
\end{thm}

\subsection{Premi\`{e}re d\'{e}monstration}
Elle repose sur la proposition suivante:

\begin{prop}
Soit $H$ un sous-groupe de $G$ et soit $S$ un $p$-Sylow de $G$.
Alors il existe $g\in G$ tel que $H\cap gSg^{-1}$ soit un
$p$-Sylow de $H$.
\end{prop}

Soit $X$ l'ensemble des classes \`{a} gauche de $G$ modulo $S$. Le
groupe $G$ (resp. $H$) agit sur $X$ par translations. Les
stabilisateurs des points de $X$ sous $G$ (resp. sous $H$) sont les
conjugu\'{e}s de $S$ (resp. les $H\cap gSg^{-1}$). Or $|X|\not\equiv
0{\pmod p}$ car $S$ est un $p$-Sylow de $G$. L'une des orbites
$\mathcal O$ de $X$ sous l'action de $H$ a un nombre
d'\'{e}l\'{e}ments premier \`{a} $p$ (sinon $|X|$ serait divisible
par $p$); soit $x\in \mathcal O$ et soit $H_x$ le stabilisateur de
$x$ dans $H$. Le groupe $H_x$ est un $p$-groupe, de la forme $H\cap
gSg^{-1}$ (pour un certain $g$) et $(H:H_x)=|{\mathcal O}|$ est
premier \`{a} $p$. Donc $H_x$ est un $p$-Sylow de $H$ de la forme
$H\cap gSg^{-1}$.~\findem

\begin{coro}
Si $G$ a des $p$-Sylow et si $H$ est un sous-groupe de $G$, alors
$H$ a aussi des $p$-Sylow.
\end{coro}

\App[Une premi\`{e}re preuve du th. \ref{sylow 1}] Soit $G$ un
groupe fini d'ordre $n$. On peut plonger $G$ dans le groupe
sym\'{e}trique $\SM_n$. D'autre part, $\SM_n$ se plonge dans
$\mathbf{GL}_n(K)$ (o\`{u} $K$ est un corps fini de
caract\'{e}ristique $p$): si $\sigma\in \SM_n$ et si
$(e_i)_{1\leqslant i\leqslant n}$ est une base de $K^n$, on associe
\`{a} $\sigma$ la transformation lin\'{e}aire $f$ d\'{e}finie par
$f(e_i)=e_{\sigma(i)}$. Donc $G$ se plonge dans $\mathbf{GL}_n(K)$.
D'apr\`{e}s l'exemple du \S\ \ref{sylow}, $\mathbf{GL}_n(K)$
poss\`{e}de un $p$-Sylow. Le corollaire ci-dessus permet de
conclure.

\subsection{Seconde d\'{e}monstration (Miller-Wielandt)}\label{Miller-Wielandt}

On suppose que $|G|=p^n m$ avec $m$ premier \`{a} $p$. On note $X$
l'ensemble des parties de $G$ \`{a} $p^n$ \'{e}l\'{e}ments et $s$ le
nombre de $p$-Sylow de $G$.

\begin{lemme}
$|X|\equiv sm{\pmod p}$.
\end{lemme}

Le groupe $G$ op\`{e}re sur $X$ par translations \`{a} gauche. Soit
$X=\coprod_{i} X_i$ la d\'{e}composition de $X$ en orbites sous
l'action de $G$. Si $A_i\in X_i$, on a $X=\coprod_{i} GA_i$. On note
$G_i$ le stabilisateur de $A_i$. On rappelle que
$|GA_i|=|G|/|G_i|$.\\
Remarque: $|G_i|\leqslant p^n$. Soit en effet $x\in A_i$. Si $g\in
G_i$, alors $gx$ appartient \`{a} $A_i$, donc peut prendre $p^n$
valeurs. On a donc au plus $p^n$ choix pour $g$. On distingue donc
deux cas:\\
$\bullet$ Si $|G_i|<p^n$, alors $|GA_i|$ est divisible par $p$.\\
$\bullet$ Si $|G_i|=p^n$, alors $G_i$ est un $p$-Sylow de $G$.\\
R\'{e}ciproquement soit $P$ un $p$-Sylow de $G$; $Pg\in X$ pour tout
$g\in G$ et le stabilisateur de $Pg$ est $P$. De m\^{e}me, si $P$
stabilise une partie $A$ de $X$, alors $PA\subset A$ donc pour tout
$a\in A$, on a $Pa\subset A$ donc $A=Pa$. (les deux ensembles ont
m\^{e}mes cardinaux). Donc $P$ stabilise exactement son orbite sous
l'action de $G$ (et le cardinal de cette orbite est $|G/P|=m$).
Finalement
$$|X|={\sum_{i\,/\,|G_i|<p^n}\!\!\!\!{|GA_i|}}\ \ +{\sum_{i\,/\,|G_i|=p^n}\!\!\!\!{|GA_i|}}$$
soit
$$|X|\equiv 0+sm{\pmod p}$$
d'o\`{u} le r\'{e}sultat.~\findem

\bigskip
Ce lemme nous donne le th. \ref{sylow 1}. En effet, d'apr\`{e}s ce
lemme, la classe de $s$ modulo $p$ ne d\'{e}pend que de l'ordre de
$G$. Or $G'=\ZM/|G|\ZM$ a un unique $p$-Sylow (qui est isomorphe
\`{a} $\ZM/p^n\ZM$). Donc $s\equiv 1{\pmod p}$; en particulier, $s$
est non nul.~\findem

\bigskip\rmq On a d\'{e}montr\'{e} en fait que le nombre de $p$-Sylow d'un groupe
$G$ est congru \`{a} $1$ modulo $p$. On retrouvera cette
propri\'{e}t\'{e} ult\'{e}rieurement.
\begin{coro}[Cauchy]
Si $p$ divise l'ordre de $G$, alors $G$ contient un \'{e}l\'{e}ment
d'ordre~$p$.
\end{coro}

En effet, soit $S$ un $p$-Sylow de $G$ (il en existe d'apr\`{e}s le
th. \ref{sylow 1}); $S$ n'est pas r\'{e}duit \`{a} $\{1\}$ car $p$
divise l'ordre de $G$. Soit $x\in S$ distinct de $\{1\}$. L'ordre de
$x$ est une puissance de $p$, soit $p^m$ ($m\geqslant 1$). Alors
$x^{p^{m-1}}$ est d'ordre $p$.~\findem

\section{Propri\'{e}t\'{e}s des $p$-Sylow}\label{2.3}

\begin{thm}[Second th\'{e}or\`{e}me de Sylow]\label{sylow 2}

\

\begin{enumerate}
 \item[(1)] Tout $p$-sous-groupe de $G$ est contenu dans un $p$-Sylow de
$G$.
 \item[(2)] Les $p$-Sylow de $G$ sont conjugu\'{e}s.
 \item[(3)] Le nombre des $p$-Sylow est congru \`{a} $1$ modulo $p$.
\end{enumerate}
\end{thm}

\begin{lemme}\label{lemme1}
Soit $X$ un ensemble fini sur lequel op\`{e}re un $p$-groupe $P$ et
soit $X^P$ l'ensemble des \'{e}l\'{e}ments de $X$ fix\'{e}s par $P$.
Alors $|X|\equiv |X^P|\pmod p$.
\end{lemme}

Les orbites \`{a} un \'{e}l\'{e}ment de $X$ sous l'action de $P$
sont celles constitu\'{e}es d'un point de $X^P$. L'ensemble $X-X^P$
est donc r\'{e}union d'orbites non triviales, de cardinal divisible
par $p$.~\findem

\bigskip
On peut alors d\'{e}montrer les points $(1)$ et $(2)$ du th.
\ref{sylow 2}: soit $S$ un $p$-Sylow de $G$ et soit $P$ un
$p$-sous-groupe de $G$. On applique le lemme \ref{lemme1} \`{a}
l'ensemble $X$ des classes \`{a} gauche de $G$ modulo $S$:
$|X|\not\equiv 0 \pmod p$ donc $|X^P|\not\equiv 0 \pmod p$. En
particulier, il existe $x\in X$ fix\'{e} par $P$. Le stabilisateur
de $x$ contient donc $P$ et est un conjugu\'{e} de $S$. Donc $P$ est
contenu dans un conjugu\'{e} de $S$
(c'est-\`{a}-dire dans un $p$-Sylow de $G$).\\
Pour le point $(2)$, on applique $(1)$ \`{a} $P=S'$ o\`{u} $S'$ est
un $p$-Sylow de $G$. Il existe $g\in G$ tel que $S'\subset
gSg^{-1}$, donc $S'=gSg^{-1}$.~\findem

\bigskip
Pour le point $(3)$, on donne une nouvelle d\'{e}monstration
bas\'{e}e sur le lemme suivant:

\begin{lemme}\label{lemme2}
Soient $S$ et $S'$ deux $p$-Sylow de $G$. Si $S'$ normalise $S$,
alors $S=S'$.
\end{lemme}

Soit $H$ le sous-groupe de $G$ engendr\'{e} par $S$ et $S'$. Le
groupe $H$ normalise $S$ qui est un $p$-Sylow de $H$. Donc $S$ est
le seul $p$-Sylow de $H$ (les $p$-Sylow de $H$ sont conjugu\'{e}s);
or $S'$ est un $p$-Sylow de $H$ donc $S=S'$.~\findem

\bigskip
Montrons le  point $(3)$. Si $X$ est l'ensemble des $p$-Sylow de
$G$, alors $S$ op\`{e}re sur $X$ par conjugaison et d'apr\`{e}s le
lemme \ref{lemme2}, $S$ est le seul \'{e}l\'{e}ment de $X$ fix\'{e}
par $S$. On r\'{e}applique le lemme \ref{lemme1} (avec $P=S$):
$|X|\equiv 1 \pmod p$.~\findem

\begin{coro}
Si $S$ est un $p$-Sylow de $G$, alors $\big(G:N_G(S)\big)\equiv 1
\pmod p$.
\end{coro}

L'application $f$ de $G/N_G(S)$ dans l'ensemble des $p$-Sylow de
$G$, d\'{e}finie par $f({\bar g})=gSg^{-1}$ (o\`{u} $g$ est un
repr\'{e}sentant quelconque de $\bar g$) est
bijective.~\findem\bigskip

On a vu que pour tout sous-groupe $H$ de $G$, il existe un
$p$-Sylow de $G$ dont l'intersection avec $H$ est un $p$-Sylow de
$H$. Ce n'est pas vrai pour tout $p$-Sylow de $G$. Mais si $H$ est
normal, on a:

\begin{prop}\label{prop2.10}
Soit $H$ un sous-groupe normal de $G$ et soit $S$ un $p$-Sylow de
$G$. Alors
\begin{enumerate}
 \item[(1)] $S\cap H$ est un $p$-Sylow de $H$.
 \item[(2)] L'image de $S$ dans $G/H$ est un $p$-Sylow de $G/H$ (et
on les obtient tous ainsi).
 \item[(3)] (\,{\bf Frattini}) Si $Q$ est un $p$-Sylow de $H$, alors $H.N_G(Q)=G$.
\end{enumerate}
\end{prop}\label{thmFrattini}

$(1)$ Evident.

$(2)$ L'image de $S$ dans $G/H$ est isomorphe \`{a} $S/(H\cap S)$.
Si $p^a$ (resp. $p^b$) est la puissance de $p$ maximale divisant
l'ordre de $H$ (resp. de $G/H$), $p^{a+b}$ est la puissance maximale
de $p$ divisant l'ordre de $G$. Par suite, $S$ a $p^{a+b}$
\'{e}l\'{e}ments. De plus, $H\cap S$ a au plus $p^a$
\'{e}l\'{e}ments donc $S/(H\cap S)$ au moins $p^b$ et donc
exactement $p^b$. Il s'ensuit que $S/(H\cap S)$ est un $p$-Sylow de
$G/H$. D'autre part, on obtient tous les $p$-Sylow par conjugaison,
d'o\`{u} $(2)$.

$(3)$ Soit $g\in G$. On a $gQg^{-1}\subset gHg^{-1}=H$ ($H$ est
normal). Or $gQg^{-1}$ est un $p$-Sylow de $H$, donc il existe
$h\in H$ tel que $gQg^{-1}=hQh^{-1}$, donc $h^{-1}g\in N_G(Q)$ et
donc $g\in H.N_G(Q)$. Ainsi $G\subset H.N_G(Q)$, donc
$H.N_G(Q)=G$.~\findem

\begin{coro}\label{2.11}
Soit $S$ un $p$-Sylow de $G$ et soit $H$ un sous-groupe de $G$
contenant $N_G(S)$. Alors $N_G(H)=H$.
\end{coro}

Le groupe $H$ est normal dans $N_G(H)$ et contient $S$ qui est donc un
$p$-Sylow de $H$. On applique le point $(3)$ de la proposition
ci-dessus: $H.N_G(S)=N_G(H)$. Donc $N_G(H)\subset H$ d'o\`{u} le
r\'{e}sultat.~\findem\bigskip

En particulier, si $S$ est un $p$-Sylow de $G$, on a
$N_G\big(N_G(S)\big)=N_G(S)$.

\section{Fusion}\label{2.4}

Soit $S$ un $p$-Sylow de $G$. On note $N$ le normalisateur de $S$
dans $G$. On se pose le probl\`{e}me de savoir si deux
\'{e}l\'{e}ments de $S$ conjugu\'{e}s dans $G$ sont conjugu\'{e}s
dans $N$. On a la: \label{Burnside1}
\begin{prop}[Burnside]
Soient $X$ et $Y$ deux parties du centre de $S$, conjugu\'{e}es dans
$G$ et soit $g\in G$ tel que $gXg^{-1}=Y$. Alors il existe $n\in N$
tel que $nxn^{-1}=gxg^{-1}$ pour tout $x\in X$. En particulier,
$nXn^{-1}=Y$.
\end{prop}

On veut trouver $n\in N$ tel que $nxn^{-1}=gxg^{-1}$ pour tout $x\in
X$ i.e. $g^{-1}nxn^{-1}g=x$ pour tout $x\in X$. Donc on cherche
$n\in N$ tel que $g^{-1}n\in A=C_G(X)$ (le centralisateur de $X$).
Or $X$ est contenu dans le centre de $S$ donc $A$ contient $S$. De
m\^{e}me, $Y=gXg^{-1}$ donc $g^{-1}Sg$ est contenu dans $A$. Les
groupes $S$ et $g^{-1}Sg$ sont des $p$-Sylow de $A$ (il suffit de
regarder leurs ordres) donc sont conjugu\'{e}s dans $A$: il existe
$a\in A$ tel que $ag^{-1}Sga^{-1}=S$. Donc $n=ga^{-1}$ appartient
\`{a} $N$ et $g^{-1}n$ appartient \`{a} $A$.~\findem

\begin{coro} Soient $x$ et $y$ deux \'{e}l\'{e}ments du centre de $S$.
S'ils sont conjugu\'{e}s dans $G$, ils sont conjugu\'{e}s dans $N$.
\end{coro}

\rmq L'hypoth\`{e}se \og $x$ et $y$ appartiennent au centre de
$S$\fg \ ne peut \^{e}tre supprim\'{e}e: si l'on prend
$G=\mathbf{GL}_3(\ZM/p\ZM)$ et
$$S=\left\{\left(
\begin{array}{ccc}
  1 & \times & \times \\
  0 & 1 & \times \\
  0 & 0 & 1 \\
\end{array}
\right)\right\}$$ alors
$$N=\left\{\left(
\begin{array}{ccc}
  \times & \times & \times \\
  0 & \times & \times \\
  0 & 0 & \times \\
\end{array}
\right)\right\}$$ et les \'{e}l\'{e}ments
$$x=\left(
\begin{array}{ccc}
  1 & 1 & 0 \\
  0 & 1 & 0 \\
  0 & 0 & 1 \\
\end{array}
\right)\;\;\;\mbox{et}\;\;\; y=\left(
\begin{array}{ccc}
  1 & 0 & 0 \\
  0 & 1 & 1 \\
  0 & 0 & 1 \\
\end{array}
\right)$$ sont conjugu\'{e}s dans $G$ et ne le sont pas dans $N$.

\

Disons que deux \'{e}l\'{e}ments $x,y$ de $S$ sont \emph{localement
conjugu\'{e}s} s'il existe un sous-groupe $U$ de $S$ les contenant
tel que $x$ et $y$ soient conjugu\'{e}s dans $N_G(U)$.

\begin{thm}[Alperin]\label{Alperin}
La relation d'\'{e}quivalence sur $S$ engendr\'{e}e par la relation
\og $x$ et $y$ sont localement conjugu\'{e}s \fg\ est la relation
\og $x$ et $y$ sont conjugu\'{e}s dans $G$ \fg.
\end{thm}

En d'autres termes: \pagebreak[2]
\begin{thm}
Si $x,y\in S$ sont conjugu\'{e}s dans $G$, il existe une suite $a_0,
\dots, a_n$ d'\'{e}l\'{e}ments de $S$ telle que:
\begin{enumerate}
 \item[(1)] $a_0=x$ et $a_n=y$.
 \item[(2)] $a_i$ est localement conjugu\'{e} de $a_{i+1}$ pour
$0\leqslant i\leqslant n-1$.
\end{enumerate}
\end{thm}

Cela r\'{e}sulte du th\'{e}or\`{e}me plus pr\'{e}cis suivant:

\begin{thm}
Soit $A$ une partie de $S$ et soit $g\in G$ tel que $A^g\subset S$.
Il existe alors un entier $n\geqslant 1$, des sous-groupes $U_1,
\dots, U_n$ de $S$ et des \'{e}l\'{e}ments $g_1, \dots, g_n$ de $G$
avec:
\begin{enumerate}
 \item[(1)] $g=g_1\cdots g_n$.
 \item[(2)]$g_i\in N_G(U_i)$ pour $1\leqslant i\leqslant n$.
 \item[(3)] $A^{g_1\cdots g_{i-1}}\subset U_i$ pour $1\leqslant
i\leqslant n$.
\end{enumerate}
\end{thm}

(Dans cet \'{e}nonc\'{e}, $A^g$  d\'{e}signe  $g^{-1}Ag$.)

\rmq Pour $i=1$, $(3)$ signifie que $A\subset U_1$. Noter que l'on
a $A^{g_1\cdots g_i}\subset U_i$ pour $1\leqslant i\leqslant n$,
comme on le voit en combinant $(2)$ et $(3)$. On a en particulier
$A^g\subset U_n$.

\bigskip Le th\'{e}or\`{e}me ci-dessus est un corollaire de celui-ci
(prendre $A$ r\'{e}duit \`{a} un \'{e}l\'{e}ment).

\bigskip
\dem Soit $T$ le sous-groupe de $S$ engendr\'{e} par $A$. On
raisonne par r\'{e}currence sur \emph{l'indice} $(S:T)$ de $T$ dans
$S$. Si cet indice est $1$, on a $T=S$, d'o\`{u} $S^g=S$ et $g\in
N_{G}(S)$. On prend alors $n=1$, $g_1=g$ et $U_1=S$.

Supposons donc $(S:T)>1$, i.e. $T\neq S$. Le groupe $T_1=N_S(T)$ est
alors distinct de $T$. C'est un $p$-sous-groupe de $N_G(T)$.
Choisissons un $p$-Sylow $\Sigma$ de $N_G(T)$ contenant $T_1$.
D'apr\`{e}s le th. \ref{sylow 2}, il existe $u\in G$ tel que
$\Sigma^u\subset S$. Posons d'autre part $V=T^g$; on a $V\subset S$
par hypoth\`{e}se. Le groupe $\Sigma^g$ est un $p$-Sylow
de $N_G(V)=\big(N_G(T)\big)^g$.\\
Comme $N_S(V)$ est un $p$-sous-groupe de $N_G(V)$, il existe $w\in
N_G(V)$ tel que $\big(N_S(V)\big)^w\subset\Sigma^g$. Posons
$v=u^{-1}gw^{-1}$. On a $g=uvw$.

On va maintenant d\'{e}composer $u$ et $v$:\\
(i) On a $T^u\subset\Sigma^u\subset S$. Comme l'indice de $T_1$ dans
$S$ est strictement inf\'{e}rieur \`{a} celui de $T$,
l'hypoth\`{e}se de r\'{e}currence montre qu'il existe des
sous-groupes $U_1,\dots,U_m$ de $S$ et des \'{e}l\'{e}ments $u_1\in
N_G(U_1),\dots,u_m\in N_G(U_m)$ avec $u=u_1\cdots u_m$ et
$T_1^{u_1\cdots u_{i-1}}\subset U_i$ pour
$1\leqslant i\leqslant m$.\\
(ii) Posons $T_2=N_S(V)$ et $T_3=T_2^{v^{-1}}=T_2^{wg^{-1}u}$. Comme
$T_2^w$ est contenu dans $\Sigma^g$, on a
$T_3\subset\Sigma^{gg^{-1}u}=\Sigma^u$. Le groupe $T_3$ est contenu
dans $S$, et $T_3^v=T_2$ aussi. Comme l'indice de $T_3$ est
strictement inf\'{e}rieur \`{a} celui de $T$, on en d\'{e}duit comme
ci-dessus l'existence de sous-groupes $V_1,\dots,V_r$ de $S$ et
d'\'{e}l\'{e}ments $v_j\in N_G(V_j)$, avec $v=v_1\cdots v_r$ et
$T_3^{v_1\cdots v_{j-1}}\subset V_j$ pour $1\leqslant j\leqslant r$.

Il reste \`{a} v\'{e}rifier que les sous-groupes
$U_1,\dots,U_m,V_1,\dots,V_r,V$ de $S$ et la d\'{e}composition
$g=u_1\cdots u_m v_1\cdots v_r w$ de $g$ satisfont aux conditions
du th\'{e}or\`{e}me.\\
On a $$u_i\in N_G(U_i),\; v_j\in N_G(V_j),\; w\in N_G(V)$$ par
construction, ainsi que $T^{u_1\cdots u_{i-1}}\subset U_i$
($1\leqslant i\leqslant m$) puisque $T$ est contenu dans $T_1$. Il
reste \`{a} voir que $$T^{u_1\cdots u_m v_1\cdots v_{j-1}}\subset
V_j$$ pour $1\leqslant j\leqslant r$.\\
Or $T^{u_1\cdots u_m}=T^u$ est contenu dans $T_3=T_2^{wg^{-1}u}$; en
effet, $V=T^g$ est normalis\'{e} par $w^{-1}$; on a donc
$T^{gw^{-1}}=V\subset N_S(V)=T_2$, d'o\`{u} $T\subset T_2^{wg^{-1}}$
et $T^u\subset T_2^{wg^{-1}u}$.\\
On d\'{e}duit de l\`{a} que
$$T^{u_1\cdots u_m v_1\cdots v_{j-1}}=T^{u v_1\cdots
v_{j-1}}\subset T_3^{v_1\cdots v_{j-1}}\subset V_j,$$ ce qui
ach\`{e}ve la d\'{e}monstration.~\findem

\chapter{Groupes r\'{e}solubles et groupes nilpotents}\label{chap3}

\section{Groupes r\'{e}solubles}\label{grpe resol}

Soit $G$ un groupe et soient $x,y$ deux \'{e}l\'{e}ments de $G$.
L'\'{e}l\'{e}ment $x^{-1}y^{-1}xy$ est appel\'{e} le
\emph{commutateur}\label{commutateur} de $x$ et $y$. On le note
$(x,y)$. On a
$$xy=yx(x,y).$$

Si $A$ et $B$ sont deux sous-groupes de $G$, on note $(A,B)$ le
groupe engendr\'{e} par les commutateurs $(x,y)$ avec $x\in A$ et
$y\in B$. Le groupe $(G,G)$ est appel\'{e} le \emph{groupe des
commutateurs} de $G$ ou encore le \emph{groupe d\'{e}riv\'{e}} de
$G$ et est not\'{e} $D(G)$. C'est un sous-groupe caract\'{e}ristique
de $G$. De sa d\'{e}finition r\'{e}sulte aussit\^{o}t la:

\begin{prop}
Soit $H$ un sous-groupe de $G$. Les propri\'{e}t\'{e}s suivantes
sont \'{e}quivalentes:
\begin{enumerate}
 \item[(1)] $H$ contient $D(G)$.
 \item[(2)] $H$ est normal et $G/H$ est ab\'{e}lien.
\end{enumerate}
\end{prop}

Ainsi $G/D(G)$ est le plus grand quotient ab\'{e}lien de $G$.
On le note parfois $G^{ab}$.

On peut it\'{e}rer le proc\'{e}d\'{e} et d\'{e}finir la suite des
sous-groupes d\'{e}riv\'{e}s de $G$:
\begin{eqnarray*}
D^0G & = & G,\\
D^nG & = & (D^{n-1}G,D^{n-1}G)\quad \mbox{pour } n\geqslant 1.
\end{eqnarray*}

On a $G\supset D^1G\supset D^2G\supset \cdots$.

\label{classeresolubilite}
\begin{defi}
Un groupe $G$ est dit \emph{r\'{e}soluble} s'il existe un entier
$n\geqslant 0$ tel que $D^nG=\{1\}$. On appelle alors \emph{classe
de r\'{e}solubilit\'{e}} de $G$ et on note $cl(G)$ le plus petit
entier $n$ positif pour lequel $D^nG=\{1\}$.
\end{defi}

Ainsi, $cl(G)=0$ \'{e}quivaut \`{a} $G=\{1\}$ et $cl(G)\leqslant 1$
\'{e}quivaut \`{a} dire que $G$ est ab\'{e}lien.

\begin{prop}
Soit $G$ un groupe et soit $n$ un entier $\geqslant 1$. Les
propri\'{e}t\'{e}s suivantes sont \'{e}quivalentes:
\begin{enumerate}
\item[(1)] $G$ est r\'{e}soluble de classe $\leqslant n$,

\item[(2)] Il existe une suite $G=G_0\supset G_1\supset \cdots
\supset G_n=\{1\}$ de sous-groupes normaux de $G$ tels que
$G_i/G_{i+1}$ soit ab\'{e}lien pour $0\leqslant i \leqslant n-1$,

\item[(2')] Il existe une suite $G=G_0\supset G_1\supset \cdots
\supset G_n=\{1\}$ de sous-groupes de $G$ tels que $G_i$ soit normal
dans $G_{i-1}$  et que $G_{i-1}/G_i$ soit ab\'{e}lien, pour
$1\leqslant i \leqslant n$,

\item[(3)] Il existe un sous-groupe ab\'{e}lien $A$ normal
 dans $G$ tel que $G/A$ soit r\'{e}soluble de classe $\leqslant n-1$.
\end{enumerate}
\end{prop}

$(1)\Rightarrow(2)$ Posons $G_i=D^iG$ pour tout $i\geqslant0$.
Puisque $D(G)$ est stable par tout automorphisme (m\^{e}me non
int\'{e}rieur!) de $G$, $D^iG$ est normal dans $G$ pour tout $i$. La
suite ${(G_i)}_{i\geqslant 0}$ ainsi d\'{e}finie v\'{e}rifie donc
$(2)$.

$(2)\Rightarrow(2')$ est trivial.

$(2')\Rightarrow(1)$ Par r\'{e}currence sur $k$ on voit que
$D^kG\subset G_k$ pour tout $k$, d'o\`{u} $D^nG=\{1\}$.

$(1)\Rightarrow(3)$ On prend $A=D^{n-1}G$.

$(3)\Rightarrow(1)$ D'apr\`{e}s l'implication $(1)\Rightarrow(2)$,
appliqu\'{e}e \`{a} $G/A$ et \`{a} $n-1$, il existe une suite
$$A_0=G\supset A_1 \cdots \supset A_{n-1}=A$$ de sous-groupes normaux de $G$ telle que la
suite des quotients
$$G/A\supset A_1/A\supset \cdots \supset
A_{n-1}/A=\{1\}$$ v\'{e}rifie la condition $(2)$. Alors la suite
$$G\supset A_1\supset \cdots \supset A_{n-1}\supset \{1\}$$ v\'{e}rifie
la condition $(2)$ et l'implication $(2)\Rightarrow(1)$
appliqu\'{e}e \`{a} $G$ et \`{a} $n$ permet de conclure.~\findem

\bigskip
\rmq Tout sous-groupe (et tout groupe quotient) d'un groupe
r\'{e}soluble de classe~$\leqslant~n$ est r\'{e}soluble de classe
$\leqslant n$.

\begin{prop}
Soit $G$ un groupe fini et soit $G=G_0\supset G_1 \supset \cdots
\supset G_n=\{1\}$ une suite de Jordan-H\"{o}lder de $G$. Pour que
$G$ soit r\'{e}soluble, il faut et il suffit que $G_i/G_{i+1}$ soit
cyclique d'ordre premier pour $0\leqslant i\leqslant n-1$.
\end{prop}

Remarquons d'abord que si un groupe est simple et r\'{e}soluble,
alors son groupe d\'{e}riv\'{e}, \'{e}tant normal, est r\'{e}duit
\`{a} $\{1\}$; le groupe est donc ab\'{e}lien et, \'{e}tant simple,
est cyclique d'ordre premier. La proposition en
r\'{e}sulte.~\findem

\bigskip\exs
$(1)$ Les groupes $\SM_n$ sont r\'{e}solubles si et seulement si
$n\leqslant 4$.

$(2)$ Un groupe simple non ab\'{e}lien n'est pas r\'{e}soluble.

$(3)$ Soit $V$ un espace vectoriel de dimension $n$ sur un corps
commutatif $K$ et soit
$$V=V_0\supset V_1\supset \cdots \supset V_n=0$$ un drapeau
complet (i.e. une suite d\'{e}croissante de sous-espaces vectoriels
de $V$ tels que $\codim(V_i)=i$). On pose
$$G=\left\{s\in \mathbf{GL}(V)\ |\ sV_i=V_i,\ 0\leqslant i\leqslant n\right\}$$
(si on choisit dans $V$ une base adapt\'{e}e au drapeau, $G$ peut
\^{e}tre
identifi\'{e} au groupe des matrices triangulaires sup\'{e}rieures).\\
On d\'{e}finit alors une suite de sous-groupes $(B_i)_{0\leqslant i
\leqslant n}$ de $G$ par
$$B_i=\{s\in G\ |\ (s-1)V_j\subset V_{i+j},\ 0\leqslant j\leqslant n-i\}.$$
En particulier, $B_0=G$.

On va d\'{e}montrer que $(B_j,B_k)\subset B_{j+k}$ pour $0\leqslant
j \leqslant n$ et $0\leqslant k \leqslant n$ avec $0\leqslant j+k
\leqslant n$. Soient en effet $s\in B_j$, $t\in B_k$ et $x\in V_i$.
Il existe $v_{i+k}\in V_{i+k}$ tel que
$$tx=x+v_{i+k},$$ puis
$$stx=sx+sv_{i+k}=x+w_{i+j}+v_{i+k}+t_{i+j+k}$$
(avec $w_{i+j}\in V_{i+j}$ et $t_{i+j+k}\in V_{i+j+k}$). De m\^{e}me
$$tsx=t(x+w_{i+j})=x+v_{i+k}+w_{i+j}+t'_{i+j+k}$$ (avec
$t'_{i+j+k}\in V_{i+j+k}$). Donc
$$stx\equiv tsx \pmod{V_{i+j+k}}$$ ou encore
$$s^{-1}t^{-1}stx \equiv x \pmod{V_{i+j+k}}$$ d'o\`{u} le r\'{e}sultat.
En particulier:\\
$\bullet$ $(B_0,B_i)\subset B_i$ pour $0\leqslant i\leqslant n$,
donc les $B_i$ sont normaux dans $B_0=G$.\\
$\bullet$ $(B_i,B_i)=D(B_i)\subset B_{2i}\subset B_{i+1}$ pour
$1\leqslant i \leqslant n$, donc les quotients $B_i/B_{i+1}$ sont
ab\'{e}liens pour $1\leqslant i \leqslant n-1$.\\
$\bullet$ Enfin, $B_0/B_1=G/B_1$ s'identifie au groupe des matrices
diagonales (ab\'{e}lien car $K$ est commutatif). Donc la suite
$B_0=G\supset B_1\supset\cdots\supset B_n=\{1\}$ v\'{e}rifie la
condition $(2)$ et $G$ est r\'{e}soluble.

$(4)$ On verra ult\'{e}rieurement (th. \ref{thmburn5.4}) que tout
groupe d'ordre $p^aq^b$ (o\`{u} $p$ et $q$ sont premiers) est
r\'{e}soluble.

$(5)$ Mentionnons aussi le (tr\`{e}s difficile) th\'{e}or\`{e}me de
Feit-Thompson\footnote{R\'{e}f\'{e}rence: W. Feit et J.G. Thompson,
{\it Solvability of groups of odd order}, Pacific J. Math. $13$
($\mit 1963$), $775-1029$.}: tout groupe d'ordre impair est
r\'{e}soluble (ou encore: l'ordre d'un groupe simple non ab\'{e}lien
est pair).

$(6)$ Les groupes r\'{e}solubles interviennent en th\'{e}orie des
corps. Soit $K$ un corps de caract\'{e}ristique $0$ et soit
$\overline{K}$ une cl\^{o}ture alg\'{e}brique de $K$. On note
$K_{rad}$ le plus petit sous-corps de $\overline{K}$ contenant $K$
tel que pour tout $x\in K_{rad}$ et tout entier $n\geqslant 1$, on
ait $x^{1/n}\in K_{rad}$. On d\'{e}montre qu'une extension
galoisienne finie de $K$ est contenue dans $K_{rad}$ si et seulement
si son groupe de Galois est r\'{e}soluble (i.e. une \'{e}quation est
\emph{r\'{e}soluble par radicaux} si et seulement si son groupe de
Galois est r\'{e}soluble. C'est de l\`{a} que provient la
terminologie \og r\'{e}soluble\fg).

\section{Suite centrale descendante}
\label{suitecentrale} Soit $G$ un groupe. On appelle \emph{suite
centrale descendante} de $G$ la suite $(C^nG)_{n\geqslant 1}$ de
sous-groupes de $G$ d\'{e}finie par r\'{e}currence par:
\begin{eqnarray*}
C^1G & = & G,\\
C^{n+1}G & = & (G,C^nG)\quad \mbox{pour } n\geqslant 1.
\end{eqnarray*}
Pour tout $n\geqslant 1$, $C^nG$ est un sous-groupe
caract\'{e}ristique de $G$.

\begin{prop}
On a $(C^iG,C^jG)\subset C^{i+j}G$ pour tout $i\geqslant 1$ et
tout $j\geqslant 1$.
\end{prop}

On raisonne par r\'{e}currence sur $i$, la proposition \'{e}tant
claire pour $i=1$ et tout $j\geqslant 1$. Soit $j\geqslant 1$; on a
$(C^{i+1}G,C^jG)=\big((G,C^iG),C^jG\big)$. Or
$\big((C^iG,C^jG),G\big)\subset (C^{i+j}G,G)$ (par hypoth\`{e}se de
r\'{e}currence), donc $\big((C^iG,C^jG),G\big)\subset C^{i+j+1}G$.
De m\^{e}me $\big((C^jG,G),C^iG\big)$ est contenu dans $C^{i+j+1}G$.
Le lemme suivant permet de conclure.\findem

\begin{lemme}
Si $X$, $Y$ et $Z$ sont des sous-groupes normaux de $G$ et si $H$
est un sous-groupe de $G$ contenant $\big((Y,Z),X\big)$ et
$\big((Z,X),Y\big)$, alors $H$ contient $\big((X,Y),Z\big)$.
\end{lemme}

On utilise l'identit\'{e} de Hall:
$$\big(x^y,(y,z)\big)\big(y^z,(z,x)\big)\big(z^x,(x,y)\big)=1,$$
(o\`{u} $x^y=y^{-1}xy$) qui s'obtient en d\'{e}veloppant les
quarante-deux termes du membre de gauche. (cf. Bourbaki, A.I, \S\
6).~\findem

\bigskip\rmq L'identit\'{e} de Hall est l'analogue pour les groupes de
l'identit\'{e} de Jacobi:
$$\big[x,[y,z]\big]+\big[y,[z,x]\big]+\big[z,[x,y]\big]=0$$
pour les alg\`{e}bres de Lie. On peut l'utiliser pour associer \`{a}
tout groupe $G$ muni d'une filtration $(G_i)$ satisfaisant \`{a}
$(G_i,G_j)\subset G_{i+j}$ une alg\`{e}bre de Lie $gr(G)$, \`{a}
savoir
$$gr(G)=\bigoplus_i{G_i/G_{i+1}}.$$
Si $\xi\in G_i/G_{i+1}$ et $\eta\in G_j/G_{j+1}$, le crochet
$$[\xi,\eta]\in G_{i+j}/G_{i+j+1}$$ est par d\'{e}finition l'image du
commutateur $(x,y)$ o\`{u} $x$ (resp. $y$) est un repr\'{e}sentant
dans $G_i$ (resp. $G_j$) de $\xi$ (resp. $\eta$). Ceci s'applique
notamment au cas o\`{u} $G_i=C^iG$ (cf. aussi Bourbaki, Lie II, \S\
4, $\mbox{n}^{\rm o}$\! 4).

\section{Groupes nilpotents}\label{grpe nilp}
\label{classenilpotence}
\begin{defi}
Un groupe $G$ est dit \emph{nilpotent} s'il existe un entier $n$
positif tel que $C^{n+1}G=\{1\}$. La \emph{classe de nilpotence}
de $G$  est alors le plus petit tel entier $n$.
\end{defi}

En particulier:
\begin{itemize}
\item Le groupe $G$ est ab\'{e}lien si et seulement s'il est nilpotent
de classe $\leqslant 1$.

\item Un produit fini de groupes nilpotents est nilpotent et la
classe de nilpotence du produit est la borne sup\'{e}rieure des
classes des groupes.

\item Un sous-groupe (resp. un groupe quotient) d'un groupe
nilpotent est nilpotent.
\end{itemize}

\begin{prop}
Tout groupe nilpotent est r\'{e}soluble.
\end{prop}

En effet, pour tout $n\geqslant 0$, on a $D^nG\subset
C^{2^n}G$.~\findem

\bigskip
La r\'{e}ciproque est fausse: le groupe $\SM_3$ est r\'{e}soluble de
classe $2$. Regardons la suite centrale descendante de $\SM_3$; on a
$C^1\SM_3=\SM_3$, $C^2\SM_3={\cal C}_3$ (groupe cyclique d'ordre 3)
puis $C^3\SM_3={\cal C}_3$, etc. La suite est stationnaire et
n'atteint pas $\{1\}$. Donc $\SM_3$ n'est pas nilpotent.

\bigskip
On peut former des groupes nilpotents de la fa\c{c}on suivante:

\begin{prop}
Un groupe $G$ est nilpotent de classe $\leqslant n+1$ si et
seulement s'il est extension centrale d'un groupe $\Gamma$ nilpotent
de classe $\leqslant n$ (i.e. s'il existe une suite exacte
$\{1\}\rightarrow A \rightarrow G \rightarrow \Gamma \rightarrow
\{1\}$ o\`{u} $A$ est contenu dans le centre de $G$).
\end{prop}

Si $G$ est nilpotent de classe $n+1$, alors $C^{n+2}G=\{1\}$ donc
$C^{n+1}G$ est dans le centre de $G$. Posons
$\Gamma=G/{C^{n+1}G}$; on a $C^{n+1}\Gamma=\{1\}$ donc $\Gamma$
est nilpotent de classe $\leqslant n$.\\
R\'{e}ciproquement, si une telle suite exacte existe et si
$C^{n+1}\Gamma=\{1\}$, alors $C^{n+1}G\subset A$, donc $C^{n+1}G$
est contenu dans le centre de $G$, et $C^{n+2}G=\{1\}$.~\findem

\begin{coro}\label{3.8}
Soit $G$ un groupe nilpotent et soit $H$ un sous-groupe de $G$
distinct de $G$. Alors $N_G(H)$ est distinct de $H$.
\end{coro}

On raisonne par r\'{e}currence sur la classe de nilpotence $n$ de
$G$. Si $n=1$, on a $N_G(H)=G$ (car $G$ est ab\'{e}lien) donc
$N_G(H)\neq
H$.\\
Si $n\geqslant 2$, choisissons un sous-groupe central $A$ de $G$ tel
que $G/A$ soit nilpotent de classe $\leqslant n-1$. Alors $N_G(H)$
contient $A$. Si $H$ ne contient pas $A$ alors $N_G(H)\neq H$. Si
$H$ contient $A$, alors $H/A$ est un sous-groupe propre de $G/A$ et
l'hypoth\`{e}se de r\'{e}currence montre que $H/A\neq N_{G/A}(H/A)$.
Comme $N_G(H)/A=N_{G/A}(H/A)$, on en d\'{e}duit que $N_G(H)$ est
distinct de $H$.~\findem

\bigskip Une autre caract\'{e}risation des groupes nilpotents est
donn\'{e}e par la

\begin{prop}
Un groupe $G$ est nilpotent si et seulement s'il existe une
filtration $(G_i)_{1\leqslant i\leqslant n+1}$ telle que
$G_1=G\supset G_2\supset \cdots \supset G_{n+1}=\{1\}$ avec
$(G,G_i)\subset G_{i+1}$ pour tout $1 \leqslant i \leqslant n$.
\end{prop}

\rmq Une condition plus forte serait $(G_i,G_j)\subset G_{i+j}$.

\bigskip \dem Si une telle filtration existe, alors $C^kG\subset
G_k$ pour tout $k\geqslant 1$, donc $G$ est nilpotent.
R\'{e}ciproquement, si $G$ est nilpotent, on prend
$G_k=C^kG$.~\findem

\bigskip \ex Soit $V$ un espace vectoriel de dimension finie $n$
sur un corps $K$ et soit $$V=V_0\supset V_1\supset \cdots \supset
V_n=0$$ un drapeau complet de $V$. Reprenons l'exemple du \S\
\ref{grpe resol} et posons
$$B_j=\{g\in \mathbf{GL}(V)\ |\ (g-1)V_i\subset V_{i+j},\ i\geqslant 0\}.$$
Alors $B_1$ est nilpotent; en effet $(B_i)_{i\geqslant 1}$ est une
filtration telle que $B_1\supset B_2\supset \cdots \supset
B_n=\{1\}$ et $(B_i,B_j)\subset B_{i+j}$ comme on l'a vu plus
haut.

\bigskip En application, on a le:
\begin{thm}[Kolchin] Soit $V$ un espace vectoriel de
dimension finie sur un corps commutatif $K$ et soit $G$ un
sous-groupe de $\mathbf{GL}(V)$. On suppose que tout \'{e}l\'{e}ment
$g$ de $G$ admet $1$ comme unique valeur propre (i.e. $g-1$ est
nilpotent). Alors il existe un drapeau complet de $V$ tel que $G$
soit contenu dans le groupe $B_1$ correspondant (cf. ci-dessus). En
particulier, $G$ est nilpotent.
\end{thm}\label{Kolchin}

On raisonne par r\'{e}currence sur la dimension $n$ de $V$, le cas
$n=0$ \'{e}tant trivial. Supposons donc $n\geqslant 1$ et montrons
qu'il existe $x\in V$ non nul tel que $gx=x$ pour tout $g\in G$. Le
probl\`{e}me \'{e}tant lin\'{e}aire, on peut \'{e}tendre les
scalaires et supposer $K$ alg\'{e}briquement clos. Soit $A_G$ le
sous-espace vectoriel de $\End(V)$ engendr\'{e} par $G$. C'est une
sous-alg\`{e}bre de $\End(V)$. Distinguons deux cas:\\
$\bullet$ $V$ est un $A_G$-module r\'{e}ductible, i.e. il existe
$V'\subset V$, stable par $G$, distinct de $0$ et de $V$.
L'hypoth\`{e}se de r\'{e}currence s'applique \`{a} $V'$ et fournit
un $x$ non
nul dans $V'$ tel que $gx=x$ pour tout $g\in G$.\\
$\bullet$ Si $V$ est irr\'{e}ductible, on a $A_G=\End(V)$
d'apr\`{e}s un th\'{e}or\`{e}me de  Burnside (cf. Bourbaki, A. VIII,
\S\ 4, $\mbox{n}^{\rm o}$\!~3). Or si $a,a'\in A_G$, on a
$n\Tr(aa')=\Tr(a)\Tr(a')$. En effet, c'est clair si $a,a'\in G$ car
alors $\Tr(aa')=\Tr(a)=\Tr(a')=n$ et c'est donc vrai dans $A_G$ par
lin\'{e}arit\'{e}. Si $n>1$, les \'{e}l\'{e}ments $a$ et $a'$ de
matrices
respectives $$\left(%
\begin{array}{cc}
  1 & 0 \\
  0 & 0 \\
\end{array}%
\right)\;\;\; \mbox{et}\;\;\; \left(%
\begin{array}{cc}
  0 & 0 \\
  0 & 1 \\
\end{array}%
\right)$$ sont dans $\End(V)$ donc dans $A_G$. Or
$\Tr(a)=\Tr(a')=1$ et $\Tr(aa')=0$, donc $\dim(V)=1$ et tout $x$
non nul de $V$ convient.\\
Une fois l'existence de $x$ prouv\'{e}e, on note $V_{n-1}$ la droite
engendr\'{e}e par $x$. L'hypoth\`{e}se de r\'{e}currence,
appliqu\'{e}e \`{a} $V/V_{n-1}$ fournit un drapeau complet pour
$V/V_{n-1}$, stable par $G$, d'o\`{u} aussit\^{o}t un drapeau
complet pour $V$, stable par $G$ et il est clair que $G$ est contenu
dans le sous-groupe $B_1$ correspondant.~\findem

\section{Groupes nilpotents finis}\label{grpe nilp fini}
Soit $p$ un nombre premier.
\begin{prop}
Tout $p$-groupe est nilpotent.
\end{prop}

On va en donner deux d\'{e}monstrations.\\
$\bullet$ Soit $P$ un $p$-groupe; alors $P$ peut se plonger dans
$\mathbf{GL}_n(\ZM/p\ZM)$ pour un entier $n$ suffisamment grand.
Donc $P$ est contenu dans un $p$-Sylow de $\mathbf{GL}_n(\ZM/p\ZM)$,
qui est conjugu\'{e}, comme on l'a d\'{e}j\`{a} vu, \`{a} l'ensemble
$B_1$ des matrices triangulaires sup\'{e}rieures \`{a}
\'{e}l\'{e}ments diagonaux \'{e}gaux \`{a} $1$. D'apr\`{e}s
l'exemple du \S\ \ref{grpe nilp}, $B_1$ est nilpotent, donc aussi
$P$.

$\bullet$ On peut supposer $P\neq\{1\}$. Faisons op\'{e}rer $P$ sur
lui-m\^{e}me par automorphismes int\'{e}rieurs; l'ensemble des
points fixes est le centre $C(P)$ de $P$. Comme $P$ est un
$p$-groupe, on a d'apr\`{e}s le lemme \ref{lemme1}
$$|P|\equiv |C(P)| \pmod{p},$$
donc $C(P)\neq\{1\}$. Ainsi $P/C(P)$ est d'ordre strictement
inf\'{e}rieur \`{a} celui de $P$. Une r\'{e}currence permet de
conclure.~\findem

\begin{coro}\label{corogrpenilp}
Soit $G$ un groupe d'ordre $p^n$ (avec $p$ premier et $n\geqslant
1$). Alors:
\begin{enumerate}
\item[(1)] Tout sous-groupe de $G$ d'ordre $p^{n-1}$ est normal.

\item[(2)] Si $H$ est un sous-groupe de $G$, il existe une suite
de sous-groupes $(H_i)_{1\leqslant i\leqslant m}$ telle que
$H=H_1\subset H_2\subset \cdots \subset H_m=G$ avec
$(H_i:H_{i-1})=p$ pour $2\leqslant i\leqslant m$.

\item[(3)] Tout sous-groupe de $G$ distinct de $G$ est contenu
dans un sous-groupe d'ordre $p^{n-1}$.
\end{enumerate}
\end{coro}

$(1)$ Soit $H$ un sous-groupe de $G$ d'ordre $p^{n-1}$. Alors $H$
est distinct de $N_G(H)$ car $G$ est nilpotent. Donc l'ordre de
$N_G(H)$ est $p^n$ et $H$ est normal dans $G$.

$(2)$ On raisonne par r\'{e}currence sur l'ordre de $G$. Soit donc
$H$ un sous-groupe de $G$, qu'on peut supposer distinct de $G$. Soit
$H'$ un sous-groupe de $G$, distinct de $G$, contenant $H$, d'ordre
maximal. Alors $H'$ est distinct de $N_G(H')$, donc $N_G(H')=G$ et
$H'$ est normal dans $G$. En particulier, $G/H'$ est un
$p$-sous-groupe d'ordre $p^k$ pour un certain entier $k\geqslant 1$.
Si $k>1$, il existe un sous-groupe de $G/H'$, distinct de $\{1\}$ et
de $G/H'$, donc un sous-groupe de $G$, distinct de $G$ contenant
$H'$ strictement, ce qui est absurde. Donc $k=1$ et $(G:H')=p$.
L'hypoth\`{e}se de r\'{e}currence appliqu\'{e}e \`{a} $H'$ donne le
r\'{e}sultat.

$(3)$ d\'{e}coule imm\'{e}diatement de $(2)$.~\findem

\begin{coro}\label{coronilp}
Tout produit fini de $p$-groupes est nilpotent.
\end{coro}

On va d\'{e}montrer une r\'{e}ciproque.

\begin{thm}\label{caractnilp}
Soit $G$ un groupe fini. Les assertions suivantes sont
\'{e}quivalentes:
\begin{enumerate}
\item[(1)] $G$ est nilpotent.

\item[(2)] $G$ est produit de $p$-groupes.

\item[(3)] Pour tout $p$ premier, $G$ a un unique $p$-Sylow.

\item[(4)] Soient $p$ et $p'$ deux nombres premiers distincts et
soit $S_p$ (resp. $S_{p'}$) un $p$-Sylow de $G$ (resp. un
$p'$-Sylow); alors $S_p$ et $S_{p'}$ se centralisent mutuellement
(i.e. tout \'{e}l\'{e}ment de $S_p$ commute \`{a} tout
\'{e}l\'{e}ment de $S_{p'}$).

\item[(5)] Deux \'{e}l\'{e}ments de $G$ d'ordres premiers entre eux
commutent.
\end{enumerate}
\end{thm}

$(2)\Rightarrow(1)$ C'est le cor. \ref{coronilp} ci-dessus.

$(1)\Rightarrow(3)$ Soit $S$ un $p$-Sylow de $G$ et soit $N$ son
normalisateur. Alors $N$ est son propre normalisateur (cf. cor.
\ref{2.11}). Puisque nous supposons $G$ nilpotent, cela entra\^{i}ne
$N=G$, cf. cor. \ref{3.8}, donc $S$ est normal. Les $p$-Sylow de $G$
\'{e}tant conjugu\'{e}s, $G$ a un unique $p$-Sylow.

$(3)\Rightarrow(4)$ Pour tout nombre premier $p$, soit $S_p$
l'unique $p$-Sylow de $G$: il est normal dans $G$. Si $p$ et $p'$
sont deux premiers distincts, $S_p\cap S_{p'}$ est r\'{e}duit \`{a}
$\{1\}$ car c'est \`{a} la fois un $p$-groupe et un $p'$-groupe. Or
si $x\in S_p$ et $y\in S_{p'}$, on a $x^{-1}y^{-1}xy\in S_p\cap
S_{p'}$, d'o\`{u} $x^{-1}y^{-1}xy=1$. Donc $S_p$ et $S_{p'}$ se
centralisent mutuellement.

$(4)\Rightarrow(2)$ Pour tout $p$ premier, choisissons un $p$-Sylow
de $G$ qu'on note $S_p$. Le groupe engendr\'{e} par les $S_p$ ($p$
d\'{e}crivant l'ensemble des nombres premiers) est $G$ tout entier
(car son ordre est divisible par celui de $G$). D\'{e}finissons
alors une application
$$\varphi:
\left\{\!\!\begin{array}{rcl}
\prod_p{S_p} & \longrightarrow & G\\
{(s_p)}_p & \longmapsto & \prod_p{s_p}.
\end{array}\right.$$

Par hypoth\`{e}se, $\varphi$ est un homomorphisme, car les $S_p$ se
centralisent mutuellement. D'autre part, $\varphi$ est surjective,
car $G$ est engendr\'{e} par les $S_p$. Enfin $G$ et $\prod_p{S_p}$
ont m\^{e}me cardinal, donc $\varphi$ est un isomorphisme, ce qui
d\'{e}montre $(2)$.

$(2)\Rightarrow(5)$ Supposons que $G$ soit un produit de $p$-groupes
$G_p$. Soient $x$ et $y$ deux \'{e}l\'{e}ments de $G$ d'ordres
premiers entre eux. Alors $x={(x_p)}_p$ et $y={(y_p)}_p$ et pour
tout $p$ on a $x_p=1$ ou $y_p=1$. En effet, on a $x_p=1$ si et
seulement si $p$ ne divise pas l'ordre de $x$. Donc
$x^{-1}y^{-1}xy={(x_p^{-1}y_p^{-1}x_py_p)}_p=(1)$, donc $xy=yx$.

$(5)\Rightarrow(4)$ C'est \'{e}vident.

En r\'{e}sum\'{e}, nous avons montr\'{e} les implications suivantes:
$$
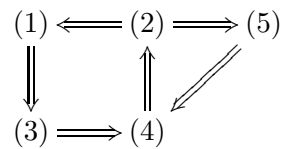
\xymatrix{(1) \ar@{=>}[d] & (2) \ar@{=>}[l] \ar@{=>}[r]& (5) \ar@{=>}[dl]\\
(3) \ar@{=>}[r] & (4) \ar@{=>}[u] & }$$
 Le th\'{e}or\`{e}me en r\'{e}sulte.~\findem

\section{Cas des groupes ab\'{e}liens}

Tout groupe ab\'{e}lien est nilpotent; d'apr\`{e}s le \S\ \ref{grpe
nilp fini}, tout groupe ab\'{e}lien est produit de $p$-groupes
ab\'{e}liens. On peut alors poursuivre la d\'{e}composition
gr\^{a}ce au

\begin{thm}
Tout $p$-groupe ab\'{e}lien est produit de groupes cycliques.
\end{thm}

\rmq Si $G$ est un $p$-groupe ab\'{e}lien fini, il existe un entier
$n$ tel que $p^nx=0$ pour tout $x\in G$, et l'on peut consid\'{e}rer
$G$ comme un module sur $\ZM/p^n\ZM$ (pour $n$ assez grand). Il
suffit donc de d\'{e}montrer le th\'{e}or\`{e}me suivant:

\begin{thm}
Tout module (non n\'{e}cessairement fini) sur $\ZM/p^n\ZM$ est somme
directe de modules monog\`{e}nes (isomorphes \`{a} $\ZM/p^i\ZM$ pour
$i\leqslant n$).
\end{thm}

Soit $G$ un module sur $\ZM/p^n\ZM$ et soit $V$ le module quotient
$G/pG$: c'est un espace vectoriel sur $\ZM/p\ZM$. Notons $G_i$
l'ensemble $\{x\in G\ |\ p^ix=0\}$. Les $G_i$ d\'{e}finissent une
filtration de $G$:
$$G_0=0\subset G_1\subset \cdots \subset G_n=G.$$
Soit $V_i$ l'image de $G_i$ dans $V$, alors
$$V_0=0\subset V_1\subset \cdots \subset V_n=V.$$
On peut choisir une base $S$ de $V$ adapt\'{e}e \`{a} cette
d\'{e}composition
(i.e. $S_i=S\cap V_i$ est une base de $V_i$).\\
Si $s\in S$, on note $i(s)$ le plus petit $i$ pour lequel $s\in S_i$
et on choisit un repr\'{e}sentant $\bar{s}$ de $s$ dans $G_{i(s)}$.
On a $p^{i(s)}\bar{s}=0$. Soit $G'=\bigoplus_{s\in
S}{\ZM/p^{i(s)}\ZM}$; on d\'{e}finit un homomorphisme $\varphi$ de
$G'$ dans $G$ de la fa\c{c}on suivante: si ${(n_s)}_{s\in S}\in G'$,
on pose $\varphi\big((n_s)\big)=\sum_{s\in S}{n_s\bar{s}}$. On va
montrer que $\varphi$ est un isomorphisme.\\
$\bullet$ $\varphi$ est surjectif: il suffit de prouver que le
sous-groupe $H$ de $G$ engendr\'{e} par les $\bar{s}$ est $G$ tout
entier. Or l'application projection de $H$ dans $V$ est surjective
(par d\'{e}finition, $\bar{s}$ s'envoie sur $s$ qui parcourt une
base de $V$). Donc $H+pG=G$ ou encore $G=pG+H$. En it\'{e}rant,
$$G=p(pG+H)=p^2G+H=\cdots=p^nG+H=H.$$
$\bullet$ $\varphi$ est injectif: soit $(n_s)\in G'$ tel que
$\sum_{s\in S}{n_s\bar{s}}=0$ dans $G$. Montrons que $n_s$ est
divisible par $p^i$ pour tout $i$ (donc $n_s=0$). En effet, montrons
par r\'{e}currence sur $i$ que $n_s\in
p^i\left(\ZM/p^{i(s)}\ZM\right)$. Pour $i=0$, il n'y a rien \`{a}
d\'{e}montrer. Si c'est vrai jusqu'\`{a} l'ordre $k$, on a $n_s\in
p^k\left(\ZM/p^{i(s)}\ZM\right)$. En particulier $n_s=0$ si
$i(s)\leqslant k$. On regarde donc les $s$ pour lesquels
$i(s)\geqslant k+1$. Par hypoth\`{e}se,
$$\sum_{i(s)\geqslant k+1}\!\!\!{n_s\bar{s}}=0.$$
Posons $n_s=p^km_s$ avec $m_s\in \ZM/p^{i(s)}\ZM$; on a
$$p^k\!\!\!\sum_{i(s)\geqslant k+1}\!\!\!{m_s\bar{s}}=0,$$
c'est-\`{a}-dire
$$\sum_{i(s)\geqslant k+1}\!\!\!{m_s\bar{s}}\ \in\,G_k,$$
d'o\`{u}, par projection, $$\sum_{i(s)\geqslant k+1}\!\!\!{m_s s}\
\in\,V_k.$$ Vu le choix de $S$, on en d\'{e}duit que $m_s\equiv 0
\pmod{p}$. Donc $p^{k+1}$ divise $n_s$ et on en d\'{e}duit que $n_s$
appartient \`{a} $p^{k+1}\!\left(\ZM/p^{i(s)}\ZM\right)$.~\findem

\bigskip\exo Si $G$ est un module sur $\ZM/p^n\ZM$, tout
sous-module de $G$ isomorphe \`{a} $\ZM/p^n\ZM$ est facteur direct.

\bigskip \noindent{\it Application des $2$-groupes: constructions par r\`{e}gle et compas.}
Soit $K$ un corps de caract\'{e}ristique quelconque et soit $L$ une
extension galoisienne de $K$ dont le groupe de Galois $G$ est un
$2$-groupe. D'apr\`{e}s le cor. \ref{corogrpenilp}, il existe $G'$
normal dans $G$, d'indice $2$, donc un corps interm\'{e}diaire $L'$,
fix\'{e} par $G'$, qui est une extension quadratique de $K$. Si on
it\`{e}re, $L/K$ appara\^{i}t comme une tour d'extensions
quadratiques.\\
R\'{e}ciproquement, si $L/K$ est une telle tour, l'extension
galoisienne engendr\'{e}e a pour groupe de Galois un $2$-groupe. On
a ainsi une caract\'{e}risation des extensions galoisiennes dont le
groupe est un $2$-groupe. Si la caract\'{e}ristique de $K$ est
distincte de $2$, une extension quadratique est de la forme
$K\big[\sqrt{a}\, \big]\simeq K[X]/(X^2-a)$ o\`{u} $a\in
K^*\mathbf{-} {K^*}^2$. (Si $\car
K=2$ on remplace $X^2-a$ par $X^2+X+a$).\\
Ceci rejoint le probl\`{e}me des nombres constructibles par
r\`{e}gle et compas: ce sont les nombres (alg\'{e}briques sur $\QM$)
contenus dans une extension galoisienne de $\QM$ dont le groupe de
Galois est un $2$-groupe.

\bigskip\ex (Impossibilit\'e de la duplication du cube) Le nombre$\sqrt[3]{2}$ n'est pas constructible par r\`{e}gle et compas, car $X^3-2$ est irr\'{e}ductible, et son degr\'e n'est pas une puissance de $2$.

\section{Sous-groupe de Frattini}
\label{sgFrattini} Soit $G$ un groupe fini. On appelle
\emph{sous-groupe de Frattini} de $G$ et on note $\Phi(G)$
l'intersection des sous-groupes maximaux\footnote{Un sous-groupe $H$
de $G$ est dit \emph{maximal} s'il est distinct de $G$ et maximal
pour cette propri\'{e}t\'{e}; on dit alors que l'action de $G$ sur
$G/H$ est \emph{primitive}.} de $G$. C'est un sous-groupe
caract\'{e}ristique de~$G$.

Un probl\`{e}me int\'{e}ressant est de savoir \`{a} quelles
conditions une partie $S$ de $G$ engendre le groupe $G$. On a la
proposition suivante:

\begin{prop}
Soit $S$ une partie de $G$ et soit $H$ le sous-groupe engendr\'{e}
par $S$. On a $H=G$ si et seulement si $H.\Phi(G)=G$, i.e. si $S$
engendre $G/\Phi(G)$.
\end{prop}

En effet, si $H.\Phi(G)=G$ et $H\neq G$, il existe $H'$ maximal
contenant $H$ et distinct de $G$; $\Phi(G)$ est aussi contenu dans
$H'$ par d\'{e}finition: $G=H.\Phi(G)$ est donc contenu dans $H'$,
ce qui est absurde.~\findem

\begin{thm}
Le groupe $\Phi(G)$ est nilpotent.
\end{thm}

On va utiliser la caract\'{e}risation $(3)$ du th. \ref{caractnilp}.
Soit $S$ un $p$-Sylow de $\Phi(G)$ (pour un $p$ premier quelconque).
On a vu dans l'\'{e}tude des groupes de Sylow (prop. \ref{prop2.10})
que $G=\Phi(G).N_G(S)$. D'apr\`{e}s la proposition ci-dessus, on a
$N_G(S)=G$, donc $S$ est normal dans $G$, donc dans $\Phi(G)$, donc
est l'unique $p$-Sylow de $\Phi(G)$ qui est alors nilpotent.~\findem

\bigskip Dans le cas o\`{u} $G$ est un $p$-groupe, on a une
caract\'{e}risation simple de $\Phi(G)$:

\begin{thm}
Si $G$ est un $p$-groupe, $\Phi(G)$ est le sous-groupe engendr\'{e}
par les commutateurs et les puissances $p$-i\`{e}mes de $G$, i.e.
$\Phi(G)=(G,G).G^p$.
\end{thm}

Si $G$ est d'ordre $p^n$, ses sous-groupes $H$ maximaux sont d'ordre
$p^{n-1}$, donc sont normaux, donc sont des noyaux d'homomorphismes
surjectifs de $G$ dans $G/H\simeq \ZM/p\ZM$ (et r\'{e}ciproquement
un tel homomorphisme d\'{e}finit un sous-groupe maximal de $G$).
Donc $\Phi(G)$ est l'intersection des noyaux d'homomorphismes
surjectifs de $G$ dans $\ZM/p\ZM$. Or un tel
homomorphisme est trivial sur $(G,G)$ et $G^p$. Donc $(G,G).G^p\subset\Phi(G)$.\\
R\'{e}ciproquement $V=G/(G,G).G^p$ est un espace vectoriel sur
$\ZM/p\ZM$, dans lequel $0$ est intersection d'hyperplans. Or un
hyperplan est le noyau d'un homomorphisme de $V$ dans $\ZM/p\ZM$. Le
groupe $(G,G).G^p$ contient donc $\Phi(G)$. Finalement
$\Phi(G)=(G,G).G^p$.~\findem

\bigskip\App $G/\Phi(G)$ est donc le plus grand quotient de $G$ qui soit
ab\'{e}lien \'{e}l\'{e}mentaire (i.e. produit de groupes cycliques
d'ordre $p$).\label{elementaire}

\begin{coro}
Une partie $S$ de $G$ engendre $G$ si et seulement si son image
dans $G/(G,G).G^p$ engendre ce groupe.
\end{coro}

En effet, on a alors $\langle S\rangle .\Phi(G)=G$ donc $\langle S
\rangle=G$.~\findem

\bigskip Ainsi le cardinal minimum d'une partie $S$ g\'{e}n\'{e}ratrice de
$G$ est $\dim_{\FM_p}\big(G/\Phi(G)\big)$ (dans le cas o\`{u} $G$
est un $p$-groupe).

\bigskip \noindent{\it Caract\'{e}risations par les sous-groupes \`{a} deux
g\'{e}n\'{e}rateurs.} On peut dans la m\^{e}me veine \'{e}tudier si
les propri\'{e}t\'{e}s de certains sous-groupes de $G$ permettent de
d\'{e}montrer des propri\'{e}t\'{e}s analogues pour le groupe tout
entier.

\begin{prop}
Soit $G$ un groupe (resp. un groupe fini). Supposons que tout
sous-groupe de $G$ engendr\'{e} par deux \'{e}l\'{e}ments soit
commutatif (resp. nilpotent). Alors $G$ est commutatif (resp.
nilpotent).
\end{prop}

Soient $x,y\in G$. Le groupe $\langle x,y\rangle$ est commutatif,
donc $xy=yx$. Pour les groupes nilpotents finis, on utilise la
caract\'{e}risation $(5)$ du th. \ref{caractnilp}.~\findem

\bigskip On peut se demander si un th\'{e}or\`{e}me analogue est vrai pour
les groupes r\'{e}solubles. On va d'abord d\'{e}finir la notion de
\emph{groupe simple minimal}. Soit $G$ un groupe fini simple non
ab\'{e}lien. On dit que $G$ est \emph{minimal}\label{minimal} si
tout sous-groupe de $G$ distinct de $G$ est r\'{e}soluble.

\begin{lemme}\label{3.22}
Si $G$ n'est pas r\'{e}soluble, il existe un sous-groupe $H$ de $G$
et un sous-groupe normal $K$ de $H$, tels que $H/K$ soit simple
minimal.
\end{lemme}

\ex $G=\AM_6$ est simple (non minimal); on peut prendre $H=\AM_5$
et $K=\{1\}$.

\bigskip\dem Soit $H$ un sous-groupe non r\'{e}soluble minimal de $G$ et
soit $K$ un sous-groupe normal de $H$, distinct de $H$, et maximal.
Le groupe $K$ est r\'{e}soluble (car strictement contenu dans $H$).
Le quotient $H/K$ est simple (car $K$ est maximal) et non
ab\'{e}lien (sinon $H$ serait r\'{e}soluble); il est de plus minimal
(tout sous-groupe s'obtient comme quotient par $K$ d'un sous-groupe
$H'$ contenant $K$ et contenu dans $H$, donc est r\'{e}soluble).

\bigskip On a alors une r\'{e}ponse partielle \`{a} notre question:

\begin{prop}
Les deux assertions suivantes sont \'{e}quivalentes:
\begin{enumerate}
\item[(1)] Tout groupe simple minimal peut \^{e}tre engendr\'{e} par deux
\'{e}l\'{e}ments.

\item[(2)] Tout groupe tel que ses sous-groupes engendr\'{e}s par deux
\'{e}l\'{e}ments soient r\'{e}solubles est r\'{e}soluble.
\end{enumerate}
\end{prop}

Supposons $(2)$. Soit $G$ simple minimal; si $\langle x,y\rangle\neq
G$ pour tout couple $(x,y)\in G\times G$, alors $\langle x,y\rangle$
est r\'{e}soluble comme sous-groupe de $G$ strict et d'apr\`{e}s
$(2)$, $G$ est aussi r\'{e}soluble, ce qui est impossible.

Supposons $(1)$. Si $G$ n'est pas r\'{e}soluble, il existe $H$ et
$K$ comme dans le lemme \ref{3.22}. Le groupe $H/K$ est simple
minimal donc engendr\'{e} par deux \'{e}l\'{e}ments $\bar{x}$ et
$\bar{y}$. Soit $H'$ le sous-groupe de $H$ engendr\'{e} par $x$ et
$y$ (repr\'{e}sentants respectifs de $\bar{x}$ et $\bar{y}$). C'est
un groupe r\'{e}soluble par hypoth\`{e}se; or $H/K$ est l'image par
projection de $H'$ donc est aussi r\'{e}soluble, ce qui est
absurde.~\findem

\bigskip On est donc ramen\'{e} au probl\`{e}me suivant: trouver la liste
de tous les groupes simples minimaux et chercher s'ils sont
engendr\'{e}s par deux \'{e}l\'{e}ments. Ce probl\`{e}me a
\'{e}t\'{e} r\'{e}solu par Thompson\label{thmThompson}, qui a
montr\'{e}\footnote{R\'{e}f\'{e}rences: J.G. Thompson, {\it Non
solvable finite groups all of whose local subgroups are solvable, I,
II, $\dots$, VI}, Bull. A.M.S. $74$ ($\mit 1968$), $383-437$; Pac.
J. Math. $33$ ($\mit 1970$),
$451-536$; $\dots$; Pac. J. Math. $51$ ($\mit 1974$), $573-630$.\\
\indent\ \ Il existe une d\'{e}monstration ind\'{e}pendante du
th\'{e}or\`{e}me de classification de Thompson: P. Flavell, {\it
Finite groups in which two elements generate a solvable group},
Invent. math. $121$ $({\mit 1995})$, $279-285$.} que tout groupe
simple minimal est isomorphe \`{a} l'un des suivants (que l'on peut
engendrer par deux \'{e}l\'{e}ments):
\begin{itemize}
 \item $\mathbf{PSL}_2(\FM_p)$, $p\geqslant 5$, $p\not\equiv \pm 1
\pmod{5}$,
 \item $\mathbf{PSL}_2(\FM_{2^p})$, $p$ premier $\geqslant 3$,
 \item $\mathbf{PSL}_2(\FM_{3^p})$, $p$ premier $\geqslant 3$,
 \item $\mathbf{PSL}_3(\FM_3)$,
 \item groupes de Suzuki $Sz(2^p)$, $p$ premier $\geqslant 3$.
\end{itemize}

\bigskip Signalons un probl\`{e}me: est-il vrai qu'un groupe simple non
ab\'{e}lien qui est \emph{minimal} au sens na\"{i}f (i.e. qui ne
contient pas de sous-groupe propre non ab\'{e}lien qui soit simple)
est minimal au sens d\'{e}fini plus haut?\footnote{R\'{e}ponse: oui,
d'apr\`{e}s la classification des groupes simples.}

\chapter{Cohomologie et extensions}\label{chap4}

\section{D\'{e}finitions}

Soient $G$ un groupe (not\'{e} multiplicativement) et $A$ un
$G$-module (autrement dit un groupe ab\'{e}lien not\'{e}
additivement sur lequel $G$ op\`{e}re par automorphismes). On note
$sa$ le transform\'{e} de l'\'{e}l\'{e}ment $a\in A$ par
l'\'{e}l\'{e}ment $s\in G$. On a
$$\begin{array}{rcl}
(st)a & = & s(ta),\\
1a & = & a,\\
s(a_1+a_2) & =& sa_1+sa_2,
\end{array}$$
si $s,t\in G$ et $a,a_1,a_2\in A$.

Des exemples d'une telle situation sont donn\'{e}s par:\\
$(1)$ Action triviale d'un groupe $G$ sur un groupe ab\'{e}lien $A$:
$sa=a$ pour $s\in G$ et $a\in A$.

$(2)$ Si $L$ est une extension galoisienne du corps $K$, de groupe
de Galois $G$, alors $G$ op\`{e}re par automorphismes sur $L$ muni
de l'addition ou $L^*$ muni de la multiplication.

\begin{defi}
Soit $n$ un entier positif ou nul. On appelle
\emph{$n$-cocha\^{i}ne}, ou \emph{cocha\^{i}ne de degr\'{e} $n$} sur
$G$ \`{a} valeurs dans $A$ toute fonction de $n$ variables de $G$
\`{a} valeurs dans $A$:
$$f:\left\{\!\!
\begin{array}{rcl}
G\times G\times \cdots \times G & \longrightarrow & A\\
(s_1,s_2,\dots,s_n) & \longmapsto & f(s_1,s_2,\dots,s_n).
\end{array}
\right.$$
\end{defi}

L'ensemble des cocha\^{i}nes, muni de l'addition induite par celle
de $A$, forme un groupe ab\'{e}lien not\'{e} $C^n(G,A)$.

\bigskip\exs $n=0$: par convention, une fonction de $0$ variable \`{a} valeurs
dans $A$ est un \'{e}l\'{e}ment de $A$. D'o\`{u} $C^0(G,A)=A$. On
note $f_a$ l'\'{e}l\'{e}ment de $C^0(G,A)$ correspondant \`{a}
l'\'{e}l\'{e}ment $a\in A$.

$n=1$: $C^1(G,A)=\{f:G\rightarrow A\}$.

$n=2$: $C^2(G,A)=\{f:G\times G\rightarrow A\}$.

\begin{defi}
Si $f\in C^n(G,A)$, on appelle \emph{cobord} de $f$ et on note $df$
l'\'{e}l\'{e}ment de $C^{n+1}(G,A)$ d\'{e}fini par la formule:
$$\begin{array}{rcl}
df(s_1,\dots,s_n,s_{n+1}) & = &
s_1f(s_2,\dots,s_{n+1})+\displaystyle\sum_{i=1}^{n}{(-1)^if(s_1,\dots,s_{i-1},s_i s_{i+1},\dots)}\\
{} & {} & \hfill{{} +(-1)^{n+1}f(s_1,\dots,s_n).}\end{array}$$
\end{defi}

Regardons ce qu'est $d$ pour les petites valeurs de $n$.\\
$\bullet$ $d:\ C^0(G,A)\longrightarrow C^1(G,A)$. Soit $a\in A$;
on cherche $df_a$. On a $df_a(s)=sa-a$. Noter que $df_a=0$
si et seulement si $a$ est fix\'{e} par $G$.\\
$\bullet$ $d:\ C^1(G,A)\longrightarrow C^2(G,A)$. Soit $f$ une
$1$-cocha\^{i}ne; on a $$df(s,t)=sf(t)-f(st)+f(s).$$ $\bullet$ $d:\
C^2(G,A)\longrightarrow C^3(G,A)$. Soit $f$ une $2$-cocha\^{i}ne; on
a
$$df(u,v,w)=uf(v,w)-f(uv,w)+f(u,vw)-f(u,v).$$

\begin{thm}[Formule fondamentale]\label{fondam}
On a $d\circ d=0$. Autrement dit, le compos\'{e}
$$\xymatrix{\relax C^n(G,A) \ar[r]^-{d} & C^{n+1}(G,A) \ar[r]^d & C^{n+2}(G,A)}$$ est nul.
\end{thm}

Nous allons faire la v\'{e}rification seulement dans les cas $n=0$
et $n=1$, laissant en exercice la v\'{e}rification g\'{e}n\'{e}rale.

Soit $a\in A$. On a $df_a(s)=sa-a$ d'o\`{u}
\begin{eqnarray*}
d\circ d\,(f_a)(s,t) & = & sdf_a(t)-df_a(st)+df_a(s)\\
{} & = & s(ta-a)-(sta-a)+(sa-a)\\
{} & = & 0.
\end{eqnarray*}

Regardons maintenant $f\in C^1(G,A)$: $$\begin{array}{rcl}
d\circ d\,(f)(u,v,w) & = & u\,df(v,w)-df(uv,w)+df(u,vw)-df(u,v)\\
{} & = & u\big(vf(w)-f(vw)+f(v)\big)-\big(uvf(w)-f(uvw)+f(uv)\big)\\
{} & {} & \hspace{1.4cm} {} +\big(uf(vw)-f(uvw)+f(u)\big)-\big(uf(v)-f(uv)+f(u)\big)\\
{} & = & 0.
\end{array}\eqno\square$$

\begin{defi}
Une $n$-cocha\^{i}ne $f$ est dite un \emph{$n$-cocycle} si $df=0$.
Elle est dite un \emph{$n$-~cobord} s'il existe une
$(n-1)$-cocha\^{i}ne $g$ telle que $f=dg$.
\end{defi}
\label{cobord} D'apr\`{e}s le th. \ref{fondam}, tout $n$-cobord est
un $n$-cocycle. On note $Z^n(G,A)$ le groupe des $n$-cocycles et
$B^n(G,A)$ le
groupe des $n$-cobords.\\
On note $H^n(G,A)$ le groupe quotient $Z^n(G,A)/B^n(G,A)$ et on
l'appelle le \emph{$n$-i\`{e}me groupe de cohomologie de $G$ \`{a}
valeurs dans $A$}.

\bigskip\exs $(1)$ Pour $n=0$, on convient que $B^0=\{0\}$. En
notant $A^G$ le groupe des \'{e}l\'{e}ments de $A$ fix\'{e}s par
$G$, on a vu que
$$df_a=0\Longleftrightarrow a\in A^G,$$
et donc $H^0(G,A)=A^G$.

$(2)$ Pour $n=1$, un \'{e}l\'{e}ment de $Z^1(G,A)$ est une
application $f$ de $G$ dans $A$ telle que $df(s,t)=0$ pour tous
$s,t\in G$, ce qui donne $$f(st)=sf(t)+f(s).$$ On dit que $f$ est un
\emph{homomorphisme crois\'{e}}. Si l'action de $G$ sur $A$ est
triviale, on a $sf(t)=f(t)$ et $f$ est un homomorphisme de $G$ dans
$A$; comme $B^1(G,A)=\{0\}$, on a alors
$$H^1(G,A)=\Hom(G,A),$$
o\`{u} $\Hom(G,A)$ est le groupe des homomorphismes de groupes de
$G$ dans $A$.

$(3)$ Pour $n=2$, une $2$-cocha\^{i}ne $f$ est un $2$-cocycle si
$$uf(v,w)-f(uv,w)+f(u,vw)-f(u,v)=0$$
pour tous $u,v,w\in G$. Une telle cocha\^{i}ne s'appelle aussi un
\emph{syst\`{e}me de facteurs}.

\section{Extensions}\label{extensions}

\begin{defi}
Soient $A$ et $G$ deux groupes. On dit que $E$ est une
\emph{extension} de $G$ par $A$ si l'on a une suite exacte
$$\xymatrix{\{1\} \ar[r] & A \ar[r] & E \ar[r] & G \ar[r] & \{1\}}$$
avec $A$ normal dans $E$.
\end{defi}
\rmq Dans ce \S, on suppose $A$ commutatif.

\bigskip Toute extension $E$ de $G$ par $A$ d\'{e}finit une action de
$G$ sur $A$ de la mani\`{e}re suivante: remarquons d'abord que $E$
agit sur $A$ par automorphismes int\'{e}rieurs (puisque $A$ est
normal dans $E$); on a un homomorphisme
$$\left\{\!\!\begin{array}{rcl}
E & \longrightarrow & \Aut(A)\\
e & \longmapsto & \Int(e)_{|A},\\
\end{array}\right.$$
qui passe au quotient $G$: en effet, si $s\in G$, on choisit $e\in
E$ qui rel\`{e}ve $s$; alors $\Int(e)$ ne d\'{e}pend pas du choix du
rel\`{e}vement de $s$; changer $e$ en $e'$ au dessus de $s$ revient
en effet \`{a} le multiplier par un \'{e}l\'{e}ment $a$ de $A$, or
$a$ agit trivialement sur $A$ par automorphismes int\'{e}rieurs
puisque $A$ est ab\'{e}lien. Donc $G$ agit sur $A$: $$\xymatrix{ E
\ar[rr] \ar[rd] && \Aut(A) \\ & G \ar[ur]}$$

On va donc consid\'{e}rer $A$ comme un $G$-module; les lois de
groupe \'{e}tant \'{e}crites multiplicativement, la loi d'action de
$G$ sur $A$ sera \'{e}crite ${}^s\!a$ pour $a\in A$ et $s\in G$. On
va associer \`{a} toute extension de $G$ par $A$ une classe de
cohomologie de $H^2(G,A)$ qui d\'{e}termine cette extension \`{a}
isomorphisme pr\`{e}s. Et l'on verra que tout \'{e}l\'{e}ment de
$H^2(G,A)$ peut \^{e}tre obtenu ainsi, cf. th. \ref{th4.3}.

Soit $E$ une extension de $G$ par $A$; on a une surjection $\pi$
de $E$ sur $G$.

\begin{defi}
Une \emph{section} $h$ de $\pi$ est une application de $G$ dans
$E$ telle que $\pi\circ h=\Id_G$.
$$\xymatrix{E \ar[d]_\pi \\ G \ar@/_0.5cm/[u]_h}$$
\end{defi}\label{section}

Au-dessus de $s\in G$, on choisit un point dans la fibre
$\pi^{-1}(s)$. Tout \'{e}l\'{e}ment $e\in E$ s'\'{e}crit alors de
mani\`{e}re unique $ah(x)$, avec $a\in A$ et $x\in G$ (en fait
$x=\pi(e)$).

Cherchons \`{a} mettre sous la forme $ch(z)$ l'\'{e}l\'{e}ment
$ah(x)bh(y)$. On a
$$ah(x)bh(y)=ah(x)bh(x)^{-1}h(x)h(y).$$
L'action de $x\in G$ sur $A$ est donn\'{e}e par l'action de
l'automorphisme int\'{e}rieur d'un \'{e}l\'{e}ment de $E$ au-dessus
de $x$, par exemple $h(x)$. Donc $h(x)bh(x)^{-1}={}^xb$ (qui est
dans $A$, puisque $A$ est normal). Posons
$$h(x)h(y)=f_h(x,y)h(xy).$$ On a
$f_h(x,y)\in A$ puisque $h(x)h(y)$ et $h(xy)$ ont m\^{e}me image
dans $G$ par $\pi$. On a finalement obtenu:
$$ah(x)bh(y)=a\,{}^xbf_h(x,y)h(xy)$$ avec $a\,{}^xbf_h(x,y)\in A$.

\bigskip Nous allons maintenant voir comment $f_h$ varie avec $h$.
Soient donc $h$ et $h'$ deux sections de $\pi$ ($h,h':G\rightarrow
E$). Alors $h(s)$ et $h'(s)$ diff\`{e}rent par un \'{e}l\'{e}ment de
$A$. Posons $h'(s)=l(s)h(s)$; l'application $l$ est une
$1$-cocha\^{i}ne de $G$ \`{a} valeurs dans $A$. Calculons $f_{h'}$
\`{a} l'aide de $l$ et de $f_h$. On a
$$h'(s)h'(t)=f_{h'}(s,t)h'(st)=f_{h'}(s,t)l(st)h(st),$$
mais
\begin{eqnarray*}
h'(s)h'(t) & = & l(s)h(s)l(t)h(t)\\
{} & = & l(s)h(s)l(t)h(s)^{-1}h(s)h(t)\\
{} & = & l(s)\,{}^sl(t)f_h(s,t)h(st),
\end{eqnarray*}
d'o\`{u} l'on tire
\begin{eqnarray*}
f_{h'}(s,t) & = & l(s)\,{}^sl(t)f_h(s,t)l(st)^{-1}\\
{} & = & f_h(s,t)\,{}^sl(t)l(s)l(st)^{-1},
\end{eqnarray*}
car $A$ est commutatif. Or, en notation multiplicative, on a
$$dl(s,t)={}^sl(t)l(s)l(st)^{-1}.$$
D'o\`{u}
$$f_{h'}=f_h\,dl.$$
Donc, quand $h$ varie, $f_h$ ne change que par multiplication par un
cobord. On peut donc associer \`{a} $E$ la classe de cohomologie de
$f_h$ dans $H^2(G,A)$; appelons $e$ cette classe. Quand trouve-t-on
$e=0$? Cela signifie (en notation multiplicative) qu'il existe une
section $h$ telle que $f_h(s,t)=1$ pour tous $s,t\in G$, i.e. que
$h$ est un homomorphisme.

\begin{defi}
Une extension $E$ de $G$ par $A$ est dite \emph{triviale} s'il
existe un homomorphisme $h:G\rightarrow E$ telle que $\pi\circ
h=\Id_G$ (ou, de fa\c{c}on \'{e}quivalente, si $e=0$).
\end{defi}

Examinons une telle extension: tout \'{e}l\'{e}ment de $E$
s'\'{e}crit $ah(s)$ de mani\`{e}re unique et
$ah(s)bh(t)=a\,{}^sbh(st)$. Donc on conna\^{i}t $E$ d\`{e}s qu'on
conna\^{i}t $A$, $G$ et l'action de $G$ sur $A$. Le groupe $E$ est
isomorphe au groupe des couples $(a,s)$ avec $a\in A$ et $s\in G$,
muni de la loi
$$(a,s)(b,t)=(a\,{}^sb,st).$$

On appelle un tel $E$ un \emph{produit
semi-direct}\label{semidirect1} de $G$ par $A$. On vient de voir: la
classe nulle de $H^2(G,A)$ correspond \`{a} l'extension triviale de
$G$ par $A$, qui est le produit semi-direct de $G$ par $A$
d\'{e}fini par l'action de $G$ sur $A$.

\begin{thm}
L'application $f_h$ est un $2$-cocycle de $G$ \`{a} valeurs dans
$A$.
\end{thm}

Il faut v\'{e}rifier que $f_h$ appartient au noyau de $d$,
l'homomorphisme de cobord. L'\'{e}criture est ici multiplicative; il
faut donc voir que
$$df_h(u,v,w)=1$$ pour tous $u,v,w\in G$; or $df_h$ s'\'{e}crit
$$df_h(u,v,w)={}^u\!f_h(v,w)f_h(u,vw)f_h(uv,w)^{-1}f_h(u,v)^{-1}.$$
Nous allons \'{e}crire $h(u)h(v)h(w)$ sous la forme $ah(uvw)$ avec
$a\in A$ de deux mani\`{e}res diff\'{e}rentes en utilisant
l'associativit\'{e}
de la loi de groupe dans $E$.\\
On a $$\big(h(u)h(v)\big)h(w)=f_h(u,v)f_h(uv,w)h(uvw)$$ et
$$h(u)\big(h(v)h(w)\big)={}^u\!f_h(v,w)f_h(u,vw)h(uvw)$$
d'o\`{u}
$${}^u\!f_h(v,w)f_h(u,vw)=f_h(u,v)f_h(uv,w)$$
ce qui est bien
$$df_h(u,v,w)=1.\eqno\square$$

Nous allons enfin voir:

\begin{thm}\label{th4.3}
Toute classe de cohomologie de $H^2(G,A)$ correspond \`{a} une
extension de $G$ par $A$.
\end{thm}

On va reconstruire la situation pr\'{e}c\'{e}dente: soit $f\in
Z^2(G,A)$. D\'{e}finissons $E$ ensemblistement par $E=A\times G$. On
d\'{e}finit la loi de $E$ par
$$(a,s)(b,t)=\big(a\,{}^sbf(s,t),st\big).$$
Tout d'abord $E$ est un groupe:\\
$\bullet$ La loi est associative: le calcul fait ci-dessus pour voir
que $f_h$ est un $2$-cocycle \`{a} partir de l'associativit\'{e} de
la loi de $E$ se reprend \`{a} l'envers.\\
$\bullet$ Si $\varepsilon=f(1,1)^{-1}$, alors l'\'{e}l\'{e}ment
$(\varepsilon,1)$ est \'{e}l\'{e}ment neutre. En effet,
$$(a,s)(\varepsilon,1)=\big(a{}^s\!\varepsilon f(s,1),s\big)$$
or $f$ est un $2$-cocycle donc $df=1$ et
$$df(s,1,1)={}^s\!f(1,1)f(s,1)^{-1}f(s,1)f(s,1)$$
donc
$$1=df(s,1,1)={}^s\!\varepsilon^{-1}f(s,1)^{-1}$$
et $(\varepsilon,1)$ est bien \'{e}l\'{e}ment neutre.\\
$\bullet$ On fait de m\^{e}me le calcul de l'inverse.\\
On a un homomorphisme surjectif \'{e}vident de $E$ dans $G$:
$$\left\{\!\!
\begin{array}{rcl}
E & \longrightarrow & G\\
(a,s) & \longmapsto & s
\end{array}\right.$$
et l'application
$$\left\{\!\!
\begin{array}{rcl}
A & \longrightarrow & E\\
a & \longmapsto & (a\varepsilon ,1)
\end{array}\right.$$
est un homorphisme (car $A$ est ab\'{e}lien) \'{e}videmment
injectif.

Finalement on a bien:
$$\xymatrix{\{1\} \ar[r] & A \ar[r] & E \ar[r] & G \ar[r] & \{1\}.}\eqno\square$$

\bigskip{\it Interpr\'{e}tation de $H^1(G,A)$ en termes
d'extensions.} Soit $E$ une extension triviale de $G$ par $A$.
Choisissons une section $h:G\rightarrow E$ qui soit un homomorphisme
(ce qui identifie $E$ au produit semi-direct $G.A$). Soit $h'$ une
autre section; on peut \'{e}crire $h'$ de fa\c{c}on unique comme
$h'=l.h$, o\`{u} $l$ est une $1$-cocha\^{i}ne $G\rightarrow A$. On a
$f_{h'}=f_h.dl=dl$ puisque $f_h=1$. Pour que $h'$ soit un
homomorphisme, il faut et il suffit que $f_{h'}=1$, i.e. que $dl=1$,
autrement dit que $l$ soit un $1$-cocycle.

D'autre part, si on conjugue $h$ par un \'{e}l\'{e}ment $a$ de $A$,
on obtient une section qui est un homomorphisme. Soit $h'$ cette
section. A quoi cela correspond-il en termes de $l$? On a
$$h'(x)=ah(x)a^{-1}=l(x)h(x)$$
avec $l(x)=a{}^x\!a^{-1}$. Donc $l=df_a$ (o\`{u} $f_a$ est
l'\'{e}l\'{e}ment de $C^0(G,A)$ correspondant \`{a} $a$). Donc $l$
doit \^{e}tre un cobord. D'o\`{u}:

\begin{thm}
Les classes de conjugaison (par les \'{e}l\'{e}ments de $A$, ou de
$G$) des sections de $E$ qui sont des homomorphismes correspondent
bijectivement aux \'{e}l\'{e}ments du groupe de cohomologie
$H^1(G,A)$.
\end{thm}

[Noter que cette correspondance \emph{d\'{e}pend} du choix de $h$.
Une fa\c{c}on plus intrins\`{e}que de s'exprimer consiste \`{a} dire
que l'ensemble des classes de sections-homomorphismes est un espace
principal homog\`{e}ne (\og torseur\fg) sous l'action de
$H^1(G,A)$.]

\begin{coro}
Pour que les sections de $\pi$ qui sont des homomorphismes soient
conjugu\'{e}es, il faut et il suffit que $H^1(G,A)=\{0\}$.
\end{coro}

\section{Groupes finis: un crit\`{e}re de nullit\'{e}}\label{4.3}

Soit $G$ un groupe \`{a} $m$ \'{e}l\'{e}ments et soit $A$ un
$G$-module.

\begin{thm}
Soient $n\geqslant 1$ et $x\in H^n(G,A)$. On a $mx=0$.
\end{thm}

Soit $f\in Z^n(G,A)$ un $n$-cocycle repr\'{e}sentant $x$. Il faut
construire $F\in
C^{n-1}(G,A)$ tel que $dF=mf$.\\
Prenons $F_1(s_1,\dots,s_{n-1})=\sum_{s\in
G}{f(s_1,\dots,s_{n-1},s)}$. Comme $f\in Z^n(G,A)$, on a $df=0$.
Or
$$\begin{array}{rcl}
df(s_1,\dots,s_{n+1}) & = &
\displaystyle s_1f(s_2,\dots,s_{n+1})-f(s_1s_2,s_3,\dots,s_{n+1})+\cdots\\
{} & {} & {}\\
{} & {} &
\hfill{{} +(-1)^nf(s_1,\dots,s_ns_{n+1})+(-1)^{n+1}f(s_1,\dots,s_n)} \\
{} & = & 0.
\end{array}$$
Donc
$$\begin{array}{rcl}
\displaystyle \sum_{s_{n+1}\in G}{df(s_1,\dots,s_{n+1})} & = &
s_1F_1(s_2,\dots,s_n)-F_1(s_1s_2,\dots,s_n)+\cdots\\
{} & {} & \hfill{{}
+(-1)^nF_1(s_1,\dots,s_{n-1})+(-1)^{n+1}mf(s_1,\dots,s_n).}
\end{array}$$

On a utilis\'{e} le fait que si $s_{n+1}$ parcourt $G$, $s_ns_{n+1}$
aussi ($s_n$ \'{e}tant fix\'{e}). On a ainsi obtenu
$$(-1)^nmf(s_1,\dots,s_n)=dF_1(s_1,\dots,s_n).$$
On pose donc $F=(-1)^n F_1$ qui v\'{e}rifie $dF=mf$, d'o\`{u} le
r\'{e}sultat.~\findem

\begin{coro}
Si l'application $a\mapsto ma$ est un automorphisme de $A$ ($m$
\'{e}tant l'ordre de $G$) alors $H^n(G,A)=\{0\}$ pour tout
$n\geqslant 1$.
\end{coro}

En effet, $x\mapsto mx$ est alors un automorphisme de $C^n(G,A)$ qui
commute \`{a} $d$. Donc c'est un automorphisme de $H^n(G,A)$ par
passage au quotient. Or c'est dans ce cas l'application nulle
d'o\`{u} $H^n(G,A)=\{0\}$.~\findem

\begin{coro}
Si $G$ et $A$ sont finis d'ordres premiers entre eux alors
$H^n(G,A)=\{0\}$ pour tout $n\geqslant 1$.
\end{coro}

En effet $a\mapsto ma$ est alors un automorphisme de $A$.~\findem

\begin{coro}
Si $G$ et $A$ sont finis d'ordres premiers entre eux alors:
\begin{enumerate}
\item[(1)] Toute extension $E$ de $G$ par $A$ est triviale.

\item[(2)] Deux homomorphismes sections de $G\rightarrow E$ sont
conjugu\'{e}s par un \'{e}l\'{e}ment de $A$.
\end{enumerate}
\end{coro}

On a $H^n(G,A)=\{0\}$ si $n\geqslant 1$. Le cas $n=2$ donne $(1)$ et
le cas $n=1$ donne $(2)$ d'apr\`{e}s l'\'{e}tude faite en
\ref{extensions}.~\findem

\section{Extensions de groupes d'ordres premiers entre
eux}\label{4.4}

Nous allons \'{e}tendre certains r\'{e}sultats sur les extensions
d'un groupe $G$ par un groupe $A$ commutatif au cas o\`{u} $A$ est
r\'{e}soluble ou m\^{e}me quelconque.

\begin{thm}[Zassenhaus]\label{Zassen}
Soient $A$ et $G$ deux groupes finis d'ordres premiers entre eux et
consid\'{e}rons une extension $\{1\}\rightarrow A\rightarrow
E\rightarrow G\rightarrow \{1\}$. Alors:
\begin{enumerate}
\item[(1)] Il existe un sous-groupe de $E$ (\emph{suppl\'{e}mentaire
de $A$}) qui se projette isomorphiquement sur $G$ ($E$ est produit
semi-direct).

\item[(2)] Si $A$ ou $G$ est r\'{e}soluble, deux tels sous-groupes
sont conjugu\'{e}s par un \'{e}l\'{e}ment de $A$ (ou de $E$, cela
revient au m\^{e}me).
\end{enumerate}
\end{thm}

On raisonne par r\'{e}currence sur $|E|$; on peut supposer $A$ et
$G$ distincts de $\{1\}$.

{\it Premier cas: $A$ est r\'{e}soluble.} On d\'{e}montre d'abord le

\begin{lemme}\label{4.4.2}
Soit $X$ un groupe r\'{e}soluble non r\'{e}duit \`{a} $\{1\}$. Il
existe un nombre premier $p$ et un $p$-sous-groupe $Y$ de $X$
distinct de $\{1\}$ tel que $Y$ soit ab\'{e}lien \'{e}l\'{e}mentaire
et caract\'{e}ristique.
\end{lemme}

On rappelle qu'un $p$-groupe ab\'{e}lien est dit
\emph{\'{e}l\'{e}mentaire} si ses \'{e}l\'{e}ments distincts de $1$
sont d'ordre $p$ et qu'un sous-groupe d'un groupe $X$ est
caract\'{e}ristique s'il est stable par tout automorphisme de $X$.

\bigskip {\it D\'{e}monstration du lemme. } Soient $D^i(X)$
les d\'{e}riv\'{e}s successifs de $X$. Comme $X$ est r\'{e}soluble,
il existe $i$ tel que $D^i(X)$ est distinct de $\{1\}$ et
$D^{i+1}(X)$ est r\'{e}duit \`{a} $\{1\}$. Alors $D^i(X)$ est un
sous-groupe de $X$ ab\'{e}lien et diff\'{e}rent de $\{1\}$. De plus,
il est caract\'{e}ristique. Soit alors $p$ divisant l'ordre de
$D^i(X)$ et soit $Y$ le groupe des \'{e}l\'{e}ments de $D^i(X)$
d'ordre divisant $p$. Alors $Y$ est ab\'{e}lien, diff\'{e}rent de
$\{1\}$, caract\'{e}ristique (un automorphisme de $X$ transforme un
\'{e}l\'{e}ment d'ordre $p$ en un autre de m\^{e}me ordre) et est un
$p$-groupe \'{e}l\'{e}mentaire.~\findem

\bigskip {\it Retour \`{a} la d\'{e}monstration du th\'{e}or\`{e}me. }
Appliquons le lemme avec $X=A$ et $Y=A'$ et remarquons que $A'$ est
normal dans $E$: un automorphisme int\'{e}rieur de $E$ restreint
\`{a} $A$ est un automorphisme de $A$ (car $A$ est normal dans $E$)
et
laisse donc $A'$, qui est caract\'{e}ristique, invariant.\\
Si $A=A'$, alors $A$ est ab\'{e}lien et le th\'{e}or\`{e}me est
connu. Sinon, comme $A'$ est normal dans $E$, on peut passer au
quotient par $A'$ et on obtient la suite exacte
$$\xymatrix{\{1\} \ar[r] & A/A' \ar[r] & E/A' \ar[r] & G \ar[r] & \{1\}}.$$
La situation se d\'{e}crit par le diagramme suivant:
$$\xymatrix{& E \ar[d]\\ & E/A' \ar[d]\\ G \ar@{.>}[uur] \ar@{.>}[ur] \ar[r] & E/A}$$
Comme $E/A'$ est de cardinal strictement inf\'{e}rieur \`{a} celui
de $E$, l'hypoth\`{e}se de r\'{e}currence entra\^{i}ne que $G$ se
rel\`{e}ve en un sous-groupe $G'$ de $E/A'$. Soit $E'$ l'image
r\'{e}ciproque de $G'$ par la projection $E\rightarrow E/A'$. Alors
on a la suite exacte
$$\xymatrix{\{1\} \ar[r] & A' \ar[r] & E' \ar[r] & G' \ar[r] & \{1\}}.$$
Or $A'$ est ab\'{e}lien. D'apr\`{e}s le \S\ \ref{4.3}, on peut donc
relever $G'$ en un sous-groupe de $E'$. On obtient ainsi un
rel\`{e}vement de $G$ dans $E$.

Montrons que deux tels rel\`{e}vements $G'$ et $G''$ sont
conjugu\'{e}s par un \'{e}l\'{e}ment de $A$. On a
$$E=A.G'\;\mbox{ et }\; E=A.G''.$$
L'hypoth\`{e}se de r\'{e}currence, appliqu\'{e}e \`{a} $E/A'$,
montre qu'il existe $a\in A$ tel que $aG'a^{-1}$ et $G''$ aient
m\^{e}me image dans $E/A'$. Quitte \`{a} remplacer $G'$ par
$aG'a^{-1}$, on peut donc supposer que $A'.G'=A'.G''$. La
conjugaison par un \'{e}l\'{e}ment de $A$ de $G'$ et $G''$
r\'{e}sulte alors du cas ab\'{e}lien (cf. \S\ \ref{4.3}),
appliqu\'{e} \`{a} $A'.G'=A'.G''$.

\bigskip
{\it Deuxi\`{e}me cas: assertion $(1)$ dans le cas g\'{e}n\'{e}ral.}
Soit $p$ premier divisant l'ordre de $A$ et soit $S$ un $p$-Sylow de
$A$ (cf. \S\ \ref{2.2}). Soit $E'$ le normalisateur dans $E$ de $S$.
D'apr\`{e}s le \S\ \ref{2.3}, on a $E=A.E'$. Soit $A'=E'\cap A$;
$A'$ est normal dans $E'$ et l'on a la suite exacte
$$\xymatrix{\{1\} \ar[r] & A' \ar[r] & E' \ar[r] & G \ar[r] & \{1\}}.$$
Distinguons deux cas:\\
$\bullet$ Si $|E'|<|E|$, l'hypoth\`{e}se de r\'{e}currence permet de
relever $G$ dans $E'$, donc dans $E$.\\
$\bullet$ Si $|E'|=|E|$ alors $S$ est normal dans $E$ donc aussi
dans $A$. On passe au quotient:
$$\xymatrix{\{1\} \ar[r] & A/S \ar[r] & E/S \ar[r] & G \ar[r] & \{1\}}$$
avec $E/S$ de cardinal strictement inf\'{e}rieur \`{a} celui de $E$.
Par l'hypoth\`{e}se de r\'{e}currence, $G$ se rel\`{e}ve en $G_1$ de
$E/S$. Soit $E_1$ l'image r\'{e}ciproque de $G_1$ par la projection
$E\rightarrow E/S$. On a la suite exacte
$$\xymatrix{\{1\} \ar[r] & S \ar[r] & E_1 \ar[r] & G \ar[r] &\{1\}.}$$
Or $S$ est un $p$-groupe donc est r\'{e}soluble et l'on est
ramen\'{e} au premier cas.

\bigskip
{\it Troisi\`{e}me cas: assertion $(2)$ lorsque $G$ est
r\'{e}soluble.} Soient $G$ et $G'$ deux rel\`{e}vements de $G$ dans
$E$. On a
$$E=A.G'\;\mbox{ et }\; E=A.G''.$$
Soient $p$ un nombre premier et $I$ un sous-groupe ab\'{e}lien
normal diff\'{e}rent de $\{1\}$ de $G$ (cf. lemme \ref{4.4.2}) et
soit $\widetilde{I}$ son image r\'{e}ciproque dans $E$ par la
projection $E\rightarrow G$. Soient $I'=\widetilde{I}\cap G'$ et
$I''=\widetilde{I}\cap G''$. On a
$$A.I'=A.I''\; (=\widetilde{I}).$$
Les groupes $I'$ et $I''$ sont des $p$-Sylow de $\widetilde{I}$; il
existe donc $x\in \widetilde{I}$ tel que $I''=xI'x^{-1}$; si on
\'{e}crit $x$ sous la forme $ay$ avec $a\in A$ et $y\in I'$, on a
$I''=aI'a^{-1}$. Quitte \`{a} remplacer $I'$ par $aI'a^{-1}$, on
peut
donc supposer $I''=I'$.\\
Soit $N$ le normalisateur de $I'=I''$ dans $E$. On a $G'\subset N$
et $G''\subset N$. Si $N$ est distinct de $E$, l'hypoth\`{e}se de
r\'{e}currence appliqu\'{e}e \`{a} $N$ montre que $G'$ et $G''$ sont
conjugu\'{e}s. Si $N=E$, autrement dit si $I'$ est normal dans $E$,
l'hypoth\`{e}se de r\'{e}currence appliqu\'{e}e \`{a} $E/I'$ montre
qu'il existe $a\in A$ tel que $I'.aG'a^{-1}=I'.G''$. Puisque $I'$
est normal et contenu \`{a} la fois dans $G'$ et $G''$, cela
entra\^{i}ne
$$aG'a^{-1}=G'',$$
d'o\`{u} le r\'{e}sultat.~\findem

\bigskip\rmq L'hypoth\`{e}se \og $A$ ou $G$ est r\'{e}soluble\fg\ faite
dans $(2)$ est automatiquement satisfaite d'apr\`{e}s le
th\'{e}or\`{e}me de Feit-Thompson (cf. \S\ \ref{grpe resol}) disant
que tout groupe d'ordre impair est r\'{e}soluble.

\section{Rel\`{e}vements d'homomorphismes}\label{4.5}

Soient $\{1\}\rightarrow A\rightarrow E\stackrel{\pi}{\rightarrow}
\Phi\rightarrow\{1\}$ une suite exacte , $G$ un groupe et $\varphi$
un homomorphisme de $G$ dans $\Phi$. Peut-on relever $\varphi$ en un
homomorphisme $\psi$ de $G$ dans $E$?
$$\xymatrix{\{1\} \ar[r] & A \ar[r] & E \ar[r]^{\pi} & \Phi \ar[r] & \{1\}\\
{} & {} & {} & G \ar[u]_\varphi \ar@{.>}[ul]^\psi}$$ La question
\'{e}quivaut \`{a} celle du rel\`{e}vement de $G$ dans une extension
$E_\varphi$ de $G$ par $A$ associ\'{e}e \`{a} $\varphi$, d\'{e}finie
de la fa\c{c}on suivante:
$$E_\varphi=\{(g,e)\in G\times E\ |\
\varphi(g)=\pi(e)\}$$ muni de la loi de groupe habituelle pour le
produit cart\'{e}sien. Alors $A$ se plonge dans $E_\varphi$ par
$a\mapsto (1,a)$ et $E_\varphi$ se projette sur $G$ par
$(g,e)\mapsto g$.
$$\xymatrix{E_\varphi \ar[r] \ar[d] & E \ar[d]^\pi\\ G \ar[r]^\varphi & \Phi}$$
On a la suite exacte
$$\xymatrix{\{1\} \ar[r] & A \ar[r] & E_\varphi \ar[r] & G \ar[r] & \{1\}}.$$
(on dit parfois que $E_\varphi$ est \emph{l'image r\'{e}ciproque}
(\og pull-back\fg) de l'extension $E$ par l'homomorphisme
$\varphi$).

Voyons l'\'{e}quivalence des deux probl\`{e}mes. Soit $\psi$ un
rel\`{e}vement de $\varphi$. Alors l'ensemble
$G_\psi=\{(g,\psi(g)),\ g\in G\}$
est un sous-groupe de $E_\varphi$ qui est un rel\`{e}vement de $G$.\\
Soit maintenant $G'$ un rel\`{e}vement de $G$. Alors $G'$ est
form\'{e} de couples $(g,e)$ avec $g\in G$ et $e\in E$, chaque $g\in
G$ apparaissant dans un et un seul couple. Alors $\psi$ d\'{e}fini
par $\psi(g)=e$ est un homomorphisme
qui rel\`{e}ve $\varphi$.\\
De plus, deux rel\`{e}vements $\psi'$ et $\psi''$ sont conjugu\'{e}s
par $a\in A$ si et seulement si $G_{\psi'}$ et $G_{\psi''}$ sont
conjugu\'{e}s par $(1,a)\in E_\varphi$. Le \S\ \ref{4.4} donne alors
le

\begin{thm}
Soit $\{1\}\rightarrow A\rightarrow E \rightarrow
\Phi\rightarrow\{1\}$ une suite exacte et soit $\varphi$ un
homomorphisme d'un groupe $G$ dans le groupe $\Phi$. Supposons $G$
et $A$ finis d'ordres premiers entre eux. Alors:
\begin{enumerate}
\item[(1)] Il existe un homomorphisme $\psi$ de $G$ dans $E$ qui
rel\`{e}ve $\varphi$.

\item[(2)] Si $G$ ou $A$ est r\'{e}soluble, deux tels homomorphismes
sont conjugu\'{e}s par un \'{e}l\'{e}ment de $A$.
\end{enumerate}
\end{thm}

\App On se donne un homomorphisme $\varphi: G\rightarrow
\mathbf{GL}_n(\ZM/p\ZM)$ o\`{u} $p$ ne divise pas l'ordre de $G$. On
va voir qu'\emph{on peut relever $\varphi$ en $\varphi_\alpha:
G\rightarrow \mathbf{GL}_n(\ZM/p^\alpha\ZM)$ pour tout $\alpha
\geqslant 1$}.

Commen\c{c}ons par relever $\varphi$ en $\varphi_2$. On a la suite
exacte
$$\xymatrix{\{1\} \ar[r] & A \ar[r] & \mathbf{GL}_n(\ZM/p^2\ZM) \ar[r] & \mathbf{GL}_n(\ZM/p\ZM) \ar[r] & \{1\}}$$
o\`{u} $A$ est form\'{e} des matrices de la forme $1+pX$ avec $X$
matrice $n\times n$ modulo $p$ et o\`{u} l'application de
$\mathbf{GL}_n(\ZM/p^2\ZM)$ dans $\mathbf{GL}_n(\ZM/p\ZM)$ est la
r\'{e}duction modulo $p$. Le groupe $A$ est alors isomorphe \`{a}
$\MM_n(\ZM/p\ZM)$ qui est un $p$-groupe ab\'{e}lien. On peut donc
appliquer le th\'{e}or\`{e}me pr\'{e}c\'{e}dent et relever $\varphi$
en
$\varphi_2$ de mani\`{e}re essentiellement unique.\\
Le m\^{e}me argument permet de relever $\varphi_\alpha$ en
$\varphi_{\alpha+1}$. On a la suite exacte
$$\xymatrix{\{1\} \ar[r] & A \ar[r] & \mathbf{GL}_n(\ZM/p^{\alpha +1}\ZM) \ar[r] & \mathbf{GL}_n(\ZM/p^{\alpha}\ZM) \ar[r] & \{1\}\\
&& G \ar[u]^{\varphi_{\alpha +1}} \ar[ur]_{\varphi_{\alpha}}}.$$ On
peut passer \`{a} la limite projective: comme $\varprojlim
{(\ZM/p^\alpha\ZM)}=\ZM_p$, on obtient une repr\'{e}sentation
$$\varphi_\infty: \xymatrix{G \ar[r] & \mathbf{GL}_n(\ZM_p)\ \ar@{^{(}->}[r] & \mathbf{GL}_n(\QM_p)}.$$
Or $\QM_p$ est de caract\'{e}ristique $0$: ainsi, \`{a} partir d'une
repr\'{e}sentation en caract\'{e}ristique $p$, on en obtient une en
caract\'{e}risque $0$.

\chapter{Groupes r\'{e}solubles et sous-groupes de Hall}

Nous allons essayer de g\'{e}n\'{e}raliser les th\'{e}or\`{e}mes de
Sylow. Le probl\`{e}me \'{e}tait alors le suivant: soit $G$ un
groupe d'ordre $\prod_{p}{p^{\alpha(p)}}$ (o\`{u} $p$ est premier),
existe-t-il, pour tout nombre premier $p$, un sous-groupe de $G$
d'ordre
$p^{\alpha(p)}$?\\
On peut se demander si, plus g\'{e}n\'{e}ralement, pour tout $n$
divisant l'ordre de $G$, on peut trouver un sous-groupe de $G$
d'ordre $n$. C'est vrai si $G$ est nilpotent, mais  faux sans
hypoth\`{e}se sur $G$; m\^{e}me \og $G$ r\'{e}soluble \fg\ est
insuffisant: le groupe $\AM_4$, d'ordre $12$ est r\'{e}soluble et
n'a pas de sous-groupe d'ordre $6$. On va donc faire des
hypoth\`{e}ses plus restrictives sur $n$.

\section{$\Pi$-sous-groupes}\label{Pisg}

Soit $\Pi$ un ensemble de nombres premiers et soit $\Pi'$ son
compl\'{e}mentaire. Si $n\in \NM$, on \'{e}crit $n=n_{\Pi}n_{\Pi'}$,
avec $n_{\Pi}$ (resp. $n_{\Pi'}$) divisible uniquement par des
\'{e}l\'{e}ments de $\Pi$ (resp. $\Pi'$). Un groupe $G$ est
appel\'{e} un \emph{$\Pi$-groupe} si tous les facteurs premiers de
l'ordre de
$G$ appartiennent \`{a} $\Pi$.\\
Le probl\`{e}me consiste en la recherche des $\Pi$-sous-groupes d'un
groupe donn\'{e} et de ses $\Pi$-sous-groupes maximaux tels qu'ils
sont d\'{e}finis ci-dessous.

\begin{defi}
Soit $G$ un groupe et soit $\Pi$ un ensemble de nombres premiers. On
appelle \emph{$\Pi$-Sylow} ou \emph{$\Pi$-sous-groupe de Hall} de
$G$ un sous-groupe $H$ tel que $|H|=|G|_\Pi$.
\end{defi}\label{Hall}

\rmq Si $\Pi=\{p\}$, un $\Pi$-Sylow de $G$ est un $p$-Sylow de
$G$.

\begin{thm}[P. Hall]\label{5.1.1}
Soient $G$ un groupe r\'{e}soluble et $\Pi$ un ensemble de nombres
premiers. Alors:
\begin{enumerate}
\item[(1)] $G$ poss\`{e}de des $\Pi$-Sylow.

\item[(2)] Soient $S_{\Pi}$ un $\Pi$-Sylow de $G$ et $H$ un
$\Pi$-sous-groupe de $G$. Alors $H$ est contenu dans un conjugu\'{e}
de $S_{\Pi}$.
\end{enumerate}
\end{thm}

La d\'{e}monstration de ce th\'{e}or\`{e}me sera donn\'{e}e au \S\
\ref{5.4}.

\begin{coro}
Deux $\Pi$-Sylow d'un groupe r\'{e}soluble sont conjugu\'{e}s.
\end{coro}

L'hypoth\`{e}se \og r\'{e}soluble\fg\ est essentielle:

\begin{thm}\label{5.1.2}
Si pour tout ensemble $\Pi$ de nombres premiers, $G$ poss\`{e}de un
$\Pi$-Sylow, alors $G$ est r\'{e}soluble.
\end{thm}

La d\'{e}monstration sera donn\'{e}e au \S\ \ref{5.6}.

\begin{thm}[Burnside]\label{thmburn5.4}
Soient $p$ et $q$ deux nombres premiers. Tout groupe d'ordre
$p^aq^b$ ($a,b\in \NM$) est r\'{e}soluble.
\end{thm}
\label{Burnside2} En effet, la question de l'existence de
$\Pi$-Sylow ne se pose vraiment que si $\Pi=\{p\}$ ou $\{q\}$ ou
$\{p,q\}$. Les th\'{e}or\`{e}mes de Sylow (cf. \S\ \ref{2.2})
r\'{e}pondent dans les deux premiers cas et $G$ lui-m\^{e}me
convient dans le troisi\`{e}me. Le th. \ref{5.1.2} assure alors que
$G$ est r\'{e}soluble.~\findem

\bigskip En fait le th\'{e}or\`{e}me de Burnside sera d\'{e}montr\'{e} en annexe
(cf. th. \ref{thmburn}) par la th\'{e}orie des caract\`{e}res et il
sera utilis\'{e} dans la d\'{e}monstration du th. \ref{5.1.2}.

\section{Pr\'{e}liminaires: sous-groupes permutables}\label{5.2}

Nous allons d\'{e}montrer quelques lemmes sur les produits de
sous-groupes.\\
Soient $A$ et $B$ deux sous-groupes d'un groupe $G$. Notons $A.B$
l'ensemble des produits $ab$, o\`{u} $a\in A$ et $b\in B$.

\begin{lemme}\label{5.2.1}
Il y a \'{e}quivalence entre:
\begin{enumerate}
 \item[(1)] $A.B=B.A$.
 \item[(2)] $A.B$ est un sous-groupe de $G$.
\end{enumerate}
\end{lemme}

$(1)\Rightarrow (2)$ car si $A.B=B.A$, on a $A.B.A.B\subset
A.A.B.B\subset A.B$ et $(A.B)^{-1}\subset B.A=A.B$ et $A.B$ est un
sous-groupe de $G$.

$(2)\Rightarrow (1)$ Si $A.B$ est un sous-groupe de $G$, on a
$A.B=(A.B)^{-1}=B.A$.~\findem
\bigskip

On dit que deux groupes $A$ et $B$ sont \emph{permutables} si
$A.B=B.A$.

\begin{lemme}\label{5.2.2}
Soient $A_1,\dots,A_n$ des sous-groupes de $G$ deux \`{a} deux
permutables. Alors $A_1\dots A_n$ est un sous-groupe de $G$.
\end{lemme}

La d\'{e}monstration se fait par r\'{e}currence sur $n$. Le lemme
\ref{5.2.1} donne le cas $n=2$. D'apr\`{e}s l'hypoth\`{e}se de
r\'{e}currence, $A_1\dots A_{n-1}$ est un groupe. Il est permutable
avec $A_n$ car $A_1\dots A_{n-1}.A_n=A_1\dots A_n.A_{n-1}$, d'o\`{u}
apr\`{e}s $(n-1)$ op\'{e}rations $A_1\dots A_{n-1}.A_n=A_n.A_1\dots
A_{n-1}$. D'apr\`{e}s le lemme \ref{5.2.1}, $A_1\dots A_n$ est un
sous-groupe de $G$.~\findem

\pagebreak[2]
\begin{lemme}\label{5.2.3}
Il y a \'{e}quivalence entre:
\begin{enumerate}
 \item[(1)] $A.B=G$.
 \item[(1')] $B.A=G$.
 \item[(2)] $G$ op\`{e}re transitivement sur $G/A\times G/B$.
\end{enumerate}
De plus, si $G$ est fini, ces propri\'{e}t\'{e}s sont
\'{e}quivalentes \`{a} chacune des suivantes:
\begin{enumerate}
\item[(3)] $(G:A\cap B)=(G:A).(G:B)$. \item[(3')] $(G:A\cap
B)\geqslant(G:A).(G:B)$.
\end{enumerate}
\end{lemme}

En effet:\\
$(1) \Leftrightarrow(1')$ car si $A.B=G$, $A.B$ est un sous-groupe
et d'apr\`{e}s le lemme \ref{5.2.1}, on a $A.B=B.A$, donc $B.A=G$.

$(1)\Rightarrow (2)$ Il s'agit de prouver que pour tous $g_1,g_2\in
G$, il existe $g\in G$ tel que $g\in g_1A$ et $g\in g_2B$. Or par
hypoth\`{e}se, il existe $a\in A$ et $b\in B$ tels que
$g_1^{-1}g_2=ab$, d'o\`{u} $g_1a=g_2b^{-1}$ et l'\'{e}l\'{e}ment
$g=g_1a=g_2b^{-1}$ convient.

$(2)\Rightarrow (1)$ Le groupe $G$ op\`{e}re transitivement sur
$G/A\times G/B$. Prenons donc, pour tout $g_1\in G$, un
\'{e}l\'{e}ment $g\in G$ tel que $g\in 1.A$ et $g\in g_1.B$. Cela
entraine $g_1\in A.B$, i.e. $A.B=G$.

Soit maintenant $G$ un groupe fini. Montrons $(2)\Leftrightarrow
(3)$. Soit $\dot{1}$ l'image de l'\'{e}l\'{e}ment unit\'{e} de $G$
dans $G/A$ (resp. $G/B$). Le stablisateur de $(\dot{1},\dot{1})$ par
l'action de $G$ sur $G/A\times G/B$ est $A\cap B$. Le nombre $n$
d'\'{e}l\'{e}ments de l'orbite de $(\dot{1},\dot{1})$ est donc
l'indice $(G:A\cap B)$ de $A\cap B$ dans $G$. Or
$$\begin{array}{rcl} G\ \mbox{op\`{e}re
transitivement sur}\ G/A\times G/B & \Longleftrightarrow
& n=|G/A|\times |G/B| \\
{} & \Longleftrightarrow & n=(G:A)(G:B) \\
{} & \Longleftrightarrow & (G:A\cap B)=(G:A)(G:B), \\
\end{array}$$
ce qui est bien l'\'{e}quivalence entre $(2)$ et $(3)$.

$(3')\Leftrightarrow (3)$ car $(G:A\cap B)$ est le cardinal de
l'orbite de $(1,1)$, cardinal major\'{e} par celui de $G/A\times
G/B$, qui est $(G:A)(G:B)$.~\findem

\begin{lemme}\label{5.2.4}
Les propri\'{e}t\'{e}s du lemme \ref{5.2.3} sont vraies si les
indices de $A$ et $B$ dans $G$ sont premiers entre eux.
\end{lemme}
En effet, $(G:A\cap B)$ est divisible par $(G:A)$ et $(G:B)$ donc
par leur produit, ce qui prouve $(3')$ du lemme
pr\'{e}c\'{e}dent.~\findem

\section{Syst\`{e}mes permutables de sous-groupes de Sylow}\label{5.3}

Soit $G$ un groupe. Pour tout nombre premier $p$, choisissons un
$p$-Sylow $H_p$ de $G$. Nous dirons que le syst\`{e}me $\{H_p\}$ est
\emph{permutable} si les $H_p$ sont deux \`{a} deux permutables au
sens du \S\ \ref{5.2}. Dans ce cas, si $\Pi$ est un ensemble de
nombres premiers, le groupe $H_{\Pi}=\prod_{p\in\Pi}{H_p}$ est un
$\Pi$-sous-groupe de $G$.

\begin{thm}\label{5.3.1}
Si $G$ est r\'{e}soluble, $G$ poss\`{e}de un syst\`{e}me permutable
de sous-groupes de Sylow.
\end{thm}

La d\'{e}monstration se fait par r\'{e}currence sur l'ordre de $G$.
On suppose $G\neq\{1\}$; d'apr\`{e}s le lemme \ref{4.4.2}, il existe
alors un nombre premier $p_0$ et un $p_0$-sous-groupe normal $A$ de
$G$ distinct de $\{1\}$. D'apr\`{e}s l'hypoth\`{e}se de
r\'{e}currence, le groupe $G/A$ poss\`{e}de un syst\`{e}me
permutable $\{H'_p\}$ de $p$-Sylow. Soit $H'=\prod_{p\neq
p_0}{H'_p}$; c'est un sous-groupe de $G/A$ d'ordre $\prod_{p\neq
p_0}{|H'_p|}$. Soit $G'$ son image r\'{e}ciproque dans $G$, on a une
suite exacte:
$$\xymatrix{\{1\} \ar[r] & A \ar[r] & G' \ar[r] & H' \ar[r] &
\{1\}}.$$ Comme $A$ et $H'$ sont d'ordres premiers entre eux, il
existe un sous-groupe $H$ de $G$ qui rel\`{e}ve $H'$ (cf. \S 4.4).
Si $p\neq p_0$, posons $H_p$ le sous-groupe de $H$ qui rel\`{e}ve
$H'_p$; les $H_p$ sont des $p$-Sylow de $G$ deux \`{a} deux
permutables. Pour $p=p_0$, d\'{e}finissons $H_p$ comme l'image
r\'{e}ciproque de $H'_{p_0}$ dans $G$. C'est un $p_0$-Sylow de $G$
qui permute aux $H_p$ ($p\neq p_0$). Le syst\`{e}me $\{H_p\}$
r\'{e}pond donc \`{a} la question.~\findem

\section{D\'{e}monstration du th. \ref{5.1.1}}\label{5.4}
L'assertion $(1)$ sur l'existence de $\Pi$-Sylow r\'{e}sulte du th.
\ref{5.3.1} et du lemme \ref{5.2.2}.

Prouvons l'assertion $(2)$ par r\'{e}currence sur l'ordre de $G$.
Prenons comme au \S\ \ref{5.3} un  $p_0$-sous-groupe normal $A$ de
$G$ distinct de $\{1\}$. Soient $H'$ et $S'_{\Pi}$ les images
respectives de $H$ et $S_{\Pi}$ dans $G'=G/A$. D'apr\`{e}s
l'hypoth\`{e}se de r\'{e}currence $H'$ est contenu dans un
conjugu\'{e} de $S'_{\Pi}$. Quitte \`{a} remplacer $H'$ par un de
ses conjugu\'{e}s, on peut donc supposer $H'\subset S'_{\Pi}$. Il
faut maintenant
examiner deux cas:\\
$\bullet$ $p_0\in\Pi$. Alors $A\subset S_{\Pi}$ car $S_{\Pi}$
contient un $p_0$-Sylow $S_0$ de $G$ et comme $A$ est normal,
$A\subset S_0$ (cf. les th\'{e}or\`{e}mes de Sylow). L'inclusion
$H'\subset S'_{\Pi}$ donne alors $H\subset S_{\Pi}$.\\
$\bullet$ $p_0\notin\Pi$. Alors les ordres de $A$ et $S_{\Pi}$ sont
premiers entre eux et on a $A\cap H=\{1\}$, et $A\cap
S_{\Pi}=\{1\}$. Les projections $H\rightarrow H'$ et
$S_{\Pi}\rightarrow S'_{\Pi}$ sont des isomorphismes. Soit
$\widetilde{H}$ le sous-groupe de $S_{\Pi}$ qui se projette sur
$H'$; dans $A.H$ les groupes $H$ et $\widetilde{H}$  sont des
rel\`{e}vements de $H'$; ils sont donc conjugu\'{e}s (voir le \S\
\ref{4.4}, th. \ref{Zassen}).~\findem

\section{Un crit\`{e}re de r\'{e}solubilit\'{e}}\label{5.5}
\begin{thm}[Wielandt]\label{5.5.1}
Soit $G$ un groupe fini et soient $H_1$, $H_2$, $H_3$ trois
sous-groupes de $G$. Si les $H_i$ sont r\'{e}solubles et si leurs
indices sont premiers entre eux deux \`{a} deux, alors $G$ est
r\'{e}soluble.
\end{thm}

La d\'{e}monstration se fait par r\'{e}currence sur l'ordre de $G$.
Remarquons tout d'abord que $G=H_1.H_2$. En effet, les indices
$(G:H_1)$ et $(G:H_2)$ sont premiers entre eux, donc comme chacun
d'eux divise $(G:H_1\cap H_2)$, on a $(G:H_1\cap H_2)\geqslant (G:
H_1)(G:H_2)$ et d'apr\`{e}s le lemme \ref{5.2.3}, on a $G=H_1.H_2$.\\
On peut supposer que $H_1\not=\{1\}$. D'apr\`{e}s le lemme
\ref{4.4.2}, il existe un nombre premier $p$ et un $p$-sous-groupe
normal $A$ de $H_1$ diff\'{e}rent de $\{1\}$. On peut supposer que
$p$ ne divise pas $(G:H_2)$. Alors $H_2$ contient un $p$-Sylow de
$G$ , donc un conjugu\'{e} de $A$. Comme $G=H_1.H_2$, tout
conjugu\'{e} de $A$ est de la forme $h_2^{-1}h_1^{-1}Ah_1h_2$ avec
$h_i\in H_i$ ($i=1,2$) et comme $A$ est normal dans $H_1$, et qu'un
de ses conjugu\'{e}s est contenu dans $H_2$, tous ses conjugu\'{e}s
sont dans $H_2$.\\
Soit $\widetilde{A}$ le sous-groupe de $G$ engendr\'{e} par les
conjugu\'{e}s de $A$. Alors $\widetilde{A}$ est normal dans $G$ et
contenu dans $H_2$, donc $\widetilde{A}$ est r\'{e}soluble. Soit
$H'_i$ l'image de $H_i$ dans $G'=G/\widetilde{A}$. Les indices
$(G':H'_i)$ sont premiers entre eux deux \`{a} deux (car $(G':H'_i)$
divise $(G:H_i)$) et les $H'_i$ sont r\'{e}solubles. L'hypoth\`{e}se
de r\'{e}currence montre alors que $G'$ est r\'{e}soluble; donc $G$
est r\'{e}soluble.~\findem

\section{D\'{e}monstration du th. \ref{5.1.2}}\label{5.6}
Soit $p$ un nombre premier et soit $G$ un groupe. On appelle
\emph{$p$-compl\'{e}ment} de $G$ tout sous-groupe $H$ de $G$ qui est
un $p'$-Sylow o\`{u} $p'$ est l'ensemble des nombres premiers
diff\'{e}rents de $p$.

Nous allons d\'{e}montrer le th. \ref{5.1.2} sous la forme
apparemment plus forte suivante:

\begin{thm}
Si, pour tout nombre premier $p$, le groupe $G$ a un
$p$-compl\'{e}ment, alors $G$ est r\'{e}soluble.
\end{thm}

On raisonne par r\'{e}currence sur $|G|$. On distingue deux
cas.\\
$\bullet$ Le nombre des facteurs premiers de $|G|$ est $\leqslant
2$, autrement dit $|G|$ est de la forme $p^aq^b$, o\`{u} $p$ et $q$
sont premiers. D'apr\`{e}s un th\'{e}or\`{e}me de Burnside,
(d\'{e}montr\'{e} gr\^{a}ce \`{a} la th\'{e}orie des caract\`{e}res,
cf.
\S\ \ref{Burn}), $G$ est r\'{e}soluble.\\
$\bullet$ Le nombre des facteurs de l'ordre de $G$ est sup\'{e}rieur
ou \'{e}gal \`{a} $3$. Soient $p_i$ ($i=1, 2, 3$) de tels facteurs
et $H_i$ ($i=1, 2, 3$) un $p_i$-compl\'{e}ment de $G$. Alors les
indices $(G:H_i)$ pour $i=1, 2, 3$ sont premiers entre eux deux
\`{a}
deux.\\
De plus, $H_i$ poss\`{e}de un $p$-compl\'{e}ment pour tout nombre
premier $p$. En effet, si $p=p_i$, alors $H_i$ est son propre
$p_i$-compl\'{e}ment. Sinon, soit $H_p$ un $p$-compl\'{e}ment pour
$G$; comme $(G:H_i)$ et $(G:H_p)$ sont premiers entre eux, le lemme
\ref{5.2.4} montre que
$$(G:H_i\cap H_p)=(G:H_i)(G:H_p),$$ d'o\`{u} $(H_i:H_i\cap H_p)$ est la
plus grande puissance de $p$ divisant l'ordre de $H_i$. Comme
$H_i\cap H_p$ est un $p'$-groupe, $H_i\cap H_p$ est un
$p$-compl\'{e}ment de $H_i$. D'apr\`{e}s l'hypoth\`{e}se de
r\'{e}currence, $H_i$ est r\'{e}soluble. Mais alors $G$ satisfait
aux hypoth\`{e}ses du th. \ref{5.5.1}, donc $G$ est
r\'{e}soluble.~\findem

\chapter{Groupes de Frobenius}

\section{R\'{e}union des conjugu\'{e}s d'un sous-groupe}
\label{6.1}
\begin{thm}[Jordan]\label{remplissage}
Soit $G$ un groupe fini et soit $H$ un sous-groupe de $G$ distinct
de $G$; alors $\bigcup_{g\in G}{(gHg^{-1})}\neq G$. De fa\c{c}on
plus pr\'{e}cise, on a:
$$\Big|\bigcup_{g\in G}{gHg^{-1}}\Big|\leqslant \left|G\right|-\left(\frac{|G|}{|H|}-1\right).$$
\end{thm}

On sait \'{e}videmment que $1$ appartient \`{a} $H\cap (gHg^{-1})$
pour tout $g\in G$. On raisonne sur $G\mathbf{-}\{1\}$. On a
$$\bigcup_{g\in G}{\left(gHg^{-1}\mathbf{-}\{1\}\right)}=\bigcup_{g\in G/H}{\left(gHg^{-1}\mathbf{-}\{1\}\right)}$$
d'o\`{u}
$$\Big|\bigcup_{g\in G}{\left(gHg^{-1}\mathbf{-}\{1\}\right)}\Big|\leqslant \frac{|G|}{|H|}\big(|H|-1\big)$$
puis
$$\Big|\bigcup_{g\in G}{gHg^{-1}}\Big|\leqslant
|G|-\frac{|G|}{|H|}+1.\eqno\square$$

\bigskip On va voir que cette propri\'{e}t\'{e} reste vraie sans supposer que $G$ est fini,
pourvu que $G/H$ le soit. On utilise un lemme:

\begin{lemme}\label{lem6.2}
Soit $G$ un groupe et soit $H$ un sous-groupe de $G$ d'indice fini
$n$. Il existe un sous-groupe $N$ normal dans $G$ et contenu dans
$H$, tel que l'indice $(G:N)$ divise $n!$.
\end{lemme}

En effet, le groupe $G$ agit sur $X=G/H$ qui a $n$ \'{e}l\'{e}ments.
On obtient ainsi un homomorphisme $\varphi$ de $G$ dans le groupe
$\SM_X$ des permutations de $X$ qui est de cardinal $n!$. Le groupe
$N=\ker\varphi$ r\'{e}pond \`{a} la question.~\findem

\bigskip Si l'on applique le th. \ref{remplissage} \`{a} $G/N$ et
$H/N$, on constate que la r\'{e}union des conjugu\'{e}s de $H$ ne
remplit pas $G$ modulo $N$ donc, {\it a fortiori}, que
$\bigcup_{g\in G}{gHg^{-1}}\neq G$.~\findem

\bigskip\rmq Le cas $G=SO_3(\RM)$ et $H=\mathbf{S}_1$ montre que
l'hypoth\`{e}se $(G:H)<\infty$ ne peut \^{e}tre
supprim\'{e}e.\bigskip

Mentionnons deux reformulations du th. \ref{remplissage}:

{\bf Th\'{e}or\`{e}me 6.1'}\, {\it  Si un sous-groupe $H$ de $G$
rencontre toutes les classes de conjugaison de $G$, on a $H=G$.}

(C'est un crit\`{e}re souvent utilis\'{e} en th\'{e}orie des
nombres, $G$ \'{e}tant un groupe de Galois.)

{\bf Th\'{e}or\`{e}me 6.1''}\, {\it Si $G$ op\`{e}re transitivement
sur un ensemble $X$, et si $|X|\geqslant 2$, il existe un
\'{e}l\'{e}ment de $G$ qui op\`{e}re sans point fixe.}

En effet, si $H$ est le fixateur d'un point de $X$, on choisit un
\'{e}l\'{e}ment qui n'appartient \`{a} aucun conjugu\'{e} de
$H$.\findem

\bigskip Voici deux applications du th\'{e}or\`{e}me ci-dessus:

{\it Tout corps fini est commutatif (th\'{e}or\`{e}me de
Wedderburn).}\label{Wedderburn} En effet, soit $D$ un corps fini et
soit $F$ son centre. On sait que $(D:F)$ est un carr\'{e} $n^2$ et
que tout $x\in D$ est contenu dans un sous-corps $L$ commutatif,
contenant $F$ et tel que $(L:F)=n$. Comme deux tels sous-corps sont
isomorphes, le th\'{e}or\`{e}me de Skolem-Noether montre qu'ils sont
conjugu\'{e}s. Si $L$ est l'un d'eux et si l'on pose $G=D^*$ et
$H=L^*$, on a $G=\bigcup gHg^{-1}$, d'o\`{u} $G=H$, $n=1$ et $D$ est
commutatif.

{\it Racine d'une \'{e}quation modulo $p$.} Soit
$f=X^n+a_1X^{n-1}+\cdots+a_n$ un polyn\^{o}me \`{a} coefficients
dans $\ZM$, irr\'{e}ductible sur $\QM$. Si $p$ est un nombre
premier, notons $f_p$ la r\'{e}duction de $f$ modulo $p$; c'est un
\'{e}l\'{e}ment de $\FM_p[X]$. Notons ${\mathcal P}_f$ l'ensemble
des $p$ tels que $f_p$ ait une racines (au moins) dans $\FM_p$. On
va voir que la densit\'{e} de ${\mathcal P}_f$ est strictement
inf\'{e}rieure \`{a} $1$ d\`{e}s que $n\geqslant 2$. [On dit que
$\mathcal P$ a une densit\'{e} \'{e}gale \`{a} $\rho$ si
$$\frac{|\{p\leqslant x,\
p\in {\mathcal P}\}|}{|\{p\leqslant x\}|} \longrightarrow \rho\quad
\mbox{pour $x\rightarrow +\infty$.]}$$ Soit $X=\{x_1,\dots,x_n\}$
l'ensemble des racines de $f$ dans une extension de $\QM$ et soit
$G$ le groupe de Galois de $f$. Ce groupe op\`{e}re transitivement
sur $X$; on a $X\simeq G/H$ o\`{u} $H$ est le fixateur de $x_1$. On
d\'{e}montre (th\'{e}or\`{e}me de
Chebotarev-Frobenius)\label{Chebotarev-Frobenius} que la densit\'{e}
de ${\mathcal P}_f$ existe, et est \'{e}gale \`{a}
$$\frac{1}{|G|}\Big|\bigcup_{g\in
G}{gHg^{-1}}\Big|.$$ Vu le th\'{e}or\`{e}me ci-dessus, cette
densit\'{e} est $<1$.

\bigskip {\it Corollaire.} Si $n\geqslant 2$, il existe une infinit\'{e}
de $p$ tels que $f_p$ n'ait aucune racine dans $\FM_p$.

\bigskip Pour plus de d\'{e}tails, voir J.-P Serre, {\it On a theorem of Jordan}, Bull. A.M.S.
$40$ $(\mit 2003)$, $429-440$, reproduit dans Doc.Math.$1$, seconde \'edition, SMF,$2008$.

\section{Groupes de Frobenius: d\'{e}finition}\label{frob}

On va d\'{e}sormais s'int\'{e}resser aux couples $(G,H)$ tels que
$$\Big|\bigcup_{g\in G}{gHg^{-1}}\Big|=\left|G\right|-\left(\frac{|G|}{|H|}-1\right).$$
Cela signifie que $(gHg^{-1}\mathbf{-}\{1\})$ et
$(hHh^{-1}\mathbf{-}\{1\})$ sont disjoints si $g$ et $h$ ne sont pas
congrus modulo $H$, ou encore que $H$ et $gHg^{-1}$ sont
d'intersection r\'{e}duite \`{a} $\{1\}$ si $g\notin H$. On dit que
$H$
\og ne rencontre pas ses conjugu\'{e}s\fg.\\
On s'int\'{e}resse au cas o\`{u} $H$ est un sous-groupe propre de
$G$. Soit $X=G/H$. Une propri\'{e}t\'{e} \'{e}quivalente est que
tout \'{e}l\'{e}ment de $G$ (distinct de 1) a au plus un point fixe
dans $X$ si on fait agir $G$ sur $X$, ou encore: tout
\'{e}l\'{e}ment de $G$ qui fixe deux points est l'identit\'{e}.

\bigskip\ex Soit $G$ le groupe des transformations $h$ de la forme
$h(x)=ax+b$, o\`{u} $a$ et $b$ sont dans un corps fini $\FM$, avec
$a\neq 0$. Soit $H$ le sous-groupe $\{x\mapsto ax\}$. Si $N$ est le
sous-groupe de $G$ des translations, $N$ est normal dans $G$, et $G$
est le produit semi-direct de $H$ par $N$. Alors $(G,H)$ est un
exemple de la situation pr\'{e}c\'{e}dente.

\bigskip On va  g\'{e}n\'{e}raliser:

\begin{defi}
Un groupe $G$ est dit \emph{de Frobenius} s'il poss\`{e}de un
sous-groupe $H$ distinct de $\{1\}$ et de $G$ tel que
$\big|\bigcup_{g\in
G}{gHg^{-1}}\big|=\left|G\right|-\left({|G|}/{|H|}-1\right)$. On
parle dans ce cas pour $(G,H)$ de \emph{couple de Frobenius}.
\end{defi}\label{couplefrob}

\exs $(1)$ Soient $N$ et $H$ deux groupes finis, o\`{u} $H$ agit sur
$N$: \`{a} tout $h\in H$, on associe $\sigma_h:\ N\rightarrow N$
d\'{e}fini par $\sigma_h(n)=hnh^{-1}$. On a
$\sigma_{h_1h_2}=\sigma_{h_1}\circ\sigma_{h_2}$ pour tous
$h_1,h_2\in H$. Soit $G$ le produit semi-direct correspondant.
Cherchons \`{a} quelle condition $(G,H)$ est de Frobenius. Il faut
et il suffit que $H\cap nHn^{-1}=\{1\}$ pour tout $n\in N\mathbf{-}
\{1\}$. En effet, soit $h\in H\cap nHn^{-1}$; alors $h$ s'\'{e}crit
$nh'n^{-1}$ avec $h'\in H$. On quotiente modulo $N$, ce qui fournit
$h=h'$ ($G/N\simeq H$). Donc $h=nhn^{-1}$ soit encore
$n=h^{-1}nh=\sigma_{h^{-1}}(n)$. Donc $n$ est fix\'{e} par
$\sigma_{h^{-1}}$. Si $h\neq 1$, alors
n\'{e}cessairement $n=1$.\\
Une condition n\'{e}cessaire et suffisante pour que $(G,H)$ soit de
Frobenius est qu'il n'existe pas de couple $(h,n)$ avec $h\neq 1$ et
$n\neq 1$ tel que $\sigma_h(n)=n$, ou encore que $H$ op\`{e}re
librement sur $N\mathbf{-}\{1\}$. On a $\bigcup_{g\in
G}{gHg^{-1}}=\{1\}\cup (G\mathbf{-} N)$ d'o\`{u}
$G\mathbf{-}\bigcup_{g\in G}{gHg^{-1}}=~N\mathbf{-}\{1\}$ (en effet
$\bigcup_{g\in G}{gHg^{-1}}\subset\{1\}\cup (G\mathbf{-} N)$ et il
suffit de compter les \'{e}l\'{e}ments pour conclure).

$(2)$ Soit $p$ un nombre premier et soit $\FM$ un corps fini
contenant une racine $\xi$ $p$-i\`{e}me de l'unit\'{e}. Soit $N$
l'ensemble des matrices $p\times p$ triangulaires sup\'{e}rieures
\`{a} \'{e}l\'{e}ments diagonaux \'{e}gaux \`{a} $1$. C'est un
groupe. Soit $H$ le groupe cyclique engendr\'{e} par
$$\left(%
\begin{array}{ccccc}
  1 & 0 & \ldots & \ldots & 0 \\
  0 & \xi & \ddots &  & \vdots \\
  \vdots & \ddots & \ddots & \ddots & \vdots \\
  \vdots &  & \ddots & \xi^{p-2} & 0 \\
  0 & \ldots & \ldots & 0 & \xi^{p-1} \\
\end{array}%
\right).$$ Alors $H$ normalise $N$. Le groupe $G=N.H$ est le groupe
des matrices triangulaires sup\'{e}rieures avec des $\xi^k$ sur la
diagonale. On v\'{e}rifie que l'action de $H$ est sans point
fixe.\bigskip

Ces exemples sont en fait caract\'{e}ristiques; on a en effet le

\begin{thm}[Frobenius]
Soit $(G,H)$ un couple de Frobenius. Alors l'ensemble $N$ des
\'{e}l\'{e}ments de $G$ non conjugu\'{e}s \`{a} un \'{e}l\'{e}ment
de $H$ (ou \'{e}gaux \`{a} $1$) est un sous-groupe normal et on a
$G=N.H$.
\end{thm}\label{thmfrobenius}

Le point cl\'{e} de la d\'{e}monstration est que $N$ est
effectivement un sous-groupe: il repose sur la th\'{e}orie des
caract\`{e}res; nous le d\'{e}montrerons plus tard (cf. Annexe, th.
\ref{frobannexe}). Le groupe $N$ est alors normal (car invariant par
conjugaison). D'autre part:
$$\Big|\bigcup_{g\in G}{gHg^{-1}}\Big|=\left|G\right|-\big((G:H)-1\big),$$
donc $|N|=(G:H)$. Enfin $N\cap H=\{1\}$, d'o\`{u}
$G=N.H$.~\findem\bigskip

On d\'{e}montre (nous ne le ferons pas) qu'un groupe $G$ ne peut
\^{e}tre un groupe de Frobenius que \og d'une seule mani\`{e}re\fg:
si $(G,H_1)$ et $(G,H_2)$ sont des couples de Frobenius, alors $H_1$
est conjugu\'{e} de $H_2$. En particulier le sous-groupe normal $N$
est unique.

On cherche maintenant \`{a} classer les groupes de Frobenius en
\'{e}tudiant la structure de $N$ et celle de $H$.

\section{Structure de $N$}

On suppose que $N$ et $H$ sont distincts de $\{1\}$ et qu'ils
interviennent dans le groupe de Frobenius $G$. Choisissons $x\in H$
d'ordre premier $p$. L'\'{e}l\'{e}ment $x$ d\'{e}finit un
automorphisme de $N$ d'ordre $p$ sans point fixe $\neq 1$. D'o\`{u}:
$N$ intervient dans un groupe de Frobenius si et seulement s'il
poss\`{e}de un automorphisme $\sigma$ d'ordre premier et sans point
fixe distinct de $1$.

\begin{prop}\label{prop6.4}
Soit $\sigma$ un automorphisme d'ordre $p$ (non n\'{e}cessairement
premier) d'un groupe fini $N$ sans point fixe autre que $1$. Alors:
\begin{enumerate}
 \item[(1)] L'application $x\mapsto x^{-1}\sigma(x)$ (de $N$ dans
$N$) est bijective.
 \item[(2)] Pour tout $x\in N$, on a
$x\sigma(x)\sigma^2(x)\cdots\sigma^{p-1}(x)=1$.
 \item[(3)] Si $x$ et $\sigma(x)$ sont conjugu\'{e}s dans $N$ alors $x=1$.
\end{enumerate}
\end{prop}

$(1)$ Comme $N$ est fini, il suffit de montrer que l'application est
injective. Supposons $x^{-1}\sigma(x)=y^{-1}\sigma(y)$ avec $x,y\in
N$. Alors $yx^{-1}=\sigma(yx^{-1})$ donc l'\'{e}l\'{e}ment $yx^{-1}$
est fix\'{e} par $\sigma$ donc est \'{e}gal \`{a} $1$.

$(3)$ Soit $x\in N$ et supposons qu'il existe $a\in N$ tel que
$\sigma(x)=axa^{-1}$. D'apr\`{e}s $(1)$, il existe $b\in N$ tel que
$a^{-1}=b^{-1}\sigma(b)$. Alors
$\sigma(x)=\sigma^{-1}(b)bxb^{-1}\sigma(b)$ donc
$\sigma(bxb^{-1})=bxb^{-1}$, ce qui implique $bxb^{-1}=1$ puis
$x=1$.

$(2)$ Soit $a=x\sigma(x)\sigma^2(x)\cdots\sigma^{p-1}(x)$. On a
$$\sigma(a) = \sigma(x)\sigma^2(x)\cdots\sigma^{p-1}(x)x = x^{-1}ax,$$ donc $a=1$ d'apr\`{e}s $(3)$.~\findem

\begin{coro}
Si $l$ est un nombre premier, il existe un $l$-Sylow de $N$ stable
par $\sigma$.
\end{coro}

Soit $S$ un $l$-Sylow de $N$. Le groupe $\sigma(S)$ est aussi un
$l$-Sylow de $N$, donc il existe $a\in N$ tel que
$aSa^{-1}=\sigma(S)$. On \'{e}crit $a^{-1}=b^{-1}\sigma(b)$,
d'o\`{u} $\sigma(b^{-1})bSb^{-1}\sigma(b)=\sigma(S)$, soit
$bSb^{-1}=\sigma(b)\sigma(S)\sigma(b^{-1})=\sigma(bSb^{-1})$. Donc
$bSb^{-1}$ est un $l$-Sylow de $N$ stable par~ $\sigma$.~\findem

\begin{coro}
Si $a\in N$, l'automorphisme $\sigma_a: x\mapsto a\sigma(x)a^{-1}$
est conjugu\'{e} \`{a} $\sigma$ dans $\Aut(G)$; en particulier, il
est d'ordre $p$ et sans point fixe.
\end{coro}

D'apr\`{e}s \ref{prop6.4} $(1)$, il existe $b\in G$ tel que
$a=b^{-1}\sigma(b)$; alors $\sigma_a(x)=b^{-1}\sigma(bxb^{-1})b$,
donc $b\sigma_a(x)b^{-1}=\sigma(bxb^{-1})$ i.e. le diagramme suivant
est commutatif:
$$\xymatrix{N \ar[r]^{\sigma_a} \ar[d]_{{\rm conjugaison\ par }\; b^{-1}}& N \ar[d]^{{\rm \ conjugaison\ par }\; b^{-1}}\\
N \ar[r]_{\sigma} & N}$$ Ceci donne imm\'{e}diatement le
r\'{e}sultat.~\findem

\bigskip\exs $(1)$ Si $p=2$, on a $x\sigma(x)=1$ pour tout $x\in
N$, donc $\sigma(x)=x^{-1}$. Comme $\sigma$ est un automorphisme,
$N$ est ab\'{e}lien.

$(2)$ Pour le cas $p=3$ (Burnside), posons $\sigma(x)=x'$ et
$\sigma^2(x)=x''$. La prop. \ref{prop6.4}~$(2)$ appliqu\'{e}e \`{a}
$\sigma$ et $\sigma^2$ donne $xx'x''=1$ et $xx''x'=1$, donc $x'$ et
$x''$ commutent; pour les autres c'est \'{e}vident, donc $x$, $x'$
et $x''$ commutent. De m\^{e}me $x$ et $ax'a^{-1}$ commutent pour
tout $a$, ainsi que $x$ et $ax''a^{-1}$. Donc $x'$ et $x''$
commutent \`{a} tout conjugu\'{e} de $x$. Or $x=(x'x'')^{-1}$ donc
$N$ la propri\'{e}t\'{e} suivante: deux \'{e}l\'{e}ments
conjugu\'{e}s commutent, donc aussi $x$ et $(x,y)$. Finalement
$\big(x,(x,y)\big)=1$ pour tous $x,y\in N$. Le groupe d\'{e}riv\'{e}
de $N$ est contenu dans le centre de $N$; cela entra\^{i}ne que $N$
est nilpotent de classe au plus $2$.

$(3)$ Le cas $p=5$ a \'{e}t\'{e} trait\'{e} par Higman: le groupe
$N$ est alors nilpotent de classe au plus $6$ (c'est la meilleure
borne possible).\bigskip

Thompson a g\'{e}n\'{e}ralis\'{e} ces r\'{e}sultats par le

\begin{thm}[Thompson]
$N$ est nilpotent.
\end{thm}\label{thmThompson2}

Pour une d\'{e}monstration, cf. \cite[Kap. V, Haupsatz
$8$.$14$]{Hupp}. En ce qui concerne la classe de $N$, Higman a
conjectur\'{e} que si $p$ est l'ordre de $G$, la classe de $N$ est
$\leqslant\frac{p^2-1}{4}$.

\section{Structure de $H$}\label{6.4}

Nous dirons que $H$ a la propri\'{e}t\'{e} $\mathcal F$ s'il existe
un groupe $G$ contenant $H$ et distinct de $H$ tel que le couple
$(G,H)$ soit un couple de Frobenius. D'apr\`{e}s les
th\'{e}or\`{e}mes de Frobenius et de Thompson, cela revient \`{a}
dire qu'il existe un groupe nilpotent $N\neq\{1\}$ sur lequel $H$
op\`{e}re \emph{sans point fixe} (i.e. librement sur
$N\mathbf{-}\{1\}$).\bigskip

\ex Soit $\FM$ un corps fini de caract\'{e}ristique $l$ et soit $H$
un sous-groupe de $\mathbf{SL}_2(\FM)$ d'ordre premier \`{a} $l$. Si
l'on prend pour $N$ le $\FM$-espace vectoriel $\FM^2$, on
v\'{e}rifie facilement que $H$ op\`{e}re librement sur
$N\mathbf{-}\{0\}$. Donc $H$ a la propri\'{e}t\'{e} $\mathcal F$.
(Ceci s'applique notamment au groupe icosa\'{e}dral binaire d'ordre
$120$, groupe qui n'est pas r\'{e}soluble.)

\begin{thm}
Soit $H$ un groupe fini. Les propri\'{e}t\'{e}s suivantes sont
\'{e}quivalentes:
\begin{enumerate}
\item[(1)] $H$ a la propri\'{e}t\'{e} $\mathcal F$ (i.e. $H$ intervient
dans un couple de Frobenius).

\item[(2)] Il existe un corps $K$ et une repr\'{e}sentation lin\'{e}aire
$\rho: H\rightarrow \mathbf{GL}_n(K)$, avec $n\geqslant 1$, telle
que $\rho$ soit \og sans point fixe\fg (i.e. $H$ op\`{e}re librement
sur $K^n\mathbf{-}\{0\}$).

\item[(3)] Pour tout corps $K$ dont la caract\'{e}ristique ne divise
pas $|H|$, il existe une repr\'{e}sentation lin\'{e}aire $\rho:
H\rightarrow \mathbf{GL}_n(K)$ sans point fixe.

\item[(4)] On peut faire op\'{e}rer $H$ lin\'{e}airement et librement sur
une sph\`{e}re $\mathbf{S}_{n-1}$.
\end{enumerate}
\end{thm}

[Noter que $(2)$ et $(3)$ entra\^{i}nent que $\rho$ est fid\`{e}le.]

Avant de donner la d\'{e}monstration, faisons quelques remarques sur
la propri\'{e}t\'{e} suivante d'un corps $K$, not\'{e}e $(2_K)$: il
existe une repr\'{e}sentation $H\rightarrow \mathbf{GL}_n(K)$ sans
point fixe (avec $n\geqslant 1$).

(a) \emph{Cette propri\'{e}t\'{e} ne d\'{e}pend que de la
caract\'{e}ristique $p$ de $K$.}\\ En effet, si elle est satisfaite
par $K$ et si $x$ est un vecteur non nul de $K^n$, les
transform\'{e}s par $H$ de $x$ engendrent un espace vectoriel de
dimension finie $N$ sur le corps premier $K_0$ (i.e. $\FM_p$ ou
$\QM$), d'o\`{u} une repr\'{e}sentation $H \rightarrow
\mathbf{GL}_N(K_0)$ sans point fixe. Par extension des scalaires, on
en d\'{e}duit une pour tout corps contenant $K_0$.

(b) \emph{Si $(2_K)$ est vraie, la caract\'{e}ristique $p$ de
$K$ est, soit $0$, soit un nombre premier ne divisant pas $|H|$.}\\
En effet, si $H$ op\`{e}re librement sur $\FM_p^n\mathbf{-}\{0\}$,
l'ordre de $H$ divise $p^n-1$ et n'est donc pas divisible par $p$.

(c) \emph{La propri\'{e}t\'{e} $(2_K)$ en caract\'{e}ristique $0$
entra\^{i}ne la propri\'{e}t\'{e} analogue en toute
caract\'{e}ristique $p$ ne divisant pas $|H|$.}\\ En effet,
d'apr\`{e}s (a), il existe un $\QM$-espace vectoriel $V$ de
dimension finie $\geqslant 1$ o\`{u} $H$ op\`{e}re sans point fixe.
Soit $x\in V$ non nul et soit $L$ le $\ZM$-r\'{e}seau de $V$
engendr\'{e} par les transform\'{e}s de $x$ par $H$. Le groupe $H$
op\`{e}re sans point fixe sur $L$. Il op\`{e}re aussi sur le
$\FM_p$-espace vectoriel $V_p=L/pL$. Montrons que cette action est
sans point fixe, d\`{e}s lors que $p$ ne divise pas l'ordre de $H$.
Si $s\in H$ est d'ordre $m$, l'automorphisme $s_V$ de $V$ d\'{e}fini
par $s$ est tel que $s_V^m=1$, et il n'admet pas $1$ pour valeur
propre. On a donc
$$1+s_V+s_V^2+\cdots+s_V^{m-1}=0.$$
{\it A fortiori}, la m\^{e}me \'{e}quation vaut dans $V_p$ et elle
entra\^{i}ne (vu que $m$ est premier \`{a} $p$) que $sx\neq x$ pour
tout $x\in V_p$ non nul. D'o\`{u} la propri\'{e}t\'{e} $(2_K)$ pour
$\FM_p$.

(d) \emph{La propri\'{e}t\'{e} $(2_K)$ en caract\'{e}ristique $p\neq
0$ entra\^{i}ne la propri\'{e}t\'{e} analogue en caract\'{e}ristique $0$.}\\
Soit en effet $\rho_p:H\rightarrow \mathbf{GL}_n(\ZM/p\ZM)$ une
repr\'{e}sentation lin\'{e}aire de $H$ sans point fixe. D'apr\`{e}s
(b), $p$ ne divise pas $|H|$. D'apr\`{e}s ce qu'on a vu au chapitre
\ref{chap4}, on peut relever $\rho_p$ en un homomorphisme
$\rho_{p^\infty}: H\rightarrow \mathbf{GL}_n(\ZM_p)$, o\`{u}
$\ZM_p=\displaystyle{\lim_{\longleftarrow}{\ZM/p^\nu\ZM}}$ est
l'anneau des entiers $p$-adiques. Comme $\ZM_p\subset\QM_p$, on
obtient ainsi une repr\'{e}sentation lin\'{e}aire $H\rightarrow
\mathbf{GL}_n(\QM_p)$ de caract\'{e}ristique nulle. Cette
repr\'{e}sentation est sans point fixe. En effet, si
$x=(x_1,\dots,x_n)$ est un vecteur fixe non nul, on peut supposer
(quitte \`{a} multiplier $x$ par un scalaire) que les $x_i$
appartiennent \`{a} $\ZM_p$ et que l'un d'eux n'est pas divisible
par $p$. La r\'{e}duction modulo $p$ des $x_i$ donne alors un
vecteur non nul de $\FM_p^n$ fixe par $H$, contrairement \`{a}
l'hypoth\`{e}se.\bigskip

Ceci fait, la d\'{e}monstration du th\'{e}or\`{e}me est
imm\'{e}diate. En effet, il r\'{e}sulte de (a), (b), (c) et (d) que
$(2_K)$ est ind\'{e}pendante de $K$. D'o\`{u} l'\'{e}quivalence
$(2)\Leftrightarrow (3)$ du th\'{e}or\`{e}me. Prouvons maintenant:

$(1)\Rightarrow (2)$ Si $H$ op\`{e}re sans point fixe sur le groupe
nilpotent $N\neq\{1\}$, le centre $C$ de $N$ n'est pas r\'{e}duit
\`{a} $\{1\}$. Si $p$ est un facteur premier de $|C|$, le groupe
$\mathcal{C}_p$ des \'{e}l\'{e}ments $x\in N$ tels que $x^p=1$ est
un $\FM_p$-espace vectoriel non nul o\`{u} $H$ op\`{e}re sans point
fixe.

$(3)\Rightarrow (1)$ On choisit pour $K$ un corps fini. On obtient
ainsi une action sans point fixe de $H$ sur un groupe ab\'{e}lien
\'{e}l\'{e}mentaire.

$(4)\Rightarrow (2)$ On prend $K=\RM$.

$(3)\Rightarrow (4)$ On prend $K=\RM$. On obtient une
repr\'{e}sentation lin\'{e}aire $\rho: H\rightarrow
\mathbf{GL}_n(\RM)$ sans point fixe. Comme $H$ est fini, il existe
sur $\RM^n$ une forme quadratique d\'{e}finie positive invariante
par $H$ (prendre la somme des transform\'{e}s par $H$ de la forme
quadratique standard $\sum{x_i^2}$): quitte \`{a} conjuguer $\rho$,
on peut donc supposer que $\rho(H)$ est contenu dans le groupe
orthogonal ${\mathbf O}_n(\RM)$, donc laisse stable la sph\`{e}re
$\mathbf{S}_{n-1}$ d'\'{e}quation
$$\sum_{i=1}^{n}{x_i^2}=1.$$
Ceci ach\`{e}ve la d\'{e}monstration du th\'{e}or\`{e}me.~\findem

\bigskip\rmqs $(1)$ On peut classer les groupes $H$ ayant la propri\'{e}t\'{e}
$\mathcal F$, cf. \cite{Wolf}.

$(2)$ Dans $(4)$, on ne peut pas supprimer la condition que $H$
op\`{e}re lin\'{e}airement sur $\mathbf{S}_{n-1}$. Exemple:
$\mathbf{SL}_2(\FM_p)$ pour $p\geqslant 7$.\bigskip

\exo Soit $H$ un groupe ayant la propri\'{e}t\'{e} $\mathcal F$.
Montrer
que:\\
$(i)$ Tout sous-groupe ab\'{e}lien de $H$ est cyclique.

$(ii)$ Si $p$ et $q$ sont premiers, tout sous-groupe d'ordre $pq$
de $H$ est cyclique.

Inversement, si $H$ est r\'{e}soluble et si $(i)$ et $(ii)$ sont
v\'{e}rifi\'{e}es, alors $H$ a la propri\'{e}t\'{e} $\mathcal F$
(th\'{e}or\`{e}me de Vincent). Par contre, on peut montrer que le
groupe $\mathbf{SL}_2(\FM_{17})$ a les propri\'{e}t\'{e}s $(i)$ et
$(ii)$, mais pas la propri\'{e}t\'{e} $\mathcal F$.

\chapter{Transfert}

\section{D\'{e}finition}\label{7.1}

Soient $G$ un groupe et $H$ un sous-groupe de $G$ d'indice fini.
Soit $X=G/H$ l'ensemble des classes \`{a} gauche selon $H$. Pour
tout $x\in X$, on choisit un repr\'{e}sentant $\bar{x}$ de $x$ dans
$G$. Le groupe $G$ agit sur $X$. Si $s\in G$ et $x\in X$,
l'\'{e}l\'{e}ment $s\bar{x}$ de $G$ a pour image $sx$ dans $X$. Si
$\overline{sx}$ d\'{e}signe le repr\'{e}sentant de $sx$, il existe
donc $h_{s,x}\in H$ tel que $s\bar{x}=\overline{sx}\,h_{s,x}$. On
pose:
$$\Ver(s)=\prod_{x\in X}{h_{s,x}} \pmod{(H,H)},$$ o\`{u} le produit est
calcul\'{e} dans le groupe $H^{ab}=H/(H,H)$.

\begin{thm}[Schur]
L'application $\Ver:G\rightarrow H^{ab}$ d\'{e}finie ci-dessus est
un homomorphisme et ne d\'{e}pend pas du choix du syst\`{e}me de
repr\'{e}sentants $\{\bar{x}\}_{x\in X}$.
\end{thm}

On commence par montrer que l'application $\Ver$ est bien
d\'{e}finie. Soit donc $\{\bar{x}'\}_{x\in X}$ un autre syst\`{e}me
de repr\'{e}sentants; calculons $\Ver '(s)$ le produit relatif aux
$\bar{x}'$. L'\'{e}l\'{e}ment $\bar{x}'\in X$ a pour image $x$ dans
$X$; il existe donc $h_x\in H$ tel que $\bar{x}'=\bar{x}\,h_x$.\\ On
a:
$$\begin{array}{rcl}
s\,\bar{x}' & = & s\,\bar{x}\,h_x\\
 & = & \overline{sx}\,h_{s,x}\,h_x \\
 & = & \overline{sx}\,h_{sx}\,h^{-1}_{sx}\,h_{s,x}\,h_x \\
 & = & (\overline{sx})'\,h^{-1}_{sx}\,h_{s,x}\,h_x, \\
 \end{array}$$
d'o\`{u}
$$\begin{array}{rcl}
\Ver '(s) & = & \prod{h^{-1}_{sx}\,h_{s,x}\,h_x}\pmod{(H,H)} \\
& = & \left(\prod{h_{sx}}\right)^{-1}\prod{h_{s,x}}\prod{h_x}\pmod{(H,H)},\\
\end{array}$$
car $H^{ab}$ est ab\'{e}lien\footnote{Afin d'\'{e}viter les
lourdeurs d'\'{e}criture, il est sous-entendu que les produits sont
\'{e}tendus aux $x\in X$.}. Or $\prod{h_{sx}}=\prod{h_x}$ car $sx$
parcourt $X$ quand $x$ parcourt $X$, d'o\`{u}:
$$\Ver '(s)=\textstyle\prod{h_{s,x}}=\Ver(s)\pmod{(H,H)},$$ et
l'application $\Ver$ est bien d\'{e}finie.\\
Montrons \`{a} pr\'{e}sent qu'il s'agit d'un homomorphisme. Soient
$s,t\in G$, on a: $$\begin{array}{rcl} st\,\bar{x} & = &
s\,\overline{tx}\,h_{t,x} \\
 & = & \overline{stx}\,h_{s,tx}\,h_{t,x}, \\
 \end{array}$$
d'o\`{u}:
$$\begin{array}{rcl} \Ver(st) & = &\prod{h_{s,tx}\,h_{t,x}}\pmod{(H,H)} \\
 & = & \prod{h_{s,tx}}\prod{h_{t,x}}\pmod{(H,H)},
\end{array}$$
car $H^{ab}$ est ab\'{e}lien. Or $\prod{h_{s,tx}}=\prod{h_{s,x}}$ car $tx$ parcourt $X$ quand $x$ parcourt $X$.\\
D'o\`{u} $\Ver (st)=\prod{h_{s,x}}\prod{h_{t,x}}=\Ver (s)\Ver
(t)\pmod{(H,H)}$.~\findem

\bigskip Comme $H^{ab}$ est ab\'{e}lien, l'homomorphisme $\Ver$ se
factorise en un homomorphisme de $G^{ab}$ dans $H^{ab}$ (que nous
noterons encore $\Ver$), appel\'{e} \emph{transfert}.

\bigskip
\rmq Le transfert est \emph{fonctoriel} pour les isomorphismes: si
$\sigma$ est un isomorphisme du couple $(G,H)$ sur le couple
$(G',H')$, le diagramme
$$\xymatrix{
G^{ab}\ar[r]^\sigma \ar[d]_{\Ver} & G'^{ab} \ar[d]^{\Ver} \\
H^{ab} \ar[r]^\sigma & H'^{ab} }$$ est commutatif (il suffit de
constater que si $\{\bar{x}\}$ est un syst\`{e}me de
repr\'{e}sentants de
$G/H$ alors $\{\sigma(\bar{x})\}$ en est un pour $G'/H'$).\\
En particulier, si l'on prend $G=G'$, $H=H'$ et
$\sigma(x)=gxg^{-1}$, avec $g\in N_G(H)$, ceci montre que l'image de
l'homomorphisme $\Ver~:~G^{ab}~\rightarrow~H^{ab}$ est contenue dans
l'ensemble des \'{e}l\'{e}ments de $H^{ab}$ fix\'{e}s par $N_G(H)$.

\section{Calcul du transfert}\label{7.2}
Soit $H$ un sous-groupe d'indice fini de $G$ et soit $X=G/H$.
L'\'{e}l\'{e}ment $s\in G$ op\`{e}re sur $X$; soit $C$ le
sous-groupe cyclique de $G$ engendr\'{e} par $s$. Alors $C$
d\'{e}coupe $X$ en orbites $O_\alpha$. Soient $f_\alpha=|O_\alpha|$
et $x_\alpha\in O_\alpha$. On a $s^{f_\alpha}x_\alpha=x_\alpha$. Si
$g_\alpha$ est un repr\'{e}sentant de $x_\alpha$, on a donc:
$$s^{f_\alpha}g_\alpha =g_\alpha h_\alpha,\ \mbox{avec } h_\alpha \in H.$$

\begin{prop}\label{7.2.1}
On a $\Ver(s)=\prod_{\alpha}{h_\alpha }=\prod_{\alpha}{g_\alpha
^{-1} s^{f_\alpha} g_\alpha }\pmod{(H,H)}.$
\end{prop}

On prend pour syst\`{e}me de repr\'{e}sentants de $X$ les
\'{e}l\'{e}ments $s^ig_\alpha $, avec $0\leqslant i<f_\alpha $. Si
le repr\'{e}sentant de $x\in X$ est de la forme $s^{f_\alpha
-1}g_\alpha $, l'\'{e}l\'{e}ment $h_{s,x}$ correspondant de $H$ est
$h_\alpha $; les autres $h_{s,x}$ sont \'{e}gaux \`{a} $1$. La
proposition en r\'{e}sulte. ~\findem

\begin{coro}\label{7.2.2}
Soit $\varphi$ un homomorphisme de $H^{ab}$ dans un groupe $A$. On
suppose que $\varphi(h)=\varphi(h')$ lorsque les \'{e}l\'{e}ments
$h,h'\in H$ sont conjugu\'{e}s dans $G$. Alors $$\varphi\big(\Ver
(h)\big)=\varphi(h)^n,$$ pour tout $h\in H$, o\`{u} $n=(G:H)$.
\end{coro}

En effet, on a
$\varphi\big(\Ver(h)\big)=\prod_{\alpha}{\varphi(g_\alpha
^{-1}h^{f_\alpha }g_\alpha )}$. Or les \'{e}l\'{e}ments $g_\alpha
^{-1} h^{f_\alpha }g_\alpha $ et $h^{f_\alpha }$ sont conjugu\'{e}s
dans $G$, donc on a:
$$\textstyle{\varphi\big(\Ver(h)\big)=\prod_{\alpha}{\varphi(h^{f_\alpha })}=
\prod_{\alpha}{\varphi(h)^{f_\alpha }}}.$$ Le r\'{e}sultat se
d\'{e}duit alors de la relation $\sum_{\alpha}{f_\alpha
}=\sum_{\alpha}{|O_\alpha |}=|X|=n.$ ~\findem

\bigskip
Comme $H\subset G$, on a un morphisme naturel $H^{ab}\rightarrow
G^{ab}$.

\begin{coro}
L'homomorphisme compos\'{e}
$G^{ab}\stackrel{\Ver}{\longrightarrow}H^{ab}\longrightarrow G^{ab}$
est $s\mapsto s^n$.
\end{coro}

C'est imm\'{e}diat par la proposition, car
$$g_\alpha ^{-1}s^{f_\alpha
}g_\alpha = s^{f_\alpha }\pmod{(G,G)}\hspace{0.5cm}
\mbox{et}\hspace{0.5cm} \textstyle\sum_{\alpha}{f_\alpha
}=|X|=n.\eqno\square$$

\begin{coro}
Si $G$ est ab\'{e}lien, alors $\Ver:G\rightarrow H$ est donn\'{e}
par $s\mapsto s^n.$
\end{coro}

\section{Exemples d'utilisation du transfert}\label{7.3}

\subsection{Premier exemple (Gauss)}\label{Gauss}

On fixe un nombre premier $p\neq 2$.

Soit $G=\FM_p^*$ et soit $H=\{\pm 1\}$. Alors l'indice de $H$ dans
$G$ est $(p-1)/2$ et le transfert est donn\'{e} par $\Ver
(x)=x^{(p-1)/2}$ pour $x\in \FM_p^*$. Or c'est aussi le
\emph{symbole de Legendre $\big(\frac{x}{p}\big)$}. Ceci nous
donne alors un moyen de calculer $\big(\frac{x}{p}\big)$.\\
On choisit le syst\`{e}me de repr\'{e}sentants de $X=G/H$ donn\'{e}
par ${S=\{1,2,\dots,(p-1)/2\}.}$ Soit $x\in G$ et soit $s\in S$. Si
$xs\in S$, alors le $h_{s,x}$ correspondant vaut $1$; sinon c'est
$-1$. On pose donc:

$$\varepsilon(x,s) = \left\{\!\!\begin{array}{ll}
1 & \mbox{si}\ xs\in S, \\
-1 & \mbox{si}\  xs\notin S. \\
\end{array}\right.$$

On a $\Ver (x)=\prod_{s\in S}{\varepsilon}(x,s)$\quad (\emph{lemme
de Gauss}).

\bigskip

Calculons par exemple $\big(\frac{2}{p}\big)$ pour $p\neq 2$. On a
$p=1+2m$ et:

$$\begin{array}{rcll}
\big(\frac{2}{p}\big) & = & (-1)^{m/2} & \mbox{si $m$ est pair,} \\
 & & & \\
\big(\frac{2}{p}\big) & = & (-1)^{(m+1)/2} & \mbox{si $m$ est impair.} \\
\end{array}$$

D'o\`{u}:
\begin{eqnarray*}
\textstyle p\equiv 1 \pmod{8}  \Longrightarrow \big(\frac{2}{p}\big)=+1,\\
\textstyle p\equiv 3 \pmod{8}  \Longrightarrow \big(\frac{2}{p}\big)=-1,\\
\textstyle p\equiv 5 \pmod{8}  \Longrightarrow \big(\frac{2}{p}\big)=-1,\\
\textstyle p\equiv 7 \pmod{8}  \Longrightarrow \big(\frac{2}{p}\big)=+1,\\
\end{eqnarray*}
ce qui se r\'{e}sume par: $\big(\frac{2}{p}\big)=1
\Longleftrightarrow p\equiv \pm 1 \pmod{8}$.

\subsection{Second exemple}

\begin{prop} Si $G$ est un groupe sans torsion contenant un sous-groupe
$H$ d'indice fini isomorphe \`{a} $\ZM$, alors $G$ lui-m\^{e}me est
isomorphe \`{a} $\ZM$.
\end{prop}

Quitte \`{a} remplacer $H$ par l'intersection de ses conjugu\'{e}s,
on peut supposer $H$ normal dans $G$. Le groupe $G$ agit sur $H$,
d'o\`{u} un homomorphisme $\varepsilon:G\rightarrow\Aut (H)=\{\pm
1\}.$ Soit $G'$ le noyau de $\varepsilon$; alors $H\subset G'$ car
$H$, \'{e}tant ab\'{e}lien, agit trivialement sur lui-m\^{e}me par
automorphismes int\'{e}rieurs ($H$ est m\^{e}me contenu dans le
centre de $G'$ car celui-ci agit trivialement sur $H$). Donc le
transfert $\Ver :G'^{ab}\rightarrow H^{ab}=\ZM$ est \'{e}gal \`{a}
$x\mapsto x^n$ o\`{u} $n=(G':H)$. Soit $\Phi$ le noyau de $\Ver
:G'\rightarrow H^{ab}$, alors $\Phi\cap H=\{1\}$ car $H$ est
isomorphe \`{a} $\ZM$, donc $\Phi$ est fini et m\^{e}me $\Phi=\{1\}$
puisque $G$ est sans torsion. Donc $G'$ est isomorphe
\`{a} $\ZM$. Si $G$ est \'egal \`{a} $G'$, on a fini.\\
Sinon, on a $(G:G')=2$, $G'\simeq \ZM$ et $G/G'$ op\`{e}re sur $G'$
par $y\mapsto y^{\pm 1}$. Soit donc $x\in G\mathbf{-}G'$ tel que,
pour un certain $y\in G'$, on ait $xyx^{-1}=y^{-1}$. Alors $x^2\in
G'$ car l'indice de $G'$ dans $G$ est $2$. Prenons alors $y=x^2$, il
vient $xx^2x^{-1}=x^{-2}$ d'o\`{u} $x^2=x^{-2}$ donc $x=1$ car $G$
est sans torsion. Donc $G$ est isomorphe \`{a} $\ZM$.~\findem

\bigskip
\rmq On a un r\'{e}sultat analogue (mais bien moins
\'{e}l\'{e}mentaire) pour les groupes libres non
ab\'{e}liens\label{Stallings-Swan}
(Stallings-Swan\footnote{R\'{e}f\'{e}rences: J.R. Stallings, {\it On
torsion-free groups with infinitely many ends}, Ann. Math. $88$
$(\mit 1968)$, $312-334$. R. Swan, {\it Groups of cohomological
dimension one}, J. Algebra $12$ $(\mit 1969)$, $585-610$.}).

\section{Transfert dans un sous-groupe de Sylow}\label{7.4}
\smallskip
  A partir de maintenant, on suppose que $G$  est fini.
\medskip
\begin{thm}
Soient $H$ un $p$-Sylow d'un groupe $G$ et $\varphi:H\rightarrow A$
un homomorphisme \`{a} valeurs dans  un $p$-groupe ab\'{e}lien $A$.
Alors~:
\begin{enumerate}
 \item[(1)] Pour qu'il existe un prolongement de $\varphi$ en
un homomorphisme de $G$ dans $A$, il faut et il suffit que, pour
tous $h,h'\in H$ conjugu\'{e}s dans $G$, on ait
$\varphi(h)=\varphi(h').$

 \item[(2)] Si cette condition est
v\'{e}rifi\'{e}e, alors le prolongement est unique et donn\'{e} par:
$s\mapsto \varphi\big(\Ver(s)\big)^{1/n}$, o\`{u} $n=(G:H)$, ce qui
a un sens car $n$ est premier \`{a} $p$.
\end{enumerate}
\end{thm}

$(1)$ La condition est n\'{e}cessaire car si $\widetilde{\varphi}$
est un prolongement de $\varphi$ \`{a} $G$, alors pour $h\in H$ et
$g\in G$ avec $g^{-1}hg\in H$, on a
$$\varphi(g^{-1}hg)=\widetilde{\varphi}(g)^{-1}\varphi(h)\widetilde
{\varphi}(g)=\varphi(h)$$ car $A$ est ab\'{e}lien.\\
La condition est suffisante car, $n$ \'{e}tant premier \`{a} $p$ et
$A$ \'{e}tant un $p$-groupe, $\varphi\big(\Ver(s)\big)^{1/n}$ a un
sens (pour tout $a\in A$, il existe $b\in A$ unique tel que $b^n=a$)
et l'application $s\mapsto \varphi\big(\Ver(s)\big)^{1/n}$ convient
d'apr\`{e}s le cor. \ref{7.2.2}.

$(2)$ Ce prolongement est unique car $\varphi$ est
n\'{e}cessairement \'{e}gal \`{a} $1$ sur les $p'$-Sylow de $G$
lorsque $p'\neq p.$ ~\findem

\begin{thm}
Soient $H$ un $p$-Sylow ab\'{e}lien de $G$ et $N$ son normalisateur
dans $G$. Alors l'image de l'homomorphisme $\Ver:G^{ab}\rightarrow
H^{ab}=H$ est l'ensemble des \'{e}l\'{e}ments de $H$ fix\'{e}s par
$N$ (i.e. les \'{e}l\'{e}ments de $H$ qui sont dans le centre de
$N$).
\end{thm}

Par la remarque faite \`{a} la fin du \S\ \ref{7.1}, on a
d\'{e}j\`{a} l'inclusion de l'image de $\Ver$ dans $H^N=\{h\in H|\
xhx^{-1}=h,
\forall\, x\in N\}$.\\
Montrons l'\'{e}galit\'{e}. On a $N\supset H$ et $(N:H)$ est premier
\`{a} $p$ car $H$ est un $p$-Sylow. On d\'{e}finit l'homomorphisme
$\varphi:H\rightarrow H^N$ par $\varphi(h)=\big(\prod_{x\in
N/H}{xhx^{-1}}\big)^{1/(N:H)}.$\\
On a bien $\prod_{x\in N/H}{xhx^{-1}}\in H^N$ car si $x'\in N$, on
a:
$$x'\big(\textstyle\prod_{x\in
N/H}{xhx^{-1}}\big)x'^{-1}=\prod_{x\in
N/H}{x'xhx^{-1}x'^{-1}}=\prod_{x\in N/H}{xhx^{-1}}.$$ De plus, comme
$H$ est ab\'{e}lien, si $h,h'\in H$ sont conjugu\'{e}s dans $G$, ils
le sont dans $N$ (cf. \S\ \ref{2.4}) et on a alors
$\varphi(h)=\varphi(h')$. D'apr\`{e}s le cor. \ref{7.2.2}, on a
$$\varphi\big(\Ver(h)\big)=\varphi(h)^n,$$ pour $h\in H$, o\`{u} $n$ est l'indice de $H$ dans $G$.
Or, pour $h\in H^N$, on a $\varphi(h)=h$ et comme $\Ver(h)\in H^N$,
on a, pour $h\in H^N$:
$$\Ver(h)=\varphi\big(\Ver(h)\big)=\varphi(h)^n,$$ i.e.
$\Ver(h)=h^n$ si $h\in H^N$.\\
Comme $H^N$ est un $p$-groupe et que $n$ est premier \`{a} $p$, on
obtient tous les \'{e}l\'{e}ments de $H^N$ par les puissances
$n$-i\`{e}mes d'\'{e}l\'{e}ments de $H^N$. Donc $\Im(\Ver)=H^N.$
~\findem

\begin{thm}\label{7.4.3}
Soit $H$ un $p$-Sylow de $G$. Supposons que $H$ soit ab\'{e}lien
diff\'{e}rent de $\{1\}$ et que $G$ n'ait pas de quotient cyclique
d'ordre $p$. Soit $N$ le normalisateur de $H$ dans $G$. Alors~:
\begin{enumerate}
 \item[(1)] L'ensemble $H^N$ des \'{e}l\'{e}ments de $H$ fix\'{e}s par $N$ est
r\'{e}duit \`{a} $\{1\}$.
 \item[(2)] Si $r$ est  le rang de $H$
(nombre minimum de g\'{e}n\'{e}rateurs), il existe un nombre premier
$l$, distinct de $p$, qui divise \`{a} la fois $(N:H)$ et
$\prod_{i=1}^r(p^i-1).$
\end{enumerate}
\end{thm}

$(1)$ Si $H^N\neq\{1\}$, on a alors un homomorphisme non trivial
$\Ver:G\rightarrow H^N$, o\`{u} $H^N$ est un $p$-groupe, d'o\`{u}
l'on tire un quotient cyclique d'ordre $p$ de $G$.

Pour $(2)$, soit $H_p$ le sous-groupe de $H$ form\'{e} des
\'{e}l\'{e}ments $x$ tels que $x^p=1$; c'est un $\FM\!_p$-espace
vectoriel de dimension $r$ \big(puisque
$H=\prod_{i=1}^{r}(\ZM/p^{n_i}\ZM)$\big). L'action de $N$ sur $H_p$
est non triviale par $(1)$; elle d\'{e}finit un sous-groupe $\Phi$
de $\Aut(H_p)\simeq\mathbf{GL}_r(\ZM/p\ZM).$ Si $l$ est un facteur
premier de l'ordre de $\Phi$ alors $l$ divise l'ordre de $N/H$ car
$\Phi$ est un quotient de $N/H$ (en effet, $\Phi$ est d\'{e}fini par
l'action de $N$ sur $H$, qui est en fait une action de $N/H$ car $H$
\'{e}tant ab\'{e}lien, il agit trivialement sur lui-m\^{e}me). On a
$l\neq p$, car $p$ ne divise pas l'ordre de $N/H$; comme $\Phi$ est
un sous-groupe de $\mathbf{GL}_r(\ZM/p\ZM)$, $l$ divise l'ordre de
$\mathbf{GL}_r(\ZM/p\ZM)$, qui est
$p^{r(r-1)/2}\prod_{i=1}^{r}(p^i-1)$. D'o\`{u} $(2)$. ~\findem

\begin{coro}\label{coro7.9}
Si $p=2$, le sous-groupe $H$ n'est pas cyclique.
\end{coro}

En effet, on a alors $r\geqslant 2$ par le th. \ref{7.4.3}, mais on
peut en donner une d\'{e}monstration directe: si $H$ est un
$2$-Sylow cyclique, soit $h$ un g\'{e}n\'{e}rateur de $H$. Le groupe
agit sur lui-m\^{e}me par translations. L'\'{e}l\'{e}ment $h$
d\'{e}coupe alors $G$ en $|G/H|$ orbites. On envoie $G$ dans $\{\pm
1\}$ en associant \`{a} $x\in G$ la signature de la permutation que
$x$ effectue sur $G$ par l'action de translations. Or $h$ est
form\'{e} d'un nombre impair (i.e. pr\'{e}cis\'{e}ment $|G/H|$) de
cycles de la forme $(x,hx,\dots,h^{2^n-1}x)$, o\`{u} $|H|=2^n$.
Chacun de ces cycles est de signature $(-1)$, donc la signature de
$h$ est $(-1)$. On vient d'exhiber un homomorphisme non trivial de
$G$ dans $\{\pm 1\}$, d'o\`{u} contradiction. ~\findem

\bigskip
\rmq Le th. \ref{7.4.3} montre que $N_G(H)\neq H$ quand celui-ci est
un $p$-Sylow ab\'{e}lien d'un groupe sans quotient cyclique d'ordre
$p$ (sinon $H^N=H\neq\{1\}$).

\section{Application: groupes simples d'ordre
im\-pair inf\'{e}rieur \`{a} $2000$}

En fait, on va montrer qu'il n'existe pas de groupe $G$ d'ordre
impair, avec $1<|G|<~2000$, tel que $G=(G,G)$.\\
D'apr\`{e}s le th\'{e}or\`{e}me de Burnside (cf. Annexe, th.
\ref{thmburn}), il y a au moins trois facteurs premiers dans $|G|$;
si $p^{\alpha}$ est le plus petit de ces facteurs, on a
$p^{3\alpha}<2000$, ce qui donne $5$ cas possibles: $p^{\alpha}=3,5,
7, 9$ ou $11$.

\emph{Cas $p^{\alpha}=3$}: le groupe $G$ a un $3$-Sylow d'ordre $3$,
donc cyclique, donc ab\'{e}lien. Soit $N$ son normalisateur.
D'apr\`{e}s le th. \ref{7.4.3}, il existe $l$ premier distinct de
$3$ qui divise \`a la fois $|N|$ et $(p-1)=2$. Or $|N|$ est impair,
c'est impossible.

\emph{Cas $p^{\alpha}=5$}: il s'exclut par un raisonnement
analogue.

\emph{Cas $p^{\alpha}=9$}: on remarque de m\^{e}me que le $3$-Sylow
est d'ordre $3^{2}$, donc ab\'{e}lien. Et le m\^{e}me raisonnement
(dans les deux cas possibles: $r=1$ ou $r=2$) exclut ce cas.

\emph{Cas $p^{\alpha}=7$}: par le m\^{e}me th\'{e}or\`{e}me, on doit
avoir $l$ premier divisant $|N|$ (impair) et $p-1=6$. Donc $3$
divise $|G|$. Comme les cas pr\'{e}c\'{e}dents excluent
$p^{\alpha}=3$ ou $9$, on a: $3^{3}$ divise $|G|$. Par le
th\'{e}or\`{e}me de Burnside, on a $q$ premier, distinct de $3$ et
de $7$, divisant l'ordre de $G$. Donc $|G|\geqslant
3^{3}.7.q^{\beta}$ et $q^{\beta}\geqslant 11$ (car si $q=5$, un cas
d\'{e}j\`{a} vu montre que $\beta\geqslant2$). Or c'est impossible,
car $3^3.7.11>2000$.

\emph{Cas $p^{\alpha}=11$}: le th. \ref{7.4.3}, appliqu\'{e} au
$11$-Sylow, montre qu'il existe $l$ divisant $|N|$ et $p-1=10$.
D'apr\`{e}s un cas pr\'{e}c\'{e}dent, on a alors $|G|\geqslant
11.5^{2}.q^{\beta}$ et $q^{\beta}\geqslant 13$, ce qui est
impossible.\findem

\section{Application: groupes simples non ab\'{e}liens d'ordre inf\'{e}rieur \`{a} $200$}
Dans ce \S, on suppose $|G|\leqslant 200$.
\begin{prop}

\

\begin{itemize}
 \item[(1)] On suppose $G=(G,G)$ et $G\neq\{1\}$. Alors l'ordre
de $G$ est $60,120,168$ ou $180$.
 \item[(2)] Si $G$ est simple non ab\'{e}lien, alors l'ordre de $G$ est $60$ ou $168$ et
$G$ est isomorphe \`{a} $\AM_{5}$ ou $\mathbf{PSL}_{2}(\FM_{7})$.
\end{itemize}
\end{prop}
$(1)$ D'apr\`{e}s le \S\ pr\'{e}c\'{e}dent, l'ordre de $G$ est pair
et il est m\^{e}me divisible par $4$ d'apr\`{e}s le cor.
\ref{coro7.9}
qui affirme la non-existence de $2$-Sylow cycliques.\\
{\it Cas o\`{u} le $2$-Sylow $H$ est d'ordre $4$}. C'est alors
$\ZM/2\ZM\times\ZM/2\ZM$. Soit $N=N_{G}(H)$; $N$ agit non
trivialement sur $H$ (puisque $H^{N}=\{1\}$, cf. th. \ref{7.4.3}),
d'o\`{u} un homomorphisme non trivial $N\rightarrow
\Aut(\ZM/2\ZM\times\ZM/2\ZM)$ (qui est d'ordre $6$). Si $N$ s'envoie
sur un sous-groupe d'ordre $2$ de $\Aut H$, alors $H^{N}\simeq
\ZM/2\ZM$ (il suffit pour le voir d'examiner les automorphismes de
$\ZM/2\ZM\times\ZM/2\ZM$). Donc $3$ divise
l'ordre de $N$, donc aussi l'ordre de $G$. D'o\`{u} cinq possibilit\'{e}s:\\
$\bullet$ $|G|=4.3.13$: soit $H$ un $13$-Sylow et soit $N$ son
normalisateur. Alors $(G:N)$ est le nombre des $13$-Sylow de $G$
donc $(G:N)\equiv 1\pmod{13}$; comme $(G:N)$ divise $4.3$ cela
implique $(G:N)=1$ et donc $H$ est normal: c'est impossible.\\
$\bullet$ $|G|=4.3.11$: soit $N$ le normalisateur d'un $11$-Sylow
$H$; alors $(G:N)$ divise $4.3$ et $(G:N)\equiv 1 \pmod{11}$, donc
ou $(G:N)=1$ et c'est impossible, ou $(G:N)=12$ et alors $N=H$
puis
$H^N=H$ (puisque $H$ est ab\'{e}lien) et c'est impossible (cf. \S\ \ref{7.4}).\\
$\bullet$ $|G|=4.3.7$: soit $N$ le normalisateur d'un $7$-Sylow
$H$. Alors $(G:N)$ divise $12$ et $(G:N)\equiv 1\pmod{7}$; donc $(G:N)=1$: c'est impossible.\\
$\bullet$ Restent les cas $|G|=4.3.5=60$ ou $|G|=4.3^{2}.5=180$
(les autres cas sont \'{e}limin\'{e}s car sinon $|G|>200$).\\
{\it Cas o\`{u} le $2$-Sylow $H$ est d'ordre $8$}. On a deux cas
possibles: $|G|=8.3.5$ ou  $|G|=8.3.7$, les autres cas donnent un
ordre de $G$ trop grand. De m\^{e}me, le cas $|H|>8$ est
\'{e}limin\'{e} pour des raisons d'ordre. Ceci donne l'assertion
$(1)$.

$(2)$ Etudions les cas $|G|=4.3^{2}.5$ et $|G|=8.3.5$: soit $H$ un
$5$-Sylow et soit $N$ son normalisateur; alors $(G:N)\equiv 1
\pmod{5}$ et $(G:N)$ divise $4.3^{2}$ dans un cas, $8.3$ dans
l'autre. Dans les deux cas, la seule possibilit\'{e} est $(G:N)=6$.
Soit $X$ l'ensemble des $5$-Sylow de~$G$. Le groupe $G$ s'envoie
dans le groupe des permutations de $X$, c'est-\`{a}-dire $\SM_{6}$.
Comme $G$ est simple, il s'envoie dans $\AM_{6}$. Or $\AM_{6}$ est
d'ordre $360$ et $G$ d'ordre $180$ ou $120$. Le groupe $\AM_{6}$ ne
peut avoir de sous-groupe d'indice $m$ avec $1<m<6$ (ici ce serait
$3$ ou $2$) sinon $\AM_{6}$ se plongerait dans $\SM_{m}$, ce qui est
impossible puisque $|\AM_{6}|>|\SM_{m}|$.
Les seuls ordres possibles de groupes simples non ab\'{e}liens inf\'{e}rieurs \`{a} $200$ sont donc $4.3.5=60$ et $8.3.7=168$.\\

{\it Structure des groupes simples d'ordre $60$ et $168$}.

{\it Ordre $60$}: soit $H$ un $2$-Sylow de $G$; alors $H$ ne peut
pas \^{e}tre cyclique (cor. \ref{coro7.9}) et donc (th. \ref{7.4.3})
$3$ divise $|N|$ ($N$ normalisateur de $H$). Donc $12$ divise $|N|$
et $N\neq G$ donc $|N|=12$. Ainsi $G/N$ est d'ordre $5$ et on a un
homomorphisme non trivial de $G$ dans $\SM_{5}$. Comme $G$ est
simple, cela donne un plongement de $G$ dans $\AM_{5}$ et pour des
raisons d'ordre, on a $G=\AM_{5}$.

{\it Ordre $168$}: soit $H$ un $7$-Sylow de $G$ et soit $N$ son
normalisateur. Alors $(G:N)$ divise $8.3$ et $(G:N)\equiv
1\pmod{7}$. Comme $N\neq G$, on a $(G:N)=8$. Donc $|N|=21$.
Consid\'{e}rons la suite exacte $\{1\}\rightarrow H\rightarrow N
\rightarrow N/H \rightarrow \{1\}$. Comme l'ordre de $H$ est premier
\`{a} celui de $N/H$, le groupe $N$ est produit semi-direct de $N/H$
et de $H$ (cf. th. \ref{Zassen}). Le groupe $N$ a donc deux
g\'{e}n\'{e}rateurs: $\alpha$ (g\'{e}n\'{e}rateur de $H$) avec
$\alpha^{7}=1$ et $\beta$ (g\'{e}n\'{e}rateur de $N/H$) avec
$\beta^{3}=1$. L'automorphisme $x\mapsto \beta x\beta^{-1}$ de $H$
est d'ordre $3$; c'est donc, soit $x\mapsto x^{2}$, soit $x\mapsto
x^{-2}$. Quitte \`{a} remplacer $\beta$ par $\beta^{-1}$, on peut
supposer que c'est $x\mapsto
x^{2}$. On a alors $\beta\alpha\beta^{-1}=\alpha^{2}$.\\
Soit $X$ l'ensemble des $7$-Sylow de $G$. Alors $H$ op\`{e}re sur
$X$ et fixe lui-m\^{e}me comme \'{e}l\'{e}ment de $X$; appelons cet
\'{e}l\'{e}ment $\infty$. On a $X=\{\infty\}\cup X_{0}$ avec
$|X_{0}|=7$. Le groupe $H$ op\`{e}re librement sur $X_{0}$ (car  $H$
est cyclique d'ordre $7$). L'\'{e}l\'{e}ment $\beta$ op\`{e}re sur
$X$ et fixe $\infty$ car $\beta\in N$. Comme $\beta^3=1$, il existe
$x_{0}\in X_{0}$ tel que $\beta x_{0}=x_{0}$. Alors
$$X=\{x_{0},\alpha x_{0},\dots,\alpha^{6} x_{0},\infty \}.$$ On
identifie $X$ \`{a} $\mathbf{P}_{1}(\FM_{7})$ en indexant $\alpha^{i} x_{0}$ par $i$.\\
L'\'{e}l\'{e}ment $\alpha$ op\`{e}re sur $\mathbf{P}_{1}(\FM_{7})$
par : $\alpha(i)=i+1$ si $i<6$, $\alpha(6)=0$ et
$\alpha(\infty)=\infty $. L'\'{e}l\'{e}ment $\beta$ op\`{e}re par:
$\beta(\infty)=\infty$ , $\beta(0)=0$; de plus
$\beta\alpha=\alpha^{2}\beta$ d'o\`{u} $\beta(i+1)=\beta(i)+2$ et
$\beta(i)=2i$ pour tout $i$. Ainsi $\alpha$ agit sur
$\mathbf{P}_{1}(\FM_{7})$ comme une translation et $\beta$ comme une
homoth\'{e}tie. Soit $C$ le sous-groupe cyclique de $N$ engendr\'{e}
par $\beta$ et soit $M$ son normalisateur dans $G$. Comme $C$ est un
$3$-Sylow cyclique de $G$, $2$ divise $|M|$ (cf. th. \ref{7.4.3}).
Le groupe $M$ agit de fa\c{c}on non triviale sur $C$ (cf. \S\
\ref{7.4}) et donc il existe $\gamma$ tel que $\gamma
C\gamma^{-1}=C$ et $\gamma\beta\gamma^{-1}=\beta^{-1}$. Comme
$\gamma\notin C$ et $\gamma\neq \alpha^{n}$ (car $\alpha\notin M$), $\gamma$ peut \^{e}tre choisi d'ordre $2^{n}$.\\
L'\'{e}l\'{e}ment $\gamma$ transforme une orbite de $C$ en une
orbite de $C$ donc $\gamma(\{0,\infty\})=\{0,\infty\}$. Or $\gamma$
op\`{e}re sans point fixe sur $X$, car sinon il serait conjugu\'{e}
\`{a} un \'{e}l\'{e}ment de $N$, ce qui est impossible puisque $|N|$
est impair. Donc $\gamma(0)=\infty$ et $\gamma(\infty)=0$. Comme
$\gamma^{2}$ fixe $\infty$, on a $\gamma^2\in N$. Comme $\gamma$ est
d'ordre pair, on a $\gamma^{2}=1$. Donc $\gamma$ permute les deux
orbites $\{1,2,4\}$ et $\{3,6,5\}$. Posons $\gamma(1)=\lambda$.
Alors $\lambda$ est \'{e}gal \`{a} $3,6$ ou $5$. Comme
$\gamma\beta=\beta^{-1}\gamma$, on a
$\gamma(2i)\equiv\gamma(i)/2\pmod{7}$. D'o\`{u}
$\gamma(i)=\lambda/i$ et $\gamma$ est une homographie. D'o\`{u}
$\gamma\in \mathbf{PGL}_{2}(\FM_{7})$. Comme $-\lambda$ est un
carr\'{e}, on a
$$\gamma(i)=\frac{-\mu}{\mu^{-1}i}$$
avec $\mu^{2}=-\lambda$. Le d\'{e}terminant de la matrice $$\left(
                                                                \begin{array}{cc}
                                                                  0 & -\mu \\
                                                                  \mu^{-1} & 0 \\
                                                                \end{array}
                                                              \right)
$$ \'{e}tant \'{e}gal \`{a} $1$, on a $\gamma\in \mathbf{PSL}_{2}(\FM_{7})$.\\
Or $\alpha, \beta$ et $\gamma$ engendrent $G$. En effet, soit $G'$
le sous-groupe de $G$ engendr\'{e} par ces \'{e}l\'{e}ments. Alors
$G'$ contient $N$ et $\gamma$ (d'ordre pair). Si $G'\neq G$, alors
$G'$ est d'indice $2$ ou $4$ et $G$ s'envoie dans $\AM_{4}$ ou
$\AM_{2}$, ce qui est impossible. Donc $G=G'$. On a un homomorphisme
injectif $G\rightarrow \mathbf{PSL}_{2}(\FM_{7})$. Comme ces deux
groupes ont le m\^{e}me ordre, c'est un isomorphisme.\findem

\appendix

\chapter{Th\'{e}orie des caract\`{e}res}

\section{Repr\'{e}sentations et caract\`{e}res}
Soient $G$  un groupe, $K$  un corps et $V$ un espace vectoriel de
dimension finie $n$ sur le corps $K$. A partir du th.
\ref{formuledim}, on suppose que $K=\CM$ et que $G$ est \emph{fini}.
\begin{defi}
On appelle \emph{repr\'{e}sentation lin\'{e}aire} de $G$ dans $V$ la
donn\'{e}e d'un homomorphisme $\rho$ de $G$ dans $\mathbf{GL}(V)$.
La dimension de $V$ est appel\'{e}e le \emph{degr\'{e}} de la
repr\'{e}sentation.
\end{defi}\label{representation}

\rmqs $(1)$ On d\'{e}finit ainsi une op\'{e}ration de $G$  sur $V$
par $s.x=\rho (s)(x)$ pour $x\in V$ et $s\in G$.

$(2)$ Si $\rho$ est donn\'{e}, $V$ est appel\'{e} \emph{espace de
repr\'{e}sentation} de $G$ ou simplement \emph{repr\'{e}sentation}
de $G$; on \'{e}crit souvent $\rho_{_{V}}$ au lieu de $\rho$.

\bigskip Si $V_{1}$  et $V_{2}$ sont des repr\'{e}sentations de $G$ associ\'{e}es
aux homomorphismes $\rho_{1}$ et $\rho_{2}$, on peut d\'{e}finir:\\
$\bullet$ La somme directe $V_{1}\oplus V_{2}$ de $V_{1}$ et $V_{2}$
: c'est la repr\'{e}sentation $\rho : G \rightarrow
\mathbf{GL}(V_{1}\oplus V_{2})$ telle que $\rho (s)=
\rho_{1}(s)\oplus \rho_{2}(s)$. Si on choisit une base de
$V=V_{1}\oplus V_{2}$ associ\'{e}e \`{a} la d\'{e}composition en
somme directe, la matrice associ\'{e}e dans cette base \`{a} $\rho
(s)$ est la matrice:
$$ \left(\!\!\begin{array}{cc}  A_{1}(s) & 0 \\ 0 & A_{2}(s) \end{array}\!\! \right) $$
o\`{u} $A_{i}(s)$ est la matrice associ\'{e}e \`{a} $\rho_{i}(s)$
dans la base
de $V_{i}$ correspondante.\\
$\bullet$ Le produit tensoriel $V_{1}\otimes V_{2}$: $\rho
(s)(x\otimes y)= \rho_{1}(s)(x)\otimes \rho_{2}(s)(y)$ pour tout
$x\in V_{1}$ et tout $y\in V_{2}$.\\
$\bullet$ Le dual $V_{1}^{*}$ de $V_{1}$: $
\rho(s).l(x)=l\big(\rho_{1}(s^{-1}).x\big)$ pour tout $l\in
V_{1}^{*}$ et
tout $x\in V_1$.\\
$\bullet$ $\Hom(V_{1},V_{2})$ que l'on peut identifier \`{a}
$V_{1}^{*}\otimes V_{2}$: $\rho(s).h(x)=\rho_{2}(s)
h\big(\rho_{1}(s)^{-1}.x\big)$ pour tout $x\in V_1$.\\
etc...

\subsection*{Caract\`{e}re d'une repr\'{e}sentation}
\label{caractere} Soit $V$ un espace vectoriel muni d'une base
${(e_{i})}_{1\leqslant i\leqslant n}$ et $\rho$ une application
lin\'{e}aire de $V$ dans lui-m\^{e}me de matrice $a=(a_{ij})$. On
note $\Tr( \rho)= \sum_{i}{a_{ii}}$ la trace de la matrice $a$ (elle
est ind\'{e}pendante
de la base choisie).\\
Si, maintenant $V$ est une repr\'{e}sentation d'un groupe fini $G$,
on d\'{e}finit une fonction $\chi_{_{V}}$ sur $G$ \`{a} valeurs dans
$K$ par
$$\chi_{_{V}}(s)=\Tr \big(\rho_{_{V}}(s)\big)$$ o\`{u} $\rho_{_{V}}$ est
l'homomorphisme associ\'{e} \`{a} la repr\'{e}sentation $V$. La
fonction $\chi_{_{V}}$ est appel\'{e}e le \emph{caract\`{e}re} de la
repr\'{e}sentation $V$.

\bigskip\rmq On a $\chi_{_{V}}(1)=\dim V$.
\begin{ppts}

\

$\bullet$ $\chi_{_{V}}$ est une fonction centrale \rm{\big(i.e.
$\chi_{_{V}}(sts^{-1})=\chi_{_{V}}(t)$ si $s,t\in G$\big)},\\
$\bullet$ $\chi_{_{V_{1}\oplus V_{2}}}=\chi_{_{V_{1}}}+\chi_{_{V_{2}}}$,\\
$\bullet$ $\chi_{_{V_{1}\otimes V_{2}}}=\chi_{_{V_{1}}}\chi_{_{V_{2}}}$,\\
$\bullet$ $\chi_{_{V^{*}}}(s)=\chi_{_{V}}(s^{-1})$ {\it pour tout} $s\in G$,\\
$\bullet$ $\chi_{_{\Hom(V_{1},V_{2})}}(s)=\chi_{_{V_{1}}}(s^{-1})
\chi_{_{V_{2}}}(s)$ {\it pour tout} $s\in G$.
\end{ppts}

{\it On suppose d\'{e}sormais que  $K=\CM$ et que le groupe $G$ est
fini}. Soit $V$ une repr\'{e}sentation de $G$. On pose:
$$ V^{G}=\{ x\in V |\ s.x=x,\,\forall\, s\in G\}$$
$$\mbox{et}\;\;\;\pi(x)=\frac{1}{|G|}\sum_{s\in G}{s.x}\;\;\; \mbox{pour $x\in V$.}$$
On a $\pi(x)\in V^{G}$ et $\pi(x)=x$ si $x\in V^G$. Cela montre que
$\pi$ est un projecteur de $V$ sur $V^{G}$. L'application $\pi$
commute aux \'{e}l\'{e}ments de $G$: $\pi(s.x)=s.\pi(x)$ pour tout
$s\in G$; on a donc
$$V=V^{G}\oplus\ker\pi.$$
Dans une base associ\'{e}e \`{a} cette d\'{e}composition, la matrice
de $\pi$ s'\'{e}crit:
$$\left(\begin{array}{cccccc}1 & 0 & \cdots & \cdots & \cdots & 0 \\0 & \ddots & \ddots &  &  & \vdots \\\vdots & \ddots & 1 & \ddots &  & \vdots \\\vdots &  & \ddots & 0 & \ddots & \vdots \\\vdots &  &  & \ddots & \ddots & 0 \\0 & \cdots & \cdots & \cdots & 0 & 0\end{array}\right)$$
d'o\`{u}:
\begin{thm}\label{formuledim}
$\displaystyle\dim V^{G}=\Tr(\pi)=\frac{1}{|G|} \sum_{s\in G}
{\chi_{_{V}}(s)}$.
\end{thm}
\begin{coro}\label{VG}
Soit $0\rightarrow V'  \rightarrow V \stackrel{f}{\rightarrow} V''
\rightarrow 0$ une suite exacte de repr\'{e}sentations de $G$. Alors
tout \'{e}l\'{e}ment de $V''^{G}$ est image d'un \'{e}l\'{e}ment de
$V^{G}$.
\end{coro}
Soit $x''\in V''^{G}$. La suite \'{e}tant exacte, $f$ est surjective
et il existe $x\in V$ tel que $f(x)=x''$. Alors $\pi(x)\in V^{G}$ et
$f\big(\pi(x)\big)=\pi\big(f(x)\big)=f(x)=x''$.~\findem

\enlargethispage{\baselineskip}%
\begin{coro}\label{suplm}
Si $V'$ est un sous-espace vectoriel de $V$, stable sous l'action de
$G$, il existe un suppl\'{e}mentaire de $V'$ dans $V$ stable aussi
sous l'action de $G$.
\end{coro}
Consid\'{e}rons la suite exacte: $$\xymatrix{0  \ar[r] &
\Hom(V'',V') \ar[r] &  \Hom(V'',V) \ar[r] & \Hom(V'',V'') \ar[r] &
0}$$ o\`{u} $V''=V/V'$. Soit $x=\Id_{V''}\in \Hom(V'',V'')$; $x$ est
invariant sous l'action de $G$ donc, d'apr\`{e}s le cor. \ref{VG},
il existe $\varphi :V'' \rightarrow V$ invariant sous l'action de
$G$ et s'envoyant sur $x$. Donc $\varphi$ commute avec $G$:
$$(s^{-1}.\varphi)(v)=\varphi(v)=s^{-1}.\varphi(sv)$$
si $v\in V''$ et $s\in G$. Donc $s.\varphi(v)=\varphi(s.v)$.\\
L'homomorphisme $\varphi$ est une section (car $p\circ\varphi=x$
o\`{u} $x=\Id_{V''}$ et $p:V \rightarrow V''$ projection) donc
$V=V'\oplus \Im \varphi$ et $\Im\varphi$ est stable sous l'action de
$G$.~\findem

\begin{defi}
Soit $\rho : G \rightarrow \mathbf{GL}(V)$ une repr\'{e}sentation
lin\'{e}aire de $G$. On dit qu'elle est \emph{irr\'{e}ductible} si
$V\neq 0$ et si aucun sous-espace vectoriel de $V$ n'est stable par
$G$, \`{a} part $0$ et $V$.
\end{defi}

On a alors le

\begin{thm}\label{A5}
Toute repr\'{e}sentation est somme directe de repr\'{e}sentations
irr\'{e}ductibles.
\end{thm}
On raisonne par r\'{e}currence sur la dimension de la repr\'{e}sentation $V$.\\
C'est \'{e}vident si $\dim{V}\leqslant 1$. Sinon, ou bien $V$ est
irr\'{e}ductible, ou bien il existe un sous-espace vectoriel propre
$\neq 0$ de $V$ stable sous l'action de $G$ donc, d'apr\`{e}s le
cor. \ref{suplm}, il existe une d\'{e}composition en somme directe
$V=V' \oplus V''$ avec $\dim{V}' < \dim V$ et $\dim V'' < \dim V$ et
o\`{u} $V$ et $V'$ sont stables par $G$. On applique l'hypoth\`{e}se
de r\'{e}currence \`{a} $V'$ et $V''$: ils sont sommes directes de
repr\'{e}sentations irr\'{e}ductibles donc aussi $V$.~\findem

\medskip \rmq Il n'y a pas unicit\'{e} de la d\'{e}composition: si
par exemple $G$ op\`{e}re trivialement sur $V$, d\'{e}composer $V$
en somme directe de repr\'{e}sentations irr\'{e}ductibles revient
simplement \`{a} \'{e}crire $V$ comme somme directe de droites, ce
qui peut se faire de bien des fa\c{c}ons (si $\dim{V}\geqslant 2$).

\section{Relations d'orthogonalit\'{e}}\label{A.2}
\begin{thm}[lemme de Schur]\label{schur}
Soient $\rho_{1} :G \rightarrow \mathbf{GL}(V_{1})$ et $\rho_{2} :G
\rightarrow \mathbf{GL}(V_{2})$ deux repr\'{e}sentations
irr\'{e}ductibles de $G$ et soit $f$ un homomorphisme de $V_{1}$
dans $V_{2}$ tel que $\rho_{2}(s)\circ f=f\circ\rho_{1}(s)$ pour
tout $s\in G$. Alors:
\begin{enumerate}
 \item[(1)] Si $V_{1}$ et $V_{2}$ ne sont pas isomorphes, on a $f=0$.
 \item[(2)] Si $V_{1}=V_{2}$ et $\rho_{1}=\rho_{2}$, alors $f$ est une
homoth\'{e}tie.
\end{enumerate}
\end{thm}

$(1)$ Si $x\in \ker f$, on a
$f\big(\rho_{1}(s)\big).x=\big(\rho_{2}(s)\circ f\big).x=0$ pour
tout $s\in G$. Donc $\ker f$ est stable par $G$. Comme $V_1$ est
irr\'{e}ductible, $\ker f=0$ ou $\ker f=V_{1}$. Dans le premier cas
$f$ est injective, dans le second elle est nulle. De m\^{e}me $\Im
f$ est stable sous $G$, donc $\Im f=0$ ou $V_{2}$. Donc si $f\neq
0$, on a $\ker f=0$ et $\Im f=V_{2}$: $f$ est un isomorphisme de
$V_{1}$ sur $V_{2}$, d'o\`{u} $(1)$.

$(2)$ Soient maintenant $V_{1}=V_{2}$ et $\rho=\rho_{1}=\rho_{2}$.
Soit $\lambda\in\CM$ une valeur propre de $f$. Posons
$f'=f-\lambda$; $f'$ n'est pas injective. D'autre part $\rho(s)\circ
f'=\rho(s)\circ (f-\lambda)=f'\circ\rho(s)$. Donc, d'apr\`{e}s la
premi\`{e}re partie de la d\'{e}monstration, $f'=0$: $f$ est une
homoth\'{e}tie.~\findem

\subsection*{Orthogonalit\'{e} des caract\`{e}res}
 Soient $f$ et $g$ des fonctions sur $G$. On pose:
$$ \langle f,g\rangle=\frac{1}{|G|}\sum_{s\in G}{f(s) g(s^{-1})}.$$
C'est un produit scalaire.\bigskip

\rmq Si $g$ est un \emph{caract\`{e}re}, on a
$g(s^{-1})=\overline{g(s)}$, de sorte que $\langle f,g\rangle$ se
r\'{e}crit comme un produit hermitien:
$$\langle f,g\rangle=\frac{1}{|G|}\sum_{s\in G}{f(s) \overline{g(s)}}.$$

\begin{thm}[orthogonalit\'{e} des caract\`{e}res]
Soient $V$ et $V'$ deux repr\'{e}sentations irr\'{e}ductibles de
caract\`{e}res $\chi$ et $\chi '$. Alors:
$$ \langle \chi,\chi '\rangle=
 \left\{\!\!
 \begin{array}{ll}
 1 & \mbox{si\, $V=V'$ et\, $\chi = \chi '$}, \\
 0 & \mbox{si\, $V$ et\, $V'$ ne sont pas isomorphes}.
 \end{array}
 \right.$$
\end{thm}
$\bullet$ Consid\'{e}rons $W=\Hom(V,V)$; c'est une
repr\'{e}sentation de $G$. Soit $\varphi\in W$; alors $\varphi$ est
fix\'{e} par $G$ si et seulement si $\varphi\circ s=s\circ\varphi$
pour tout $s\in G$. Soit $W^{G}$ l'ensemble des $G$-endomorphismes
de $V$ (endomorphismes fix\'{e}s par $G$). D'apr\`{e}s le lemme de
Schur, on a $\dim W^G=1$, donc
$$1=\frac{1}{|G|}\sum_{s\in G}{\chi_{_{W}}(s)},$$ d'apr\`{e}s le th.
\ref{formuledim}. Or $\chi_{_{W}}(s)=\chi(s)\chi(s^{-1})$ donc
${\langle \chi,\chi\rangle=1}$.

$\bullet$ Si $V$ et $V'$ ne sont pas isomorphes et si l'on pose
$W=\Hom(V,V')$ et $W^{G}$ l'ensemble des $G$-homomorphismes de $V$
dans $V'$, le lemme de Schur implique que $\dim W^{G}=0$, d'o\`{u}
$\langle \chi ,\chi'\rangle=0$.~\findem

\begin{coro}
Les caract\`{e}res des diff\'{e}rentes repr\'{e}sentations
irr\'{e}ductibles de $G$ sont $\CM$-lin\'{e}airement
ind\'{e}pendants (on les appelle les caract\`{e}res
irr\'{e}ductibles de $G$).
\end{coro}

Ceci va nous permettre de montrer que les caract\`{e}res \og
caract\'{e}risent\fg\ les repr\'{e}sentations de $G$.

\begin{thm}
Soit $V$ une repr\'{e}sentation de $G$, de caract\`{e}re
$\chi_{_{V}}$, et soit $V=\bigoplus_{i=1}^{p} V_{i}$ une
d\'{e}composition de $V$ en somme directe de repr\'{e}sentations
irr\'{e}ductibles $V_{i}$ (de caract\`{e}res associ\'{e}s
$\chi_{_{V_{i}}}$). Alors si $W$ est irr\'{e}ductible de
caract\`{e}re $\chi_{_{W}}$, le nombre des $V_{i}$ isomorphes \`{a}
$W$ est \'{e}gal \`{a} $\langle \chi_{_{V}},\chi_{_{W}}\rangle$.
\end{thm}
On a:
\begin{eqnarray*}
\chi_{_{V}}&=&\sum_{i=1}^{p} \chi_{_{V_{i}}}, \\
\langle \chi_{_{V}},\chi_{_{W}}\rangle&=&\sum_{i=1}^{p}{\langle
\chi_{_{V_{i}}},\chi_{_{W}}\rangle}.
\end{eqnarray*}
Or d'apr\`{e}s le th\'{e}or\`{e}me pr\'{e}c\'{e}dent, $\langle
\chi_{_{V_{i}}},\chi_{_{W}}\rangle$ vaut $1$ ou $0$ selon que
$V_{i}$ est ou n'est pas isomorphe \`{a} $W$. D'o\`{u} le
th\'{e}or\`{e}me.~\findem

\begin{coro}

\

$\bullet$ Le nombre des $V_{i}$ isomorphes \`{a} $W$ ne d\'{e}pend
pas de la d\'{e}composition choisie {\rm (en ce sens, il y a
unicit\'{e} de la
d\'{e}composition)}.\\
$\bullet$ Deux repr\'{e}sentations de m\^{e}me caract\`{e}re sont
isomorphes.\\
$\bullet$ Si ${(W_{i})}_{1\leqslant i \leqslant n}$ sont les
repr\'{e}sentations irr\'{e}ductibles de $G$, toute
repr\'{e}sentation $V$ de $G$ est isomorphe \`{a} une somme
$\bigoplus n_{i}W_{i}$ (avec la convention:
$n_{i}W_{i}=W_{i}\oplus\dots\oplus W_{i}$, $n_{i}$ fois) o\`{u}
$n_{i}=\langle \chi_{_{V}},\chi_{_{W_{i}}}\rangle$.

\end{coro}

\section{Caract\`{e}res et fonctions centrales}\label{fctcentr}

On va maintenant chercher le nombre $h$ de repr\'{e}sentations
irr\'{e}ductibles (\`{a} isomorphismes pr\`{e}s),
c'est-\`{a}-dire le nombre des caract\`{e}res irr\'{e}ductibles.\\
Soit $\mathcal{C}$ le $\CM$-espace vectoriel des fonctions centrales
sur $G$ (c'est-\`{a}-dire l'ensemble des fonctions constantes sur
chaque classe de conjugaison de $G$). La dimension de $\mathcal{C}$
est \'{e}gale au nombre de classes de conjugaison de $G$.
\begin{thm}\label{A11}
Les caract\`{e}res irr\'{e}ductibles
$(\chi_{_{1}},\dots,\chi_{_{h}})$ de $G$ forment une base de
$\mathcal{C}$ (en particulier, $h$ est le nombre de classes de
conjugaison de $G$).
\end{thm}
Pour d\'{e}montrer ce th\'{e}or\`{e}me, on va utiliser le:

\begin{lemme}\label{homothetie}
Soit $V$ une repr\'{e}sentation irr\'{e}ductible de caract\`{e}re
$\chi$ et soit $n$ la dimension de $V$ (i.e. $n=\chi (1)$). Soit
$f\in \mathcal{C}$ et soit $\pi_{f}$ l'endomorphisme de $V$
d\'{e}fini par:
$$ \pi_{f}=\sum_{s\in G}{f(s^{-1})\rho_{_{V}}(s)}.$$
Alors $\pi_{f}$ est une homoth\'{e}tie de rapport $\displaystyle
\lambda=\frac{1}{n} \sum_{s\in
G}{f(s^{-1})\chi(s)}=\frac{|G|}{n}\langle f,\chi\rangle$.
\end{lemme}

Si $t\in G$, on a $\pi_{f}.t=\sum_{s\in
G}{f(s^{-1})\rho_{_{V}}(st)}=\sum_{s\in G}{
f(u^{-1})\rho_{_{V}}(tu)}$ (avec $u=t^{-1}st$), puisque $f$ est
centrale. Donc $\pi_{f}.t=t.\pi_{f}$. D'apr\`{e}s le th.
\ref{schur}, $\pi_{f}$ est une homoth\'{e}tie et
$$\Tr(\pi_{f})=n\lambda=\sum_{s\in G}{f(s^{-1})\chi(s)},$$ d'o\`{u}
$\displaystyle \lambda=\frac{|G|}{n}\langle
f,\chi\rangle$.~\findem

\bigskip
On en d\'{e}duit la d\'{e}monstration du th. \ref{A11}: si les
$(\chi_{i})$ ne forment pas une base de $\mathcal{C}$, il existe
$f\in \mathcal{C}$ non nulle et orthogonale aux $\chi_{i}$. Donc par
le lemme ci-dessus $\pi_{f}$ est nulle pour toute repr\'{e}sentation
irr\'{e}ductible, donc (par d\'{e}composition en somme directe) pour
toute repr\'{e}sentation. On va calculer $\pi_{f}$ pour une
repr\'{e}sentation particuli\`{e}re: soit $V$ un espace de dimension
$|G|$, muni d'une base ${(e_{s})}_{s\in G}$. Soit $\rho$
l'op\'{e}ration de $G$ sur $V$ d\'{e}finie par:
$$\rho (s).e_{t}=e_{st}.$$
La repr\'{e}sentation ainsi d\'{e}finie est dite
\emph{repr\'{e}sentation r\'{e}guli\`{e}re}\label{repr-reg} de $G$
(son caract\`{e}re se note $r_{G}$ et on a $r_{G}(s)=0$ si $s\neq1$
et $r_{G}(1)=|G|$). Calculons $\pi_f$ pour cette repr\'{e}sentation.
On a $\pi_{f}=\sum_{s\in G}{f(s^{-1})\rho(s)}$ donc
$\pi_{f}(e_{1})=\sum_{s\in G}{f(s^{-1})e_{s}}$ (car $e_{1.s}=e_s$).
Si $\pi_{f}$ est nulle, alors $f(s^{-1})=0$ pour tout $s\in G$ donc
$f$ est nulle (les ${(e_{s})}_{s\in G}$ forment une base) ce qui
contredit l'hypoth\`{e}se.~\findem

\bigskip
\rmq Chaque repr\'{e}sentation irr\'{e}ductible $W$ est contenue $n$
fois dans la repr\'{e}sentation r\'{e}guli\`{e}re (avec $n=\dim W$).
En effet, si $\chi$ est le caract\`{e}re de $W$, on a:
$$\langle r_{G},\chi\rangle\,=\frac{1}{|G|}\sum_{s\in G}{r_{G}(s^{-1})\chi(s)}=\chi(1)=\dim W.$$
Donc $r_{G}=\sum_{i=1}^{h}{\chi_{i}(1)\chi_{i}}$ avec $\chi_{i}$
irr\'{e}ductible. En particulier,
$$r_{G}(1)=|G|=\sum_{i=1}^{h}{\chi_{i}(1)\chi_{i}(1)}$$ d'o\`{u}
$$|G|=\sum_{i=1}^{h}{n_{i}^{2}}$$ avec $n_{i}=\chi_{i}(1)$.

\section{Exemples de caract\`{e}res}

On va d\'{e}terminer les caract\`{e}res irr\'{e}ductibles des
groupes $\AM_{n}$ ou $\SM_{n}$ pour $n\leqslant 4$.

$(1)$ Dans le cas trivial $\SM_{1}=\AM_{1}=\AM_{2}=\{1\}$ , il y a
un seul caract\`{e}re irr\'{e}ductible: $\chi=1$.

$(2)$ Dans le cas o\`{u} le groupe est $\SM_{2}=\{1,s\}$ (avec
$s^{2}=1$), tout caract\`{e}re irr\'{e}ductible est de degr\'{e} $1$
(par exemple car $|\SM_{2}|=2=\sum{n_{i}^{2}}$ o\`{u} $n_{i}$ est le
degr\'{e} de $\chi_{i}$ irr\'{e}ductible). On a deux caract\`{e}res
irr\'{e}ductibles:$$
\begin{array}{c||c|r}
{} & 1 & s \\
\hline\hline \chi_{_{1}} & 1 & 1 \\ \hline \chi_{_{2}} & 1 &
-1\end{array}$$

$(3)$  Le groupe $\AM_{3}=\{1,t,t^2|\ t^3=1\}$ est cyclique d'ordre
$3$. On a $|\AM_{3}|=3=\sum{n_{i}^{2}}$ d'o\`{u} trois
repr\'{e}sentations irr\'{e}ductibles de degr\'{e} $1$. Elles sont
donn\'{e}es par leurs caract\`{e}res:
$$\begin{array}{c||c|r|r} & 1 & t & t^{2} \\ \hline\hline\chi_{_{1}} & 1 & 1 & 1 \\
 \hline\chi_{_{2}} & 1 & \rho & \rho^{2} \\ \hline\chi_{_{3}} & 1 & \rho^{2} & \rho\end{array}$$
o\`{u} $\rho$ est une racine cubique de l'unit\'{e} $\neq 1$.

$(4)$ Le groupe sym\'{e}trique $\SM_{3}$ est le groupe di\'{e}dral
$\mathcal{D}_{3}$. On a trois classes de conjugaison $1,s,t$ avec
$s^{2}=1, t^{3}=1$ et $sts^{-1}=t^{-1}$. Deux caract\`{e}res
irr\'{e}ductibles de degr\'{e} $1$ sont donn\'{e}s par $1$ et la
signature. D'autre part
$\chi_{_{1}}+\chi_{_{2}}+n_{3}\chi_{_{3}}=r_{\SM_{3}}$ et
$\sum{n_{i}^{2}}=|\SM_{3}|=6=2+n_{3}^{2}$ donc $n_3=2$ et $\chi_{_{3}}$ est de degr\'{e} $2$.\\
Ensuite,
$$(\chi_{_{1}}+\chi_{_{2}}+2\chi_{_{3}})(s)=r_{\SM_3}(s)=0,$$ donc $\chi_{_{3}}(s)=0$, et de m\^{e}me
$$(\chi_{_{1}}+\chi_{_{2}}+2\chi_{_{3}})(t)=r_{\SM_3}(t)=0,$$
donc $\chi_{_{3}}(t)=-1$. On a donc:
$$\begin{array}{c||c|r|r} & 1 & s & t \\ \hline\hline\chi_{_{1}} & 1 & 1 & 1 \\ \hline\chi_{_{2}} & 1 & -1 & 1 \\ \hline\chi_{_{3}} & 2 & 0 & -1\end{array}$$
On peut trouver $\chi_{_{3}}$ directement en r\'{e}alisant $\SM_{3}$
comme groupe de sym\'{e}tries d'un triangle \'{e}quilat\'{e}ral.
Alors $s$ s'interpr\`{e}te comme sym\'{e}trie par rapport \`{a} une
droite, de matrice $\left(\begin{array}{cr}1&0\\0&-\!
1\end{array}\right)$ donc de trace nulle, et $t$ comme rotation
d'angle~$\frac{2\pi}{3}$ dont la trace est
$2\cos(\frac{2\pi}{3})=-1$.

$(5)$ Le groupe $\AM_{4}$ est d'ordre $12$. On a quatre classes de
conjugaison de repr\'{e}sentants $1, s, t, t'$ o\`{u} $s$ est
d'ordre $2$ ($s=(a,b)(c,d)$), $t$ d'ordre $3$ ($t=(a,b,c)$) et
$t'=t^{2}$. On a d\'{e}j\`{a} les trois repr\'{e}sentations du
quotient d'ordre $3$ de $\AM_{4}$. Puis
$\chi_{_{1}}+\chi_{_{2}}+\chi_{_{3}}+n_{4}\chi_{_{4}}=r_{\AM_{4}}$
et $\sum{n_{i}^{2}}=|\AM_{4}|=12$. D'o\`{u} la table de
caract\`{e}res:
$$\begin{array}{c||c|r|c|c} & 1 & s & t & t' \\ \hline\hline\chi_{_{1}} & 1 & 1 & 1 & 1 \\ \hline \chi_{_{2}} & 1 & 1 & \rho &\rho^{2}  \\ \hline\chi_{_{3}} & 1 & 1 &\rho^{2}  &\rho  \\ \hline \chi_{_{4}}& 3 & -1 & 0 & 0\end{array}$$
On peut trouver $\chi_{_{4}}$ directement en regardant $\AM_{4}$
comme un groupe de sym\'{e}tries d'un t\'{e}tra\`{e}dre regulier.
Alors $s$ s'interpr\`{e}te comme la sym\'{e}trie par rapport \`{a}
une droite joignant les milieux de deux c\^{o}t\'{e}s oppos\'{e}s:
sa matrice est
$$\left(\begin{array}{crr}1 & 0 & 0 \\0 & -\! 1 & 0
\\0 & 0 & -\! 1\end{array}\right)$$ dans une base convenable, donc
sa trace est $-1$. D'autre part, $t$ s'interpr\`{e}te comme permutation circulaire donc sa trace est nulle.\\
Un autre argument  permet de retrouver ce r\'{e}sultat: si $W$ est
une repr\'{e}sentation irr\'{e}ductible de dimension $1$ et si $V$
est irr\'{e}ductible, alors $V\otimes W$ est irr\'{e}ductible. En
particulier, $\chi_{_{4}}\chi_{_{2}}$ est irr\'{e}ductible et donc
\'{e}gal \`{a} $\chi_{_{4}}$ (pour des raisons de degr\'{e}).
D'o\`{u} $\chi_{_{4}}(t)=\chi_{_{4}}(t')=0$.

\bigskip\exo On peut donner de m\^{e}me la table des caract\`{e}res du
groupe $\SM_{4}$. Il y a cinq classes de conjugaison $1, \sigma,
s, t, \tau$ avec $\sigma=(a,b)$, $s=(a,b)(c,d)$, $t=(a,b,c)$ et
$\tau=(a,b,c,d)$. On trouve:$$\begin{array}{c||c|r|r|r|r} & 1 &
\sigma & s & t & \tau
\\ \hline\hline \chi_{_{1}} & 1 & 1 & 1 & 1 & 1 \\ \hline
\chi_{_{2}} & 1 & -\! 1 & 1 & 1 & -1 \\ \hline \chi_{_{3}} & 2 & 0
& 2 & -1 & 0 \\ \hline \chi_{_{4}} & 3 & 1 & -1 & 0 & -1 \\ \hline
\chi_{_{5}} & 3 & -1 & -1 & 0 & 1\end{array} $$

\section{Propri\'{e}t\'{e}s d'int\'{e}gralit\'{e}}
\label{propintegralite} Soit $G$ un groupe fini, soit $\rho$ une
repr\'{e}sentation de $G$ et soit $\chi$ son caract\`{e}re.
\begin{prop}\label{propalg}
Les valeurs de $\chi$ sont des entiers alg\'{e}briques.
\end{prop}

[On rappelle qu'un nombre complexe $x$ est un \emph{entier
alg\'{e}brique} s'il existe un entier $n>0$ et des \'{e}l\'{e}ments
$a_{0},\dots,a_{n-1}\in \ZM$ tels que:
$$x^{n}+a_{n-1}x^{n-1}+\cdots+a_{0}=0.$$

L'ensemble des entiers alg\'{e}briques forme un sous-anneau de
$\CM$.]

\bigskip\dem Le groupe $G$ \'{e}tant fini, il existe un entier $p>0$ tel que
$s^{p}=1$ pour tout $s\in G$. Donc $\rho(s^{p})=\rho(s)^{p}=1$ pour
tout $s$. Si $\lambda$ est une valeur propre de $\rho(s)$, alors
$\lambda^{p}$ est une valeur propre de $\rho(s)^{p}$, donc
$\lambda^{p}=1$. Toute valeur propre de $\rho(s)$ est un entier
alg\'{e}brique, or $\chi(s)$ est la trace  de $\rho(s)$ donc la
somme des valeurs propres de $\rho(s)$: c'est aussi un entier
alg\'{e}brique.~\findem

\bigskip\rmq Pour tout $s\in G$, $\rho(s)$ est diagonisable. En
effet, $\rho$ est une application de $G$ dans $\mathbf{GL}_n(\CM)$
et la matrice de $\rho(s)$, pour tout $s\in G$, est semblable \`{a}
une matrice de Jordan $J_{s}$, d\'{e}composable en blocs de la
forme:
$$\left(\begin{array}{cccc}\lambda & 1 & \cdots & \times \\
 & \ddots & \ddots & \vdots\\
 &  & \ddots & 1 \\ &  &  & \lambda \end{array}\right).$$
Alors $\rho(s)^{p}$ est semblable \`{a} une matrice d\'{e}composable
en blocs
$$\left(\begin{array}{cccc}\lambda^{p} & p & \cdots & \times \\ &\ddots &
\ddots &  \vdots\\
 &  & \ddots & p \\
 &  &  & \lambda^{p}
\end{array}\right).$$ Or $\rho(s)^{p}=1$, ceci n'est possible que
si $J_{s}$ est en fait diagonale, d'o\`{u} le r\'{e}sultat.

[Variante: appliquer le th. \ref{A5} au groupe cyclique engendr\'{e}
par $s$.]

\bigskip On va montrer comment la connaissance des caract\`{e}res d'un groupe
donne des renseignements sur ce groupe. Tout d'abord:
\begin{prop}\label{A14}
Soit $G$ un groupe distinct de $\{ 1\}$. Alors $G$ est simple si et
seulement si, pour tout $s\in G\mathbf{-}\{1\}$ et tout
caract\`{e}re $\chi$ irr\'{e}ductible et distinct de $1$, on a:
$$\chi(s)\neq\chi(1).$$
\end{prop}

Soient ${(\lambda_{i})}_{1\leqslant i \leqslant n}$ les valeurs
propres de $\rho(s)$. On sait que les $\lambda_{i}$ sont des racines
de l'unit\'{e}, donc que $|\lambda_{i}|=1$. Comme
$\chi(s)=\lambda_{1}+\cdots+\lambda_{n}$ et $\chi(1)=n$, on a
$\chi(s)=\chi(1)=n$ si et seulement si $\lambda_{i}=1$ pour tout
$i$\footnote{\label{footnote} On utilise le r\'{e}sultat
\'{e}l\'{e}mentaire suivant: si $\alpha_1, \dots, \alpha_n\in \CM$
sont de module $1$, et si $|\alpha_1+\cdots+\alpha_n|=n$, alors
$\alpha_1=\cdots=\alpha_n$.}. Donc $\chi(s)=n$ si et seulement si
$s\in \ker
\rho$.\\
Supposons maintenant que $G$ soit simple. Alors $\ker\rho$ est un
sous-groupe normal de $G$, c'est donc $\{1\}$ ou $G$. Donc si
$\chi(s)=\chi(1)$ pour un $s\in G\mathbf{-}\{1\}$ alors $\rho=\Id$ et donc $\chi=1$.\\
R\'{e}ciproquement, si $G$ n'est pas simple, soit $N$ un sous-groupe
normal non trivial de $G$. Soit $\chi'$ un caract\`{e}re non trivial
de $G/N$ associ\'{e} \`{a} la repr\'{e}sentation $\rho'$. Alors
$$\xymatrix{\relax G \ar[r] & G/N \ar[r]^-{\rho'} & \mathbf{GL}_{n}(\CM)}$$
d\'{e}finit une repr\'{e}sentation $\rho$ de $G$ de caract\`{e}re
$\chi\neq 1$ et telle qu'il existe $s\in G\mathbf{-}\{1\}$ avec
$\chi(s)=\chi(1)$.~\findem

\bigskip On va maintenant g\'{e}n\'{e}raliser la prop. \ref{propalg}:
\begin{thm}\label{thmA}
Soit $\rho$ une repr\'{e}sentation irr\'{e}ductible d'un groupe $G$,
de caract\`{e}re associ\'{e} $\chi$. Soit $f$ une fonction centrale
sur $G$, dont les valeurs sont des entiers alg\'{e}briques. Posons
$n=\chi(1)$. Alors $n^{-1}\sum_{s\in G}{f(s)\chi(s)}$ est un entier
alg\'{e}brique.
\end{thm}

Il suffit de d\'{e}montrer le th\'{e}or\`{e}me pour les fonctions
valant $1$ sur une classe de conjugaison et $0$ ailleurs. Soit $s\in
G$; notons $Cl(s)$ la classe de conjugaison de $s$ et $c(s)$ le
cardinal de $Cl(s)$. Soit $f_{s}$ la fonction d\'{e}finie sur $G$,
valant $1$ sur $Cl(s)$ et $0$ ailleurs. On~a:
$$\frac{1}{n}\sum_{t\in G}{f_{s}(t)\chi(t)}=\frac{1}{n}\sum_{t\in Cl(s)}{\chi(t)}=\frac{c(s)\chi(s)}{n}$$
car $\chi$ est centrale. Il suffit donc de d\'{e}montrer le:

\begin{thm}\label{thmB}
Avec les hypoth\`{e}ses du th. \ref{thmA} et les notations
ci-dessus, $\frac{c(s)\chi(s)}{n}$ est un entier alg\'{e}brique pour
tout $s\in G$.
\end{thm}

Soit $X$ l'ensemble des fonctions centrales sur $G$ \`{a} valeurs
dans $\ZM$. C'est un $\ZM$-module libre, engendr\'{e} par les $f_s$.
Il y a sur $X$ une structure d'anneau naturelle, la
\emph{convolution}, d\'{e}finie par $f,g\mapsto f * g$ avec
$$f*g\,(s)=\sum_{uv=s}{f(u)g(v)}.$$
On v\'{e}rifie que l'anneau $X$ est associatif et commutatif, et a
pour \'{e}l\'{e}ment unit\'{e} la \og fonction de Dirac\fg\ $f_1$.\\
Si $f\in X$, on lui associe l'endomorphisme
$\rho(f)=\sum{f(s)\rho(s)}$. D'apr\`{e}s le lemme \ref{homothetie},
$\rho(f)$ est une homoth\'{e}tie, i.e. appartient \`{a} $\CM$. On
obtient ainsi une application $\widetilde{\rho}: X\rightarrow \CM$,
qui est un homorphisme d'anneaux (en vertu de la formule
$\rho(f*f')=\rho(f).\rho(f')$). L'image $\widetilde{\rho}(X)$ de $X$
est donc un sous-anneau de $\CM$ qui est un $\ZM$-module de type
fini. On d\'{e}duit de l\`{a}, par un argument standard, que les
\'{e}l\'{e}ments de cet anneau sont des entiers alg\'{e}briques.
Comme $\frac{c(s)\chi(s)}{n}$ est \'{e}gal \`{a}
$\widetilde{\rho}(f_s)$, le th\'{e}or\`{e}me en r\'{e}sulte.~\findem

\begin{coro}
Le degr\'{e} de la repr\'{e}sentation $\rho$ divise l'ordre du
groupe.
\end{coro}
En effet, soit $n=\chi(1)$ et soit $f$ d\'{e}finie sur $G$ par
$f(s)=\chi(s^{-1})$. Alors $f$ v\'{e}rifie les hypoth\`{e}ses du th.
\ref{thmA}, donc $n^{-1}\sum_{s\in G}{\chi(s^{-1})\chi(s)}$ est un
entier alg\'{e}brique. Or $n^{-1}\sum_{s\in
G}{\chi(s^{-1})\chi(s)}=|G|/n$ donc est dans $\QM$. Les seuls
\'{e}l\'{e}ments de $\QM$ qui soient entiers sur $\ZM$ sont les
\'{e}l\'{e}ments de $\ZM$, donc $n$ divise $|G|$.~\findem

\begin{coro}\label{cor2}
Soit $\chi$ un caract\`{e}re irr\'{e}ductible de degr\'{e} $n$ d'un
groupe $G$ et soit $s\in G$. Si $c(s)$ et $n$ sont premiers entre
eux, alors $\frac{\chi(s)}{n}$ est un entier alg\'{e}brique. De
plus, si $\chi(s)\neq0$, alors $\rho(s)$ (associ\'{e} \`{a}
$\chi(s)$) est une homoth\'{e}tie.
\end{coro}

Si $c(s)$ et $n$ sont premiers entre eux, il existe, d'apr\`{e}s le
th\'{e}or\`{e}me de B\'{e}zout, deux \'{e}l\'{e}ments $a,b\in \ZM$
tels que $ac(s)+bn=1$. Donc:
$$\frac{\chi(s)}{n}=\frac{ac(s)\chi(s)}{n}+b\chi(s).$$
Or $\chi(s)$ est un entier alg\'{e}brique par la prop. \ref{propalg}
et $\frac{c(s)\chi(s)}{n}$ aussi par le th. \ref{thmB}. Cela
d\'{e}montre la premi\`{e}re assertion du corollaire. Pour la
seconde, nous utiliserons le:

\begin{lemme}
Soient $\lambda_{1},\dots,\lambda_{n}$ des racines de l'unit\'{e}.
Si $(\lambda_{1}+\cdots+\lambda_{n})/n$ est un entier
alg\'{e}brique, alors, soit $\lambda_{1}+\cdots+\lambda_{n}=0$, soit
tous les $\lambda_{i}$ sont \'{e}gaux.
\end{lemme}

\emph{D\'{e}monstration du lemme.} Si $\lambda$ est une racine de
l'unit\'{e}, son \'{e}quation int\'{e}grale minimale est de la forme
$x^{p}-1=0$, donc tout conjugu\'{e} de $\lambda$ est aussi racine de
l'unit\'{e} [rappel: si $z\in \CM$ est un nombre alg\'{e}brique, on
appelle conjugu\'{e} de $z$ toute racine de l'\'{e}quation minimale
de $z$]. Soit $z=(\lambda_{1}+\cdots+\lambda_{n})/n$; alors $z$ est
un entier alg\'{e}brique par hypoth\`{e}se et tout conjugu\'{e} $z'$
de $z$ s'\'{e}crit $(\lambda_{1}'+\cdots+\lambda_{n}')/n$ o\`{u} les
$\lambda_{i}'$ sont des racines de l'unit\'{e}. Donc $|z'|\leqslant
1$. Notons $Z$ le produit de tous les conjugu\'{e}s de $z$. Alors
$|Z|\leqslant 1$. D'autre part $Z$ est rationnel (c'est le terme
constant, au signe pr\`{e}s, du polyn\^{o}me minimal de $z$). Enfin
$Z$ est entier alg\'{e}brique (comme produit d'entiers
alg\'{e}briques) donc $Z\in \ZM$. Si $Z=0$, l'un des conjugu\'{e}s
de $z$ est nul donc $z=0$; si $|Z|=1$, tous les conjugu\'{e}s de $z$
ont $1$ comme module, donc $|z|=1$ et tous les $\lambda_{i}$ sont
\'{e}gaux, puisque leur somme est de module $n$ (cf. note
\ref{footnote} \`{a} la prop. \ref{A14}).~\findem

\bigskip Revenons \`{a} la seconde assertion du corollaire. Soient
$(\lambda_{i})$ les valeurs propres de $\rho(s)$; ce sont des
racines de l'unit\'{e} et $\chi(s)=\Tr
\big(\rho(s)\big)=\sum{\lambda_{i}}$. On conclut avec le
lemme.~\findem

\section{Application: th\'{e}or\`{e}me de Burnside}\label{Burn}

On reprend les notations du \S\ pr\'{e}c\'{e}dent.
\begin{thm}Soient $s$ un \'{e}l\'{e}ment de $G\mathbf{-}\{1\}$ et $p$ un nombre premier. Supposons que $c(s)$
($=$ nombre d'\'{e}l\'{e}ments de $G$ conjugu\'{e}s \`{a} $s$) soit
une puissance de $p$. Alors il existe un sous-groupe normal $N$ de
$G$, distinct de $G$, tel que l'image de $s$ dans $G/N$ appartienne
au centre de $G/N$.
\end{thm}

Soit $r_{G}$ le caract\`{e}re de la repr\'{e}sentation
r\'{e}guli\`{e}re (cf. \S\ \ref{fctcentr}). On a $r_{G}(s)=0$ car
$s\neq1$. D'autre part, $r_{G}=\sum{n_{\chi}\chi}$ (o\`{u} l'on
somme sur l'ensemble des caract\`{e}res irr\'{e}ductibles) avec
$n_{\chi}=\chi(1)$. Donc $\sum{\chi(1)\chi(s)}=0$ ou encore
$1+\sum_{\chi\neq1}{\chi(1)\chi(s)}=0$ et donc $\sum_{\chi\neq
1}{\frac{\chi(1)\chi(s)}{p}}=-\frac{1}{p}$; or $-\frac{1}{p}$ n'est
pas un entier alg\'{e}brique. Donc il existe un caract\`{e}re $\chi$
irr\'{e}ductible et distinct de $1$ tel que
$\frac{\chi(1)\chi(s)}{p}$ ne soit pas un entier alg\'{e}brique. En
particulier, $\chi(s)\neq0$ et $p$ ne divise pas $\chi(1)$, donc
$c(s)$ et $\chi(1)$ sont premiers entre eux. D'apr\`{e}s le cor.
\ref{cor2}, si $\rho$ est la repr\'{e}sentation de $G$ attach\'{e}e
\`{a} $\chi$, $\rho(s)$ est une homoth\'{e}tie. Soit alors $N=\ker
\rho$; le groupe $N$ est un sous-groupe de $G$, normal et distinct
de $G$. D'autre part, $G/N$ s'identifie \`{a} l'image de $\rho$ dans
$\mathbf{GL}_{n}(\CM)$, or $s$ s'envoie par $\rho$ sur une
homoth\'{e}tie, donc dans le centre de $G/N$.~\findem

\bigskip On va en d\'{e}duire le:
\begin{thm}[Burnside]\label{thmburn}
Tout groupe $G$ d'ordre $p^{\alpha}q^{\beta}$ (o\`{u} $p$ et $q$
sont des nombres premiers) est r\'{e}soluble.
\end{thm}

On suppose $\alpha$ et $\beta$ non nuls (sinon $G$ est nilpotent).
On va d\'{e}montrer le th\'{e}or\`{e}me par r\'{e}currence sur $|G|$.\\
Il existe $s\in G\mathbf{-}\{1\}$ tel que $q$ ne divise pas $c(s)$.
En effet, $G=\{1\}\cup \{Cl(s)\}$ o\`{u} les $Cl(s)$ sont les
classes de conjugaison autres que 1. Donc
$|G|=1+\sum{|Cl(s)|}=1+\sum{c(s)}$. Comme $q$ divise $|G|$, il
existe $s$ tel que $q$ ne divise pas $c(s)$. Or $c(s)$ divise $|G|$,
donc $c(s)$ est une puissance de $p$. D'apr\`{e}s ce qui
pr\'{e}c\`{e}de, il existe un sous-groupe normal $N$ de $G$,
distinct de $G$, tel que l'image de $s$ dans $G/N$ soit dans le
centre de $G/N$. Si $N\neq\{1\}$, l'hypoth\`{e}se de r\'{e}currence
appliqu\'{e}e \`{a} $N$ et \`{a} $G/N$ montre qu'ils sont
r\'{e}solubles; donc $G$ est r\'{e}soluble. Si $N=\{1\}$, le centre
$C$ de $G$ est non trivial puisqu'il contient $s$; on peut appliquer
l'hypoth\`{e}se de r\'{e}currence \`{a} $C$ et \`{a} $G/C$: $G$ est
r\'{e}soluble.~\findem

\section{D\'{e}monstration du th\'{e}or\`{e}me de Frobenius}\label{demofrobenius}
Soient $G$ un groupe fini et $H$ un sous-groupe de $G$ \og ne
rencontrant pas ses conjugu\'{e}s\fg, c'est-\`{a}-dire tel que:$$
H\cap gHg^{-1}=\{1\}$$ pour tout $g\in G\mathbf{-}H$, cf. \S\
\ref{frob}.

\enlargethispage{2.5\baselineskip}%
\begin{thm}[Frobenius]\label{frobannexe}
Soit $N=\{1\}\cup\{G\mathbf{-} \bigcup_{g\in G}gHg^{-1}\}$. Alors
$N$ est un sous-groupe normal de $G$; de plus $G$ est produit
semi-direct de $H$ par $N$.
\end{thm}

La d\'{e}monstration consiste \`{a} montrer que les
repr\'{e}sentations lin\'{e}aires de $H$ se prolongent \`{a} $G$
tout entier. Tout d'abord:

\begin{lemme}
Soit $f$ une fonction centrale sur $H$. Il existe une fonction
centrale $\widetilde{f}$ sur $G$, et une seule, qui satisfait aux
deux propri\'{e}t\'{e}s suivantes:
\begin{enumerate}
 \item[(1)] $\widetilde{f}$ prolonge $f$ {\rm (i.e. $\widetilde{f}(h)=f(h)$ si
$h\in H$)},
 \item[(2)] $\widetilde{f}(x)=f(1)$ si $x\in N$ {\rm (i.e. $\widetilde{f}$ est
constante sur $N$)}.
\end{enumerate}
\end{lemme}

C'est imm\'{e}diat: si $x\notin N$, on \'{e}crit $x$ sous la forme
$ghg^{-1}$, avec $g\in G$, $h\in H$, et l'on pose
$\widetilde{f}(x)=f(h)$; cela ne d\'{e}pend pas du choix de $g,h$
car si $g'h'g'^{-1}=ghg^{-1}$, on a $g'^{-1}ghg^{-1}g'=h'$, donc
$h=h'=1$ ou bien $g'^{-1}g\in H$; donc $h$ et $h'$ sont
conjugu\'{e}s dans $H$, et, comme $f$ est centrale, on a
$f(h)=f(h')$.~\findem

\medskip Dans l'\'{e}nonc\'{e} suivant, on utilise la notation
${\langle\alpha,\beta\rangle}_G$ pour d\'{e}signer le produit
scalaire $\frac{1}{|G|}\sum_{s\in G}{\alpha(s^{-1})\beta(s)}$
relatif au groupe $G$. De m\^{e}me ${\langle\alpha,\beta\rangle}_H$
d\'{e}signe le produit scalaire relatif \`{a} $H$. Si $F$ est une
fonction sur $G$, on note $F_H$ sa restriction \`{a} $H$.

\begin{lemme}
Soient $f$ et $\widetilde{f}$ comme dans le lemme pr\'{e}c\'{e}dent,
et soit $\theta$ une fonction centrale sur $G$. On a
\begin{equation}\label{**}
{\langle\widetilde{f},\theta\rangle}_G={\langle
f,\theta_H\rangle}_H+f(1){\langle 1,\theta\rangle}_G-f(1){\langle
1,\theta_H\rangle}_H.
\end{equation}
\end{lemme}

D\'{e}montrons cette \'{e}galit\'{e}. Elle est vraie pour $f=1$ car
alors $\widetilde{f}=1$. Par lin\'{e}arit\'{e}, il suffit donc de
v\'{e}rifier (\ref{**}) pour les fonctions $f$ telles que $f(1)=0$.
Soit alors $\mathcal{R}$ un syst\`{e}me de repr\'{e}sentants des
classes \`{a} gauche modulo $H$: les diff\'{e}rents conjugu\'{e}s de
$H$ sont les $rHr^{-1}$ pour $r$ d\'{e}crivant $\mathcal{R}$ et tout
conjugu\'{e} d'un \'{e}l\'{e}ment de $H$ (distinct de $1$)
s'\'{e}crit de fa\c{c}on unique $rhr^{-1}$ avec $r\in \mathcal{R}$
et $h\in H$. Cela dit,
$$ {\langle \widetilde{f},\theta\rangle}_{G}\ =\ \frac{1}{|G|}\sum_{s\in G}{\widetilde{f}(s^{-1})\theta(s)}$$
et $\widetilde{f}$ est nulle en dehors des conjugu\'{e}s
d'\'{e}l\'{e}ments de $H$; donc
\begin{eqnarray*}
{\langle
\widetilde{f},\theta\rangle}_{G}&=&\frac{1}{|G|}\sum_{(r,h)\,\in\, \mathcal{R}\times H}{\widetilde{f}(h^{-1})\theta(h)} \\
&=&{\langle f,\theta_{H}\rangle}_{H}
\end{eqnarray*}
puisque $|G|=|H|.|\mathcal{R}|$.~\findem

\bigskip\emph{Cas particulier:} si $\theta=1$, on a ${\langle
\widetilde{f},1\rangle}_{G}={\langle f,1\rangle}_{G}$. L'application
$f\mapsto \widetilde{f}$ est donc une isom\'{e}trie, autrement dit:
${\langle \widetilde{f}_{1},\widetilde{f}_{2}\rangle}_{G}={\langle
f_{1},f_{2}\rangle}_{H}$. En effet,
$$\frac{1}{|G|}\sum_{s\in G}{\widetilde{f}_{1}(s^{-1})f_{2}(s)}=\ {\langle f_{1}^{*}\widetilde{f}_{2}^{},1\rangle}_{G}$$
o\`{u} l'on a pos\'{e} $f_{1}^{*}(s)=f_{1}(s^{-1})$, et
$${\langle \widetilde{f}_{1},\widetilde{f}_{2}\rangle}_{G}\ =\ {\langle \widetilde{f_1^*}\widetilde{f}_{2}^{},1\rangle}_{G}\ =\ {\langle f_{1}^{*}{f}_{2}^{},1\rangle}_{H}\ =\ {\langle f_{1},f_{2}\rangle}_H$$
car $\widetilde{f_{1}f_{2}}=\widetilde{f_{1}}\widetilde{f_{2}}$. On
en d\'{e}duit:
\begin{prop}
Si $\chi$ est un caract\`{e}re de $H$ et $\theta$ un caract\`{e}re
de $G$, alors ${\langle \widetilde{\chi},\theta\rangle}_{G}$ est un
entier.
\end{prop}
Il suffit de montrer que dans (\ref{**}) tous les termes du membre
de droite sont dans $\ZM$, ce qui est clair.~\findem

\bigskip Soient $\theta_{1},\dots,\theta_{n}$ les diff\'{e}rents caract\`{e}res
irr\'{e}ductibles de $G$; on a
$\widetilde{\chi}=\sum{c_{i}\theta_{i}}$ avec $c_{i}\in \ZM$ pour
tout~$i$ d'apr\`{e}s ce qui pr\'{e}c\`{e}de.

\begin{prop}
Supposons $\chi$ irr\'{e}ductible. Alors tous les $c_{i}$ sont nuls,
sauf l'un d'eux qui est \'{e}gal \`{a} $1$.\\ {\rm(Autrement dit
$\widetilde{\chi}$ est un caract\`{e}re irr\'{e}ductible de $G$.)}
\end{prop}
En effet, ${\langle
\widetilde{\chi},\widetilde{\chi}\rangle}_{G}={\langle
\chi,\chi\rangle}_{H}=1=\sum{c_{i}^{2}}$ donc tous les $c_{i}$ sont
nuls sauf un, soit $c_{i_{0}}$, dont le carr\'{e} est $1$. Si
$c_{i_{0}}$ \'{e}tait \'{e}gal \`{a} $-1$, on aurait
$\widetilde{\chi}=-\theta_{i_0}$. Or $\widetilde{\chi}(1)=\chi(1)
>0$ et $\theta_{i_0}(1)>0$, donc c'est impossible. Donc
$c_{i_{0}}=1$, i.e. $\widetilde{\chi}=\theta_{i_0}$.~\findem

\begin{coro}\label{A27}
Si $\chi$ est un caract\`{e}re de $H$, alors $\widetilde{\chi}$ est
un caract\`{e}re de $G$.
\end{coro}
 Cela r\'{e}sulte de la proposition, en d\'{e}composant
$\chi$ comme somme de caract\`{e}res irr\'{e}ductibles.~\findem

\bigskip D\'{e}montrons maintenant le th. \ref{frobannexe}. Choisissons
une repr\'{e}sentation $\rho$ de $H$ de noyau trivial, par exemple
la repr\'{e}sentation r\'{e}guli\`{e}re. Soit $\chi$ le
caract\`{e}re de $\rho$ et soit $\widetilde{\rho}$ une
repr\'{e}sentation de $G$ de caract\`{e}re $\widetilde{\chi}$, cf.
cor. \ref{A27}. Si $s$ est conjugu\'{e} \`{a} un \'{e}l\'{e}ment de
$H-\{1\}$, on a $\widetilde{\rho}(s)\neq 1$. D'autre part, si $s\in
N$, on a $\widetilde{\chi}(s)=\chi(1)=\widetilde{\chi}(1)$, d'o\`{u}
$\widetilde{\rho}(s)=1$, cf. d\'{e}monstration de la prop.
\ref{A14}. On a donc $N=\ker{\widetilde{\rho}}$, ce qui montre que
$N$ est un sous-groupe normal de $G$.~\findem

~
\bigskip
{\bf Mode d'emploi}
\bigskip
\begin{itemize}
\item Groupes r\'{e}solubles et nilpotents, $p$-groupes et th\'{e}or\`{e}mes
de Sylow: \cite{BourbAI}, \cite{Hall}, \cite{Hupp}, \cite{Jac},
\cite{Lang1}, \cite{Zass}.

\item Applications \`{a} la th\'{e}orie de Galois: \cite{Jac}.

\item Extensions, syst\`{e}mes de facteurs et cohomologie:
\cite{Ca-Fr}, \cite{Hupp}, \cite{Lang2}, \cite{Se1}, \cite{Zass}.

\item Th\'{e}or\`{e}mes de Hall: \cite{Hall}, \cite{Hupp}.

\item Groupes de Frobenius: \cite{Feit}, \cite{Hupp} (pour la
structure des groupes de type $\mathcal F$, voir \cite{Wolf}).

\item Transfert: \cite{Hupp}, \cite{Lang2}, \cite{Se1},
\cite{Zass} (voir aussi I. Schur, {\it Gesam. Abh.}, vol. 1,
p.79).

\item Repr\'{e}sentations et caract\`{e}res: \cite{BourbAVIII},
\cite{Feit}, \cite{Hupp}, \cite{Lang1}, \cite{Se2} (voir aussi
F.G. Frobenius, {\it Gesam. Abh.}, vol.3, et I. Schur, {\it Gesam.
Abh.}, vol.1).
\end{itemize}

\newpage

\

\vspace{2.2cm}

{\bf \Huge{Index}}\label{index}\thispagestyle{plain}

\

\

action (de groupe --): \ref{action}\\
Alperin (th\'{e}or\`{e}me d'--): \ref{Alperin}\\
Burnside (th\'{e}or\`{e}mes de --): \ref{Burnside1},
\ref{Burnside2},
\ref{Burn}\\
caract\`{e}re (d'une repr\'{e}sentation): \ref{caractere}\\
caract\'{e}ristique (sous-groupe --): \ref{caracteristique}\\
Cauchy (th\'{e}or\`{e}me de --): \ref{2.2}\\
centrale (fonction --): \ref{fctcentr}\\
centralisateur (d'un \'{e}l\'{e}ment): \ref{centralisateur}\\
centre (d'un groupe): \ref{centre}\\
Chebotarev-Frobenius (th\'{e}or\`{e}me de --): \ref{Chebotarev-Frobenius}\\
classe (de conjugaison): \ref{classeconjugaison}\\
classe (de nilpotence d'un groupe): \ref{classenilpotence}\\
--- (de r\'{e}solubilit\'{e} d'un groupe): \ref{classeresolubilite}\\
cobord: \ref{cobord}\\
cocha\^{i}ne: \ref{cobord}\\
cocycle: \ref{cobord}\\
cohomologie: \ref{cobord}\\
commutateur: \ref{commutateur}\\
commutateurs (groupe des --): \ref{commutateur}\\
couple de Frobenius: \ref{couplefrob}\\
degr\'{e} (d'une repr\'{e}sentation lin\'{e}aire): \ref{representation}\\
d\'{e}riv\'{e} (groupe --): \ref{commutateur}\\
descendante (suite centrale --): \ref{suitecentrale}\\
drapeau complet: \ref{commutateur}\\
\'{e}l\'{e}mentaire ($p$-groupe ab\'{e}lien --): \ref{elementaire}\\
entier alg\'{e}brique: \ref{propintegralite}\\
extension: \ref{extensions}\\
Feit-Thompson (th\'{e}or\`{e}me de --): \ref{commutateur}\\
filtration: \ref{filtrations}\\
fixateur: \ref{stabilisateur}\\
Frattini (sous-groupe de --): \ref{sgFrattini}\\
--- (th\'{e}or\`{e}me de --): \ref{thmFrattini}\\
Frobenius (th\'{e}or\`{e}me de --): \ref{thmfrobenius},
\ref{demofrobenius}\\
fusion: \ref{2.4}\\
Gauss (lemme de --): \ref{Gauss}\\
gradu\'{e}: \ref{filtrations}\\
Hall (sous-groupe de --): \ref{Hall}\\
--- (th\'{e}or\`{e}me de --): \ref{Hall}\\
homomorphisme crois\'{e}: \ref{cobord}\\
image r\'{e}ciproque (d'une extension): \ref{4.5}\\
invariant (sous-groupe --): \ref{normal}\\
Jordan (th\'{e}or\`{e}me de --): \ref{6.1}\\
Jordan-H\"{o}lder (suite de --, th\'{e}or\`{e}me de --): \ref{filtrations}\\
Kolchin (th\'{e}or\`{e}me de --): \ref{Kolchin}\\
\lhead{Index}
localement conjugu\'{e}s (\'{e}l\'{e}ments --): \ref{2.4}\\
maximal (sous-groupe --): \ref{sgFrattini}\\
Miller-Wielandt (d\'{e}monstration de --): \ref{Miller-Wielandt}\\
minimal (groupe simple --): \ref{minimal}\\
nilpotent (groupe --): \ref{grpe nilp}\\
normal (sous-groupe --): \ref{normal}\\
normalisateur (d'un sous-groupe): \ref{normalisateur}\\
orbite: \ref{action}\\
orthogonalit\'{e} (de caract\`{e}res): \ref{A.2}\\
$p$-compl\'{e}ment: \ref{5.6}\\
$p$-groupe: \ref{sylow}\\
$p$-Sylow: \ref{sylow}\\
permutable (syst\`{e}me -- de sous-groupes de Sylow): \ref{5.3}\\
permutables (sous-groupes --): \ref{5.2}\\
$\Pi$-sous-groupe: \ref{Pisg}\\
$\Pi$-Sylow: \ref{Pisg}\\
primitive (action --): \ref{sgFrattini}\\
propri\'{e}t\'{e} $\mathcal F$: \ref{6.4}\\
r\'{e}guli\`{e}re (repr\'{e}sentation --): \ref{repr-reg}\\
rel\`{e}vement (d'un homomorphisme): \ref{4.5}\\
repr\'{e}sentation (lin\'{e}aire): \ref{representation}\\
r\'{e}soluble (groupe --): \ref{grpe resol}\\
sans point fixe (action --): \ref{6.4}\\
Schur (lemme de --): \ref{A.2}\\
--- (th\'{e}or\`{e}me de --): \ref{7.1}\\
section: \ref{section}\\
semi-direct (produit --): \ref{semidirect1}, \ref{4.4}\\
simple (groupe --): \ref{grpe simple}\\
stabilisateur: \ref{stabilisateur}\\
Stallings-Swan (th\'{e}or\`{e}mes de --): \ref{Stallings-Swan}\\
suite centrale descendante: \ref{suitecentrale}\\
Sylow (th\'{e}or\`{e}mes de --): \ref{2.2}, \ref{2.3}\\
Thompson (th\'{e}or\`{e}mes de --): \ref{thmThompson}, \ref{thmThompson2}\\
transfert: \ref{7.1}\\
transitive (action --): \ref{action}\\
triviale (extension --): \ref{extensions}\\
Wedderburn (th\'{e}or\`{e}me de --): \ref{Wedderburn}\\
Wielandt (th\'{e}or\`{e}me de --): \ref{5.5}\\
Zassenhaus (th\'{e}or\`{e}me de --): \ref{4.4}\\

\hfill \footnotesize Version du $\mit 9$ juin $\mit 2008$

\hfill Ce cours est disponible en ligne sur le serveur ArXiv sous la
r\'{e}f\'{e}rence: math.GR/0503154
\end{document}